\documentclass[a4paper,11pt]{article}
\usepackage{amsmath,latexsym,amssymb,amsfonts,natbib}
\usepackage[a4paper,hmargin=1.87 cm]{geometry}
\usepackage[dvips]{graphicx}
\usepackage{pifont}
\usepackage{pstricks,pst-fill,pst-grad}
\usepackage{rotating}
\usepackage{subfigure}
\usepackage{graphics}
\usepackage{hyperref}
\usepackage{multirow}
\usepackage{bigints}
\usepackage{multicol}
\usepackage{courier}
\usepackage[utf8]{inputenc}
\usepackage{indentfirst}
\usepackage{ccaption}
\usepackage{verbatim}
\usepackage{makeidx}
\usepackage{listings}
\usepackage{dsfont}
\usepackage{inputenc}
\DeclareMathOperator*{\argmax}{arg\,max}

\newcommand{\proof}     {\paragraph{Proof}}
\newcommand{\carre}     {\hfill$\Box$}
\numberwithin{equation}{section}
\newtheorem{defi}{Definition}[section]
\newtheorem{lem}{Lemma}[section]
\newtheorem{theo}{Theorem}[section]
\newtheorem{coro}{Corollary}[section]

\newtheorem{rem}{Remark}

\newtheorem{cond}{Condition}
\newtheorem{ass}{Assumption}

\newcommand{\1}{\lfloor j t_1 \rfloor} 

\title{Propagation of chaos and large deviations in mean-field models with jumps on block-structured networks}
\author{ Donald A. Dawson\thanks{School of Mathematics and Statistics, Carleton University, 1125 Colonel By Drive
Ottawa, Ontario K1S 5B6, Canada. ddawson@math.carleton.ca}, Ahmed Sid-Ali\thanks{School of Mathematics and Statistics,
Carleton University, 1125 Colonel By Drive Ottawa, Ontario K1S 5B6, Canada. ahmedsidali@cunet.carleton.ca}, Yiqiang Q. Zhao\thanks{School of Mathematics and Statistics, Carleton University, 1125 Colonel By Drive Ottawa, Ontario K1S 5B6, Canada. zhao@math.carleton.ca} \\
{\small }\\
{\small }\\
}
\date{}
\begin{document}
\maketitle

\begin{abstract}
A system of interacting  multiclass finite-state jump processes is analyzed. The model under consideration consists of a block-structured network with dynamically changing multi-colors nodes. The interaction is local and described through local empirical measures. Two levels of heterogeneity are considered: between and within the blocks where the nodes are labeled into two types. The central nodes are those connected only to the nodes of the same block whereas the peripheral nodes are connected to both the nodes of the same block and to some nodes from other blocks. The limits of such systems as the number of particles tends to infinity are investigated. Under regularity conditions on the peripheral nodes, propagation of chaos and law of large numbers are established in a multi-population setting. In particular, it is shown that, as the number of nodes goes to infinity, the behavior of the different classes of nodes can be represented by the solution of a McKean-Vlasov system. Moreover, we prove large deviation principles for the vectors of empirical measures and the empirical processes. Our large deviations results extend classical results of \cite{Daw+Gart87} and \cite{Leonard95}.
\end{abstract}
{\it 2020 Mathematics Subject Classification:} Primary 60K35 60J27 60F10 05C82 60J74.
\\
\\
{\it Keywords:} Inhomogeneous graphs; Interacting particle systems; Mean-field limit; McKean-Vlasov equation; Propagation of chaos; Large deviations; Multi-class populations; jump processes.

\section{Introduction}
Since McKean's seminal paper \cite{McKean66}, the mean-field theory has been widely used for the study of  large stochastic interacting particle systems arising from various domains such as  statistical physics  \cite{McKean66,McKean66(2),Daw83,Gart88}, biological systems  \cite{Daw2017,Mel+Ben2015}, communication networks \cite{Gra+Mel93,Gra+Mel95,Ben+Leboud2008,Graha2000}, mathematical finance \cite{Kley+Klu2015,Giesecke+al2015}, etc. This theory, first initiated in connection with a  mathematical foundation of the Boltzmann equation, aims for a mathematically rigorous treatment of the time evolution of stochastic systems with long-range weak interaction where the interaction between particles is realized via the \textit{empirical measure} of the particle configuration. For such systems, it is then natural to investigate the behavior of the \textit{empirical process} instead of considering the particle configuration itself. In particular, one is interested in the investigation of limit theorems such as laws of large numbers and large deviations for the empirical process in the limit as the particle number tends to infinity. Another concept that plays an important role is  \textit{the propagation of chaos} introduced by  Mark Kac in the context of kinetic theory \cite{Kac56} and largely commented on in the literature. See  \cite{Sznitman91} and \cite{Gart88} for detailed developments on the subject.

In the classical case, the studied systems are homogeneous with complete interaction graphs, that is, the particles are exchangeable and each particle interacts with every other particle. In such a setting, the big picture is well understood and various asymptotic results have been established for a variety of models. One can consult \cite{Daw93}, \cite{Gart88} or \cite{Sznitman91} for an overview. However, though such assumptions are reasonable in statistical physics and well describe a variety of phenomena, it may no longer be the case when considering other applications. In this regard, many researchers studied a bunch of new applications of interacting particle systems where the homogeneity or the complete interaction assumption is not relevant.

One direction towards heterogeneity is the study of systems where it  is inherent to the particles due to their different backgrounds. In such cases, one cannot presume the particles to be identically nor symmetrically distributed. Instead, one relies on additional conditions to establish limiting results. For instance, in \cite{Finnoff93,Finnoff94,Giesecke+al2015}, models for the activities of heterogeneous economical agents were proposed and laws of large numbers  were proved under some regularity conditions. Within the same spirit, limiting results were established in \cite{Naga+Tana87} for a system of reflected diffusions segregated into two groups of blue and red particles and subject to a reflection condition. These results were extended in \cite{Naga+Tana87(2)} to the case of drift coefficients not of average form. Among other examples, we can cite the more recent works in \cite{Buck+Li2014,Car+Zhu2016}  where mean-field game models are considered with one single major player and statistically identical minor players, and the propagation of chaos was proved for the minor players conditioned on the major player.

Another natural extension in the direction of heterogeneous systems, with which our current work is in line, is the multi-population paradigm where the particles are divided into different groups within which they are homogeneous or partially homogeneous. The motivation is that many systems in statistical physics, chemistry, communication networks, biology, finance, etc., involve varied classes of similar particles (see references below). We propose to take a step forward in this direction by studying some asymptotic results for large interacting particle systems with jumps on block-structured networks.

In this paper, we will set up a model for block-structured networks with dynamically changing multi-color nodes. The evolution of node colors is described by a sequence of finite-state pure-jump processes interacting through local empirical measures describing the neighborhood of each node. The nodes of the network are divided into finite number of blocks. In addition, the nodes within each block are divided into two subgroups: central and peripheral nodes. The central nodes are those connected only with the nodes of the same block whereas the peripheral nodes interact with both the particles of the same block and with some particles from other blocks. Thus, our model describes  two levels of heterogeneity: between blocks  and within blocks.

Our idea is in a continuation of several existing works. In \cite{Collet2014,Collet+al2016},  a bi-populated Curie-Weiss model was studied, where the authors established, via a large deviation approach, the propagation of chaos and the asymptotic dynamics of the pair of group magnetizations in the infinite volume limit. The laws of large numbers and  a central limit theorem were proved in \cite{Kirsh+To2020} for an extension of this model to the case of heterogeneous coupling within and between groups. A related work \cite{Lowe+Schub2018}  studied the high temperature  fluctuations for block spin Ising models and  a central limit theorem was established. A variant of this model was analyzed in \cite{Kno+Lo+schu2020} where the vertices are divided into a finite number of blocks and pair interactions are given according to their blocks. The authors proved large deviation principles and central limit theorems. Other recent related works are the two-community noisy Kuramoto model studied in \cite{Mey2020} and the opinion dynamics with the Lotka-Volterra type interactions model studied in \cite{Ale+Min2019}. However, the closest to the current work is the model recently proposed in \cite{Bay+Wu2019} where the authors studied systems of weakly interacting jump processes on time-varying random graphs with dynamically changing multi-color edges. The node dynamics depends on the joint empirical distribution of all the nodes and edges connected to it, while the edge dynamics depends only on the corresponding nodes it connects. The law of large numbers, propagation of chaos, and central limit theorems were established by the authors for this model. Our model differs from the one in \cite{Bay+Wu2019} in several aspects. In the current work, the interaction between the nodes are local whereas it is global in \cite{Bay+Wu2019}. Moreover, the interacting system we consider is on static block-structured graphs whereas the one in \cite{Bay+Wu2019} is considered on time-varying random graphs with dynamically changing edges colors. Finally, the analysis  carried out and results  obtained in the current work are established in a multi-population context, which allows overcoming the heterogeneity due to the structure of the graph. Notice further that the current paper falls into the topic of interacting particle systems on large (random) networks which, in recent years, has attracted increasing attention. See e.g.,  \cite{Bham+Budh+Wu2019,Bay+Wu2019,Bay+Chak+Wu2020} and the references therein.

Along with the papers listed above, the multi-population framework was also considered for systems of interacting diffusions. We can mention for instance \cite{Kley+Klu2015} where a system of interacting Ornstein-Uhlenbeck processes on a heterogeneous network of credit-interlinked agents was analyzed, \cite{Toub2018,Boss+Faug+Tal2015,Budh+Wu2016} and the references therein for the study of neuronal networks composed of separate populations, or \cite{Chong+Klu2019,Nguyen+al2020} and the references therein for mean-field multi-class interacting diffusions models in a general setting.

The goal of the current work is the development of limiting results for interacting finite-state pure-jump processes on a class of block-structured networks. Our first main result, Theorem~\ref{theo2} and its Corollary \ref{coro-1}, gives propagation of chaos and law of large numbers under some regularity conditions on the degrees of the peripheral nodes. We show that in the mean-field limit, the asymptotic behavior of the node colors can be represented by the solution of a McKean-Vlasov system. Due to the lack of symmetry, we make use of the extension to multi-class systems of the notion of chaoticity and of Sznitman coupling methods developed in \cite{Graha2008,Graha+Rob2009}. The existence and uniqueness results for the limiting system are established in Theorem \ref{theo1}. The regularity conditions, which we impose on (see Condition $\ref{Cond-prin}$), can be compared to the uniform degree property introduced in \cite{Delattre+al2016} for a model of interacting diffusions on random graphs and the one introduced in \cite{Budh+Mukh+Wu2019} for a model of interacting pure-jump processes on sparse graphs. The difference lies in that we only impose regularity on the degrees of peripheral nodes and allow the blocks of the network to have very different number of nodes.

Another aspect which we are interested in is the large deviations property of the system. For this purpose and for simplicity, we will restrict ourselves to the case of a complete peripheral sub-graph, that is, the case where all peripheral nodes of the system are connected with each other. We then state our next main results in Theorem $\ref{large-dev-meas}$ which establishes the large deviation property for the \textit{empirical measure vector} over finite time duration followed by Theorem $\ref{large-dev-emp-proc}$ which gives the large deviation property for the \textit{empirical process vector}. These results generalize \cite{Leonard95} and \cite{Bork+Sund2012} to the multi-population context. Also, different from \cite{Leonard95} and similar to \cite{Bork+Sund2012}, we do not impose chaotic initial conditions but only impose converging initial conditions. The proofs of the large deviation results follow from the classical approach developed in \cite{Daw+Gart87} and adapted to the context of jump processes in \cite{Leonard95},  which provide tools for handling the technicalities arising from the multi-population context.

In summary, the current work is a contribution to the multi-population paradigm and a move towards heterogeneity for mean-field models and their large deviation behavior. The rest of this paper is organized as follows. The detailed model for interacting finite-state pure-jump processes on block-structured graphs is introduced in Section \ref{model}. In Section \ref{LLN-POC-sec} we introduce the regularity conditions on the sub-peripheral graph under which the propagation of chaos and law of large numbers results hold (see Condition \ref{Cond-prin}). These results are given in Theorem \ref{theo2} and Corollary \ref{coro-1}. Next in Section \ref{large-dev-sec} we present the large deviations principles for the empirical measure vector in Theorem $\ref{large-dev-meas}$ and for the empirical process vector in Theorem $\ref{large-dev-emp-proc}$.

\section{The model}
\label{model}

We describe in this section the model under investigation and the related notations. 
\subsection{The setting}
\subsubsection*{A block-structured network}
\begin{itemize}
\item Consider a block-structured graph $\mathcal{G}=(\mathcal{V},\Xi)$, where $\mathcal{V}$ is the set of nodes and $\Xi$ the set of edges, composed of $r$ blocks (communities) $C_1,\ldots,C_r$ of sizes $N_1,\ldots,N_r$, respectively. Denote by $|\mathcal{V}|=N_1+\cdots+N_r=N$ the total number of nodes in the network.
\item Each block $C_j$ is a clique, i.e. all the $N_j$ nodes are connected to each other.
\item The nodes of each block $C_j$ are divided into two categories:
      \begin{itemize}
       \item \textbf{Central nodes $C^c_j$:} connected to all the other nodes of the same block but not to any node from the other blocks. We set $|C_j^c|=N_j^c$.
       \item \textbf{Peripheral nodes $C^p_j$:} connected to all the other nodes of the same block and to some nodes from the other blocks. We set $|C^p_j|=N^p_j$.
      \end{itemize}
\end{itemize}

\subsubsection*{Multi-color nodes}
 Let $\mathcal{Z}=\{1,2,\ldots,K\}\subset\mathbb{N}$ be a set of $K$ colors. Suppose that each node of the graph $\mathcal{G}=(\mathcal{V},\Xi)$ is colored by one of the $K$ colors at each time. Define by $(\mathcal{Z},\mathcal{E})$ the directed graph where $\mathcal{E}\subset\mathcal{Z}\times\mathcal{Z}\backslash \{(z,z)| z \in\mathcal{Z}\}$ describes the set of admissible jumps for each particle. Moreover, whenever $(z,z')\in\mathcal{E}$, a particle colored by $z$ is allowed to move from $z$ to $z'$ at a rate that depends on the current state of the node and on the state of its neighbors (adjacent nodes).

For each $1\leq j\leq r$ and $n\in C_j^c$ (resp. $n\in C_j^p$), define by $(X_{n,j}(t),t\geq 0)$ the stochastic process that describes the color of the central (resp. peripheral) node $n$ at time $t$. Given the structure of the graph, rather than a global empirical measure, we introduce, for each $1\leq j\leq r$, the local empirical measures $\mu^N_j(t)$ describing the state of the block $j$ at time $t$ as
\begin{align*}
\mu_j^N(t)=\frac{1}{N_j}\sum_{n\in C_j}\delta_{X_{n,j}(t)}\in\mathcal{M}_1(\mathcal{Z}),
\end{align*}
 where $\mathcal{M}_1(\mathcal{Z})$ is the set of all probability measures over $\mathcal{Z}$, endowed with the topology of weak convergence. Using the heterogeneity of the nodes within each block, the empirical measure $\mu_j^N(t)$ can be decomposed as follows:
\begin{equation}
\begin{split}
\mu_j^N(t)&=\frac{1}{N_j}\sum_{n\in C_j}\delta_{X_{n,j}(t)}\\
          &=\frac{1}{N_j}\bigg(\sum_{n\in C^c_j}\delta_{X_{n,j}(t)}+\sum_{n\in C^p_j}\delta_{X_{n,j}(t)}\bigg)\\
          &=\frac{N_j^c}{N_j}\frac{1}{N_j^c}\sum_{n\in C^c_j}\delta_{X_{n,j}(t)}+\frac{N_j^p}{N_j}\frac{1}{N_j^p}\sum_{n\in C^p_j}\delta_{X_{n,j}(t)}\\
          &=\frac{N_j^c}{N_j}\mu_j^{c,N}(t)+\frac{N_j^p}{N_j}\mu_j^{p,N}(t),
\label{block-decomp}
\end{split}
\end{equation}
 where  $\mu_j^{c,N}(t)=\frac{1}{N_j^c}\sum_{n\in C^c_j}\delta_{X_{n,j}(t)}$ (resp. $\mu_j^{p,N}(t)=\frac{1}{N_j^p}\sum_{n\in C^p_j}\delta_{X_{n,j}(t)}$) is the empirical measure describing the state of the central (resp. peripheral) nodes of the $j$-th block at time $t$. The fractions $\frac{N_j^c}{N_j}$ (resp. $\frac{N_j^p}{N_j}$) represents the proportion of central (resp. peripheral) nodes in the block $j$.

Note that, given the symmetry between the central nodes of the same block, the neighborhood of any central node $n\in C_j^c$ at time $t$ is fully described by the empirical measure $\mu_j^N(t)$. However, this is not the case for the peripheral nodes since each peripheral node is connected to all nodes of the same block and  also with some peripheral nodes from the other blocks. Therefore, in order to describe the neighborhoods of the peripheral nodes, we introduce a set of local empirical measures which, for each peripheral node $n\in C_j^p$, describe the state of the star-shaped sub-graph centered at $n$ and composed of the nodes connected to $n$. Thus denoting by $deg(n)$ the degree of the peripheral node $n\in C_j^p$, we define the following empirical measure

\begin{equation}
\begin{split}
\mu_{n,j}^N(t)= & \frac{1}{deg(n)+1}\bigg(\sum_{\substack{m\in C^p_1\\(m,n)\in\Xi}}\delta_{X_{m,1}(t)}+\cdots+\sum_{m\in C^c_j}\delta_{X_{m,j}(t)}+\sum_{m\in C^p_j}\delta_{X_{m,j}(t)}+\cdots+\sum_{\substack{m\in C^p_r\\(m,n)\in\Xi}}\delta_{X_{m,r}(t)}\bigg)
 \\
= & \frac{M_1^n}{deg(n)+1}\frac{1}{M_1^n}\sum_{\substack{m\in C^p_1\\(m,n)\in\Xi}}\delta_{X_{m,1}(t)}+\cdots+\frac{N_j^c}{deg(n)+1}\frac{1}{N_j^c}\sum_{m\in C^c_j}\delta_{X_{m,j}(t)}+ \\
    &  \frac{N_j^p}{deg(n)+1}\frac{1}{N_j^p}\sum_{m\in C^p_j}\delta_{X_{m,j}(t)}+\cdots 
 +\frac{M_r^n}{deg(n)+1}\frac{1}{M_r^n}\sum_{\substack{m\in C^p_r\\(m,n)\in\Xi}}\delta_{X_{m,r}(t)},
\end{split}
\label{periph-decomp}
\end{equation}
where $M_i^n$ represents the number of peripheral nodes of the $i$-th block, which the peripheral node $n$ is connected with. In particular, since each block is a clique, $M_j^n+1=N_j^p$ for all $1\leq j\leq r$ and $n\in C_j^p$.

\subsubsection*{The random dynamics} 

The processes $(X_{n,j}(t),t\geq 0)$ are continuous-time Markov chains with state space $\mathcal{Z}$. The transition rate of each node depends on its current state and on the state of its neighbors described by the corresponding local empirical measure. Thus, the processes interact only through the dependence of their transition rates on the current empirical measures.
 We describe the random dynamic in each block $1\leq j\leq r$ as follows.
\begin{itemize}
\item \textbf{The central nodes dynamic.} For each central node $n\in C^c_j$, its color $X_{n,j}(t)$ goes from $z$ to $z'$, with $(z,z')\in (\mathcal{Z},\mathcal{E})$, at rate
\begin{align}
 \lambda_{z,z'}^c\bigg(\mu^{c,N}_j(t),\mu^{p,N}_j(t),\frac{N_j^c}{N_j},\frac{N_j^p}{N_j}\bigg),
\label{lamb-c}
\end{align}
 which depends on its current state and on the states of its neighbors through the empirical measures $\mu^{c,N}_j(t)$ and $\mu^{p,N}_j(t)$. The proportions $\frac{N_j^c}{N_j}$ and $\frac{N_j^p}{N_j}$ quantify the influence of the central and peripheral nodes on the transition rates.

\item \textbf{The peripheral nodes dynamic.} For each peripheral node $n\in C^p_j$, its color $X_{n,j}(t)$ transits from $z$ to $z'$, with $(z,z')\in (\mathcal{Z},\mathcal{E})$ at rate

\begin{equation}
\begin{split}
\lambda^p_{z,z'}\bigg(\mu_j^{c,N}(t),&\frac{1}{M_1^n}\sum_{\substack{m\in C^p_1\\(m,n)\in\Xi}}\delta_{X_{m,1}(t)},\ldots,\mu_j^{p,N}(t),\ldots,\frac{1}{M_r^n}\sum_{\substack{m\in C^p_r\\(m,n)\in\Xi}}\delta_{X_{m,r}(t)},\\
&\qquad\frac{N_j^c}{deg(n)+1},\frac{M_1^n}{deg(n)+1},\ldots,\frac{N_j^p}{deg(n)+1},\ldots,\frac{M_r^n}{deg(n)+1}\bigg),
\label{lamb-p}
\end{split}
\end{equation}
which also depends on its state and on the states of its neighbors through the empirical measures
\begin{align*}
\mu_j^{c,N}(t),\frac{1}{M_1^n}\sum_{\substack{m\in C^p_1\\(m,n)\in\Xi}}\delta_{X_{m,1}(t)},\ldots,\mu_j^{p,N}(t),\ldots,\frac{1}{M_r^n}\sum_{\substack{m\in C^p_1\\(m,n)\in\Xi}}\delta_{X_{m,r}(t)}.
\end{align*}

Again the proportions $\frac{N_j^c}{deg(n)+1},\frac{M_1^n}{deg(n)+1},\ldots,\frac{N_j^p}{deg(n)+1},\ldots,\frac{M_r^n}{deg(n)+1}$ quantify the influence of each group of nodes on the transition rate. We will introduce in Condition $\ref{Cond-prin}$ explicit forms for the rate functions. Not to clutter our notation and to facilitate the reading, let us introduce the following vectors

\begin{align}
\upsilon_j^{N}(t)=\bigg(\mu^{c,N}_j(t),\mu^{p,N}_j(t),\frac{N_j^c}{N_j},\frac{N_j^p}{N_j}\bigg),
\end{align}
\begin{align}
\upsilon_{n,j}^{N}(t)= \bigg(\mu_j^c(t),&\frac{1}{M_1^n}\sum_{\substack{m\in C^p_1\\(m,n)\in\Xi}}\delta_{X_{m,1}(t)},\ldots,\mu_j^p(t),\ldots,\frac{1}{M_r^n}\sum_{\substack{m\in C^p_r\\(m,n)\in\Xi}}\delta_{X_{m,r}(t)},\nonumber\\
&\frac{N_j^c}{deg(n)+1},\frac{M_1^n}{deg(n)+1},\ldots,\frac{N_j^p}{deg(n)+1},\ldots,\frac{M_r^n}{deg(n)+1}\bigg).
\end{align}
Thus, we will write  $ \lambda_{z,z'}^c\left(\upsilon_j^{N}(t)\right)$ instead of $(\ref{lamb-c})$ and $\lambda^p_{z,z'}\left(\upsilon_{n,j}^N(t)\right)$ instead of $(\ref{lamb-p})$.
\end{itemize}

\subsection{The infinitesimal generator}
For any $T\in (0,+\infty)$, denote by  $X^c_{n,j}: [0, T ] \rightarrow\mathcal{Z}$ for $n\in C_j^c$  and $X^p_{m,j}: [0, T ] \rightarrow\mathcal{Z}$  for $m\in C_j^p$ the processes that describe the evolution of the central and the  peripheral particles $n$  and $m$, respectively, over the time interval $[0,T]$. These are c\'adl\'ag paths and thus are elements of the Skorokhod space $\mathcal{D}([0,T],\mathcal{Z})$ equipped with the Skorokhod topology. Let
\begin{align*}
X^N=\big(X^c_{n,j},X^p_{m,j},n\in C_j^c,m\in C_j^p, 1\leq j\leq r\big)\in\mathcal{D}([0,T],\mathcal{Z}^N)
\end{align*}
denote the full description of paths of all $N$ particles. Thus the process $X^N$ is a Markov process with c\'adl\'ag paths,  the state space $\mathcal{Z}^N$, and the infinitesimal generator $\mathcal{L}^N$ acting on the bounded measurable functions $\phi$ on $\mathcal{Z}^N$ according to 

\begin{align*}
\mathcal{L}^N\phi(x^N) =& \sum_{j=1}^r\bigg[ \sum_{n\in C_j^c}\sum_{z':(z,z')\in\mathcal{E}}\lambda^c_{z,z'}\bigg(\frac{1}{N_j^c}\sum_{n\in C^c_j}\delta_{x_{n,j}},\frac{1}{N_j^p}\sum_{n\in C^p_j}\delta_{x_{n,j}},\frac{N_j^c}{N_j},\frac{N_j^p}{N_j}\bigg)\left( \phi(x^N_{n,z,z'})-\phi(x^N) \right)+
\\ 
    & \sum_{n\in C_j^p}\sum_{z':(z,z')\in\mathcal{E}}\lambda^p_{z,z'}\bigg(\frac{1}{N_j^c}\sum_{m\in C^c_j}\delta_{x_{m,j}},\frac{1}{M_1^n}\sum_{\substack{m\in C^p_1\\(m,n)\in\Xi}}\delta_{x_{m,1}},\ldots,\frac{1}{N_j^p}\sum_{m\in C^p_j}\delta_{x_{m,j}},\ldots,\frac{1}{M_r^n} \sum_{\substack{m\in C^p_r\\(m,n)\in\Xi}}\delta_{x_{m,r}}, 
\\
    &\qquad\qquad\qquad\qquad\qquad \frac{N_j^c}{deg(n)+1},\frac{M_1^n}{deg(n)+1},\ldots,\frac{N_j^p}{deg(n)+1},\ldots,\frac{M_r^n}{deg(n)+1}\bigg)\\    
&\qquad\qquad\qquad\qquad\qquad\qquad\qquad\qquad\qquad\qquad\qquad\qquad\qquad\qquad\qquad\times\left( \phi(x^N_{n,z,z'})-\phi(x^N) \right)\bigg],
\end{align*}
where $x^N=\big(x_{n,j},x_{m,j},n\in C_j^c,m\in C_j^p, 1\leq j\leq r\big)\in\mathcal{Z}^N$ and $x^N_{n,z,z'}$ describes the new configuration of the particles when we change the state of the $n$-th node from $z$ to $z'$.

\subsection{Stochastic differential equation representation}
 
Recall that, for each central node $n\in C^c_j$ (resp. peripheral node $n\in C^p_j$) at a given block $1\leq j\leq r$, the process $(X_{n,j}^{c}(t),t\geq 0)$ (resp. $(X_{n,j}^{p}(t),t\geq 0)$) is a continuous-time finite-state Markov chain with the time-dependent transition rate matrix $\left(\lambda^c_{z,z'}(\upsilon_{j}^N(t))\right)_{(z,z')\in\mathcal{E}}$ (resp. $\left(\lambda^p_{z,z'}(\upsilon_{n,j}^N(t))\right)_{(z,z')\in\mathcal{E}}$) and  the state space $\mathcal{Z}$. Using a classical approach (see e.g. \cite[p. $104$]{Skoro2009}), the Markov chains $X_{n,j}^{c}$ and $X_{n,j}^{p}$ can be represented, at least weakly, by the following system of stochastic differential equations
\begin{equation}
\begin{split}
X^{c}_{n,j}(t)&=X^{c}_{n,j}(0)+\int_{[0,t]\times\mathbb{R}_+}\sum_{(z,z')\in\mathcal{E}}\mathds{1}_{X^{c}_{n,j}(s-)=z}(z'-z)\mathds{1}_{\left[0,\lambda^c_{z,z'}(\upsilon_{j}^N(s-))\right]}(y)\mathcal{N}_{n,j}^c(ds,dy),\\
X^{p}_{n,j}(t)&=X^{p}_{n,j}(0)+\int_{[0,t]\times\mathbb{R}_+}\sum_{(z,z')\in\mathcal{E}}\mathds{1}_{X^{p}_{n,j}(s-)=z}(z'-z)\mathds{1}_{\left[0,\lambda^p_{z,z'}(\upsilon_{n,j}^N(s-))\right]}(y)\mathcal{N}_{n,j}^p(ds,dy),
\label{SDE-rep}
\end{split}
\end{equation}
where $\{\mathcal{N}_{n,j}^c, n\in C_j^c, 1\leq j\leq r\}$ and $\{\mathcal{N}_{n,j}^p, n\in C_j^p, 1\leq j\leq r\}$ are collections of Poisson random measures on $\mathbb{R}^2$ whose intensity measures are Lebesgue measures on $\mathbb{R}^2_+$. We will use the representation $(\ref{SDE-rep})$ in the analysis of the asymptotic behavior of the system when the total number of nodes $N$ goes to infinity.

\subsection{Examples}
As mentioned in the introduction of the present paper, the mean-field block models have been proposed to investigate various phenomena arising in different fields such as physics, engineering, biology, etc... This section aims to expose some examples of applications together with the references, of the model presented in the current paper and the established results. Our goal is to illustrate the usefulness of the proposed model and its flexibility to capture various phenomena. Of course, it remains a toy model that should probably be appropriately adapted to real contexts, but we believe that the insights from the current study are of great interest for both theoretical and practical purposes.

\subsubsection{Load balancing networks}

 Load balancing protocols are often used in queuing networks to improve the system performance by shortening the queue length, reducing the waiting time, and increasing the system throughput. In this regard, the mean-field approach has been proven to be useful, see e.g. \cite{Vved+Dobrushin+Kar96,Mitz96,Vved+Suhov97}. In particular, interesting work in this direction was proposed in \cite{Daw+Tang+Zhao2005}, where the authors considered a queuing network with $N$ nodes in which queue lengths are balanced through mean-field interaction using an {\it interaction function}. We propose a summary of their model and then we expose how our current model can be used to generalize the ideas in \cite{Daw+Tang+Zhao2005}. 

Consider a system consisting of $N$ queues with a mean-field interaction. At $t = 0$, for $1 \leq n \leq N$, the arrival rate to the $n$-th queue occurs according to $\zeta_{X_n(0)}$, and the service rate at queue $n$ is $\vartheta_{X_n(0)}$. Let $h(x):\mathbb{R}_+\times\mathbb{R}_+\rightarrow \mathbb{R}$ be a continuous nondecreasing {\it interaction function} satisfying certain regularity conditions (see \cite[p. 339]{Daw+Tang+Zhao2005}). This function allows to capture the mean-field interaction between queues as follows: for each queue $n= 1, 2,\ldots,N$, the arrival rate at time $t$ is given by  $\zeta_{X_j(t)}- h(X_j (t),\langle \mu^N (t)(dx),x\rangle)$, where $\mu^N (t)=\frac{1}{N}\sum_{j=1}^N\delta_{X_j(t)}$ is the empirical measure corresponding to the $N$ queues at time $t$. Note that $\langle \mu^N (t)(dx),x\rangle=\frac{1}{N}\sum_{j=1}^NX_j(t)$ is the mean queue length of the $N$ queues at time $t$. Roughly speaking, the arrival rate at each queue depends on the current size of the queue and on the mean size of its neighbors (which is the entire set of queues in this case). The authors then studied the performance of such a network in terms of limiting results as $N$ goes to infinity.

The model proposed in the current  work can be seen as a generalization of the model in \cite{Daw+Tang+Zhao2005} to heterogeneous queuing networks, namely, to block-structured networks. To see this, let consider the graph $\mathcal{G}=(\mathcal{V},\Xi)$ as a queuing network where the particles (nodes) are finite-buffer server queues of maximum size $K$ (arbitrary large), and the corresponding states $(X_{n,j}(t),X_{n,j}(t),n\in C_j^c,m\in C_j^p,1\leq j\leq r, t\geq 0)$ represent the number of customers waiting in each queue at each time $t$. Again, at $t = 0$, for $1 \leq n \leq N$, the arrival rate to the $n$-th queue occurs according to $\zeta_{X_{n,j}(0)}$, and the service rate at queue $n$ is $\vartheta_{X_{n,j}(0)}$. Since the network now is heterogeneous, the mean-field interaction is local. Thus, the arrival rate at a central node queue $n\in C_j^c$ at time $t$ is given by $\zeta_{X_{n,j}(t)}-h(X_{n,j}(t),\langle \mu^N_j (t)(dx),x\rangle)$ whereas the arrival rate at a peripheral node queue $n\in C_j^p$ at time $t$ is given by $\zeta_{X_{n,j}(t)}-h(X_{n,j}(t),\langle \mu_{n,j}^N (t)(dx),x\rangle)$, with $\mu^N_j (t)$ and $\mu_{n,j}^N (t)$ are the local empirical measures respectively given by $(\ref{block-decomp})$ and $(\ref{periph-decomp})$. The service rate $\vartheta_{X_{n,j}(t)}$ at each queue $1\leq n\leq N$ depends only on the queue size $X_{n,j}(t)$ at time $t$. Hence, the transition rates $\lambda^c_{z,z'}$ and $\lambda^p_{z,z'}$ are specified as follows: 

\begin{itemize}
\item The size $X_{n,j}(t)$  of each central queue $n\in C^c_j$ at time $t$ goes from $z$ to $z'$ at rate

\begin{align*}
\lambda^c_{z,z'}=\left\{\begin{tabular}{l l}
                         $\zeta_{X_{n,j}(t)}-h^c\big(X_{n,j}(t),\frac{1}{N_j}\sum\limits_{n\in C_j}X_{n,j}(t)\big)$&\mbox{if $z'=z+1$ and $z'\leq K$ }\\
                         $\vartheta_{X_{n,j}(t)}$& \mbox{if $z'=z-1$ and $X_{n,j}(t)\geq 1$}\\
                         $-\sum\limits_{y\neq z}\lambda^c_{z,y}$&$\mbox{if $z'=z$}$\\
                         $0$&\mbox{otherwise}
                         \end{tabular}\right.    
\end{align*}

\item The size $X_{n,j}(t)$  of each peripheral queue $n\in C^p_j$ at time $t$ goes from $z$ to $z'$ at rate
\begin{align*}
\lambda^p_{z,z'}=\left\{\begin{tabular}{l l}
                         $\zeta_{X_{n,j}(t)}-h^p\big(X_n(t),\frac{1}{deg(n)+1}\sum\limits_{\substack{m:(n,m)\in \Xi\\ 1\leq k\leq r}}(X_{m,k}(t))\big)$&\mbox{if $z'=z+1$ and $z'\leq K$ }\\
                         $\vartheta_{X_{n,j}(t)}$& \mbox{if $z'=z-1$ and $X_{n,j}(t)\geq 1$}\\
                         $-\sum\limits_{y\neq z}\lambda^p_{z,y}$&$\mbox{if $z'=z$}$\\
                         $0$&\mbox{otherwise}
                         \end{tabular}\right.    
\end{align*}

 \end{itemize} 
Note that the sparse graph topologies have been considered in applications as responses to some issues encountered when trying to implement load balancing protocols. In particular, many service systems are geographically constrained, therefore, when a task arrives at any specific server, it might be impossible to collect instantaneous state information from all the servers. Besides, executing a task commonly involves the use of some data, and storing such data for all possible tasks on all servers requires an excessive amount of storage capacity. The use of sparser graph topologies is then considered such that tasks that arrive at a specific server can only be forwarded, following a specific load balancing scheme, to the servers that possess the data required to process the tasks. In other words, a specific server can only interact with its neighbors in a suitable sparse topology. See, e.g., \cite{Budh+Mukh+Wu2019} and the references therein for more insights about the subject.   

The results obtained in the current work allow us to understand the behavior of such systems when the size $N$ of the network goes to infinity. In particular, the multi-chaotic property established in Theorem $\ref{theo2}$ tells us that the queue lengths at any finite collection of tagged servers are statistically asymptotically independent, and the queue-length process
for each server converges in distribution to the corresponding McKean-Vlasov process given by $(\ref{limit-syst})$. Also,  Condition \ref{Cond-prin} and  Remark \ref{rem-regul-cond} tell us that the  multi-chaotic result holds even when the peripheral sub-graph is not complete, which translates that we can achieve the same asymptotic performance with much fewer connections between the peripheral nodes than when all the peripheral nodes are connected.

\subsubsection{Multi-population SIS Epidemics}

 The susceptible-infected-susceptible  (SIS)  epidemic model, happened to be also useful to model the spread of information in networks given that the two phenomena are closely related. The SIS model can be summarized as follows. Consider an information/infection diffusion across a population. A member that has a copy of the information/disease is said to be {\it infected} and a member that does not have a copy of the information/infection is said to be {\it susceptible}. When an infected member comes in contact with a susceptible one, the former transmits a copy of the content (disease) to the latter, and the latter gets infected. Moreover, an infected member may spontaneously get rid of the content, a phenomenon called {\it curing}, to become susceptible again. 

In both epidemiology and networks information diffusion, the population is often formed by isolated sub-populations whose members are highly interacting, connected between them by few members. One might think, e.g., of countries as isolated communities connected by tourists traveling across the globe, or of interactions in social media which often happen in almost closed communities with a few influential members interacting across groups. Our model allows studying the spreading dynamics of an information/disease among the members of a population structured as separate communities.

Consider a population consisting of $r$ isolated communities and a "mobile" community. The members of each isolated community interact only among themselves and with members of the mobile community. Thus, there is no direct interaction between members of different communities. However, an indirect inter-community interaction holds via the set of mobile members. This idea was used in \cite{Akhi+Alt+Sund2019} where the authors considered an optimal control problem to find the optimal resource allocation strategy to maximize information spread over the multi-community population. Their objective was to obtain a good tradeoff between the information spread in the network and the use of system resources. 

Let $\mathcal{Z}=\{0,1\}$ be the state space that indicates  whether  the  particle  is  {\it susceptible} $(=0)$ or {\it infected} $(=1)$. Recall that for a given block,  a central node interacts with all other central nodes and with the peripheral nodes of the same block. However, there is no direct interaction between the central nodes of a given block with nodes from other blocks. One might then think of the central nodes of each block as an isolated community that interacts with other communities only through the peripheral nodes, which in turn can be seen as the mobile community. Note that in contrast to \cite{Akhi+Alt+Sund2019}, the central nodes of a given block interact only with the peripheral (mobile) nodes of the same block, and not with all the peripheral/mobile nodes, as stipulated in \cite{Akhi+Alt+Sund2019}. Also, the assumptions introduced in Condition \ref{Cond-prin} and  Remark \ref{rem-regul-cond} bring us to a situation where not all the peripheral nodes interact with each other and thus, the interaction graph for the peripheral members is not complete. Nevertheless, the fact that the multi-chaotic property holds under Condition \ref{Cond-prin} tells us that the systems with full connections between the peripheral components and those with fewer connections, illustrated by \ref{Cond-prin} and  Remark \ref{rem-regul-cond}, are asymptotically equivalent. This is of interest for example in resource allocation problems where we attribute a cost to each connection. This is nonetheless beyond the scope of the present paper. 

Denote by $X_{n,j}(t)$, for $n\in C_j^c$ (resp. $n\in C_j^p$), the state ("susceptible" or "infected") of the $n$-th central (resp. peripheral) particle in the $j$-th community. Two central members of the same community $j$ come in contact with each other at rate $\gamma_j$. Peripheral and central nodes from the same community interact with each other at a rate $\nu_j$. Two connected peripheral nodes come in contact with each other at a rate $\eta$. Finally, an infected node in community $j$ spontaneously gets rid of the content at a rate $\zeta_j$. Therefore, the transition rates $\lambda^c_{z,z'}$ and $\lambda^p_{z,z'}$ are specified as follows, which sums up the dynamics we are interested in:

\begin{itemize}
\item The state $X_{n,j}(t)$  of each central member $n\in C^c_j$ at time $t$ goes from $z$ to $z'$ at rate 
 \begin{align*}
\lambda^c_{z,z'}=\left\{\begin{tabular}{l l}
                         $N^c_j(\mu_j^{c,N}(t)(1))\gamma_j+N^p_j(\mu_j^{p,N}(t)(1))\nu_j$&\mbox{if $z=0$ and $z'=1$}\\
                         $\zeta_j$& \mbox{if $z=1$ and $z'=0$}\\
                         $-\sum\limits_{y\neq z}\lambda^c_{z,y}$&$\mbox{if $z'=z$}$\\
                         $0$&\mbox{otherwise}
                         \end{tabular}\right.    
\end{align*}

\item The state $X_{n,j}(t)$  of each peripheral (mobile) member $n\in C^p_j$ at time $t$ goes from $z$ to $z'$ at rate  
\begin{align*}
\lambda^p_{z,z'}=\left\{\begin{tabular}{l l}
                         $N^c_j(\mu_j^{c,N}(t)(1))\nu_j+\sum\limits_{\substack{m\in\cup_{k} C_k^p\\(n,m)\in\Xi }}X_{m,k}(t)\eta$&\mbox{if $z=0$ and $z'=1$}\\
                         $\zeta_j$& \mbox{if $z=1$ and $z'=0$}\\
                         $-\sum\limits_{y\neq z}\lambda^p_{z,y}$&$\mbox{if $z'=z$}$\\
                         $0$&\mbox{otherwise}
                         \end{tabular}\right.    
\end{align*}

\end{itemize}

Notice that the large deviations properties established in Section $\ref{large-dev-sec}$ is a step forward to study the large time behavior of such systems. Indeed, the large deviations of the empirical measure established in Theorem \ref{large-dev-emp-proc} can be used to investigate the large deviations of the invariant measure, from which one can study the large time behavior of the system and the related phenomena such as metastability and convergence to the invariant measure. This will be part of future research. The interested reader  can consult, e.g. \cite{Yaso+Sund2019},\cite{Hwang+Sheu90} and \cite{Freid+Wentz2012}.

\subsection{Notations and conventions}

Let $(\mathbb{S}, d)$ be a Polish space. For any $x\in\mathcal{D}([0,T],\mathbb{S})$ we denote $\|x\|_T=\sup_{0\leq s\leq T}\|x(s)\|$. For any $y\in\mathbb{S}^d$ we denote $\|y\|=\max (y_1, \ldots, y_d)$. Given two measures $\mu,\nu\in\mathcal{M}(\mathbb{S})$, the bounded-Lipschitz metric $d_{BL}(\cdot,\cdot)$ is defined by
\begin{align}
d_{BL}(\mu,\nu)= \sup_{g\in Lip(\mathbb{S})} \big| \langle \mu,g\rangle-\langle\nu,g\rangle \big|,
\label{bound-lip}
\end{align}
where
\begin{align*}
Lip(\mathbb{S})=\left\{g\in C_b(\mathbb{S}):\sup_{x\in\mathbb{S}} |g(x)|\leq 1,\sup_{x\neq y} \frac{|g(x)-g(y)|}{d(x,y)} \leq 1\right\}.
\end{align*}

Recall that  the bounded-Lipschitz metric metrizes the weak convergence of probability measures on $\mathbb{S}$ with respect to bounded continuous test functions $C_b(\mathbb{S})$. For $p\geq 1$, let $\mathcal{M}_p(\mathbb{S})$ be the collection of all probability measures on $\mathbb{S}$ with finite  $p$-th moment. Then, for any $\mu$  and  $\nu$  in $\mathcal{M}_p(\mathbb{S})$, the $p$-th Wasserstein distance between $\mu$ and $\nu$ is defined as
\begin{align*}
     \mathcal{W}_p( \mu, \nu)= \left(\inf_{\gamma\in\Gamma (\mu,\nu)} \int_{E\times E} d(x,y)^p d \gamma(x,y)\right)^{1/p},
\end{align*}
where  $\Gamma (\mu ,\nu )$ denotes the collection of all measures on $\mathbb{S}\times \mathbb{S}$  with marginals $\mu$  and $\nu$. Moreover, for $M_1,M_2$ in $\mathcal{M}_1\big(\mathcal{D}([0,T],\mathbb{S})\times\cdots\mathcal{D}([0,T],\mathbb{S})\big)$, the $p$-th Wasserstein distance between $M_1$ and $M_2$ is given by
\begin{equation*}
\mathcal{W}_{p,T}(M_1, M_2)= \inf\bigg\{\big[\mathbb{E}\|Y_1-Y_2\|^p_T\big]^{1/p}: Y_1, Y_2 \in \mathcal{D}([0,T],\mathbb{S})\times\cdots\mathcal{D}([0,T],\mathbb{S}), M_1= \mathcal{L}(Y_1), M_2= \mathcal{L}(Y_2)\bigg\}.
\end{equation*}

\section{Law of large numbers and propagation of chaos}
\label{LLN-POC-sec}

We study in this section the behavior of the system when the number of particles $N$ tends to infinity. We use the convention that $N$ goes to infinity when both $\min_{1\leq j\leq r}N_j^c$ and $\min_{1\leq j\leq r}N_j^p$ goes to infinity. In particular, we investigate the law of large numbers and the propagation of chaos properties. Since the system is heterogeneous, we describe its state at each time $t$ using the following \textit{empirical measure vector}
\begin{align*}
\mu^N(t)=\left(\mu_1^{c,N}(t),\mu_1^{p,N}(t),\cdots,\mu_r^{c,N}(t),\mu_r^{p,N}(t)\right),
\end{align*}
where for each $1\leq j\leq r$, $\mu_j^{c,N}(t)$ (resp. $\mu_j^{p,N}(t)$) is the empirical measure describing the states of the central (resp. peripheral) nodes of the $j$-th block at time $t$. The symmetry between the central nodes within the same block suggests the convergence of the empirical measure $\mu_j^{c,N}$ towards the distribution $\mu_j^c\in\mathcal{M}_1(\mathcal{D}([0,T],\mathcal{Z}))$, the solution of an appropriate limiting process. Also, despite the lack of symmetry between the peripheral nodes of the same block, we will introduce some conditions under which $\mu_j^{p,N}$ weakly converge  when $N$ gets large. Hence, the empirical vector $\mu^N$ should converge weakly to $\mu$ where
\begin{align*}
\mu&=\left(\mu_1^{c},\mu_1^{p},\cdots,\mu_r^{c},\mu_r^{p}\right)=\left(\mathcal{L}( \bar{X}^c_{n,1}),\mathcal{L}( \bar{X}^p_{n,1}),\ldots, \mathcal{L}( \bar{X}^c_{n,r}),\mathcal{L}( \bar{X}^p_{n,r})\right)\in\big(\mathcal{M}_1(\mathcal{D}([0,T],\mathcal{Z}))\big)^{2r},
\end{align*}
and $\big((\bar{X}^c_{n,j}(t), \bar{X}^p_{n,j}(t),t\geq 0),1\leq j\leq r\big)$ is the solution of the following system of stochastic differential equations
\begin{equation}
\begin{split}
 \bar{X}^{c}_{n,j}(t)&= \bar{X}^{c}_{n,j}(0)+\int_{[0,t]\times\mathbb{R}_+}\sum_{(z,z')\in\mathcal{E}}\mathds{1}_{ \bar{X}^{c}_{n,j}(s-)=z}(z'-z)\mathds{1}_{\left[0,\lambda^c_{z,z'}\left(\upsilon_j(s-)\right)\right]}(y)\mathcal{N}_{n,j}^c(ds,dy),\\
 \bar{X}^{p}_{n,j}(t)&= \bar{X}^{p}_{n,j}(0)+\int_{[0,t]\times\mathbb{R}_+}\sum_{(z,z')\in\mathcal{E}}\mathds{1}_{ \bar{X}^{p}_{n,j}(s-)=z}(z'-z)\mathds{1}_{\left[0,\lambda^p_{z,z'}\left(\upsilon_{n,j}(s-)\right)\right]}(y)\mathcal{N}_{n,j}^p(ds,dy).
\label{limit-syst}
\end{split}
\end{equation}
Here,  $\upsilon_j(t)$ and $\upsilon_{n,j}(t)$ are vectors defined by
\begin{equation}
\begin{split}
\upsilon_j(t)&=(\mu^c_j(t),\mu^p_j(t),p_j^c,p_j^p), \\
\upsilon_{n,j}(t)&=(\mu^c_j(t),\mu^p_1(t),\ldots,\mu^p_r(t),\alpha_j^c,q_{j,1},\ldots,q_{j,r}),
\label{ups-funct}
\end{split}
\end{equation}
and $p_j^c,p_j^p,\alpha_j^c, q_{j1},\ldots,q_{jr}\in (0,1)$ are  parameters  satisfying,
\begin{align*}
 p_j^c+p_j^p=1\quad\text{and}\quad \alpha_j^c+q_{j1}+\cdots+q_{jr}=1\quad\text{for each} \quad 1\leq j\leq r,
\end{align*}
which will later  be chosen appropriately (see Condition $\ref{Cond-prin}$). The link between the initial conditions of the systems $(\ref{SDE-rep})$ and $(\ref{limit-syst})$ will be introduced in the sequel. Observe that the solution of $(\ref{limit-syst})$  depends on the distribution of the process itself and not only on its sample path up to time $t$. Thus, the system $(\ref{limit-syst})$ is McKean-Vlasov.

We now introduce the conditions under which the results of this section hold.
\begin{cond}
\begin{enumerate}
\item There exist, for all $(z,z')\in\mathcal{E}$, some measurable functions $\gamma^c_{z,z'}:\mathcal{Z}\rightarrow\mathbb{R}^+$ and $\gamma^p_{z,z'}:\mathcal{Z}\rightarrow\mathbb{R}^+$, such that:
  \begin{itemize}
     \item For any  probability measures $\nu,\mu\in\mathcal{M}_1(\mathcal{Z})$ and  any real numbers $0<a_1,a_2<1$ with $a_1+a_2=1$ we have
\begin{equation}
\lambda_{z,z'}^c(\nu,\mu,a_1,a_2)=a_1\int_{\mathcal{Z}}\gamma_{z,z'}^c(x)\nu (dx)+ a_2\int_{\mathcal{Z}}\gamma_{z,z'}^p(x)\mu (dx).
 \label{lambda-c-func}
\end{equation}
     \item  For any $\nu,\mu_1,\ldots,\mu_r\in\mathcal{M}_1(\mathcal{Z})$ and any real numbers $0<a,b_1,\ldots,b_r<1$ such that $a+b_1+\cdots+b_r=1$ we have
\begin{align}
\lambda^p_{z,z'}(\nu,\mu_1,\ldots,\mu_r,a,b_1,\dots,b_r)= a\int_{\mathcal{Z}}\gamma_{z,z'}^c(x)\nu (dx)+b_1\int_{\mathcal{Z}}\gamma_{z,z'}^p(x)\mu_1 (dx)+\cdots+b_r\int_{\mathcal{Z}}\gamma_{z,z'}^p(x)\mu_r (dx).
\label{lambda-p-func}
\end{align}
\end{itemize}
\item For each block $1\leq j\leq r$, there exist $p_j^c,p_j^p\in(0,1)$ such that, as $N\rightarrow\infty$,
\begin{align}
\frac{N_j^p}{N_j}\rightarrow p_j^p,\quad \frac{N_j^c}{N_j}\rightarrow p_j^c\quad\text{and}\quad p^p_j+p^c_j=1 .
\label{p-regul}
\end{align}

\item For each $1\leq j\leq r$, there exist $\alpha_j^c,q_{j1},\ldots,q_{jr}\in (0,1)$ with $\alpha_j^c+q_{j1}+\cdots+q_{jr}=1$ such that, for each peripheral node $n\in C_j^p$, the following conditions hold for all $1\leq i\leq r$ with $i\neq j$

\begin{align}
\lim_{N\rightarrow\infty}\bigg|\frac{N_j^p}{deg(n)+1}-q_{jj}\bigg|= 0,\quad\lim_{N\rightarrow\infty}\bigg|\frac{N_j^c}{deg(n)+1}-\alpha_{j}^c\bigg|= 0\quad\text{and}\quad\lim_{N\rightarrow\infty}\bigg|\frac{M_i^n}{deg(n)+1}-q_{ji}\bigg|= 0.
\label{cond-regul}
\end{align}
\end{enumerate}
\label{Cond-prin}
\end{cond}

\begin{figure}
\center
\includegraphics[scale=0.3]{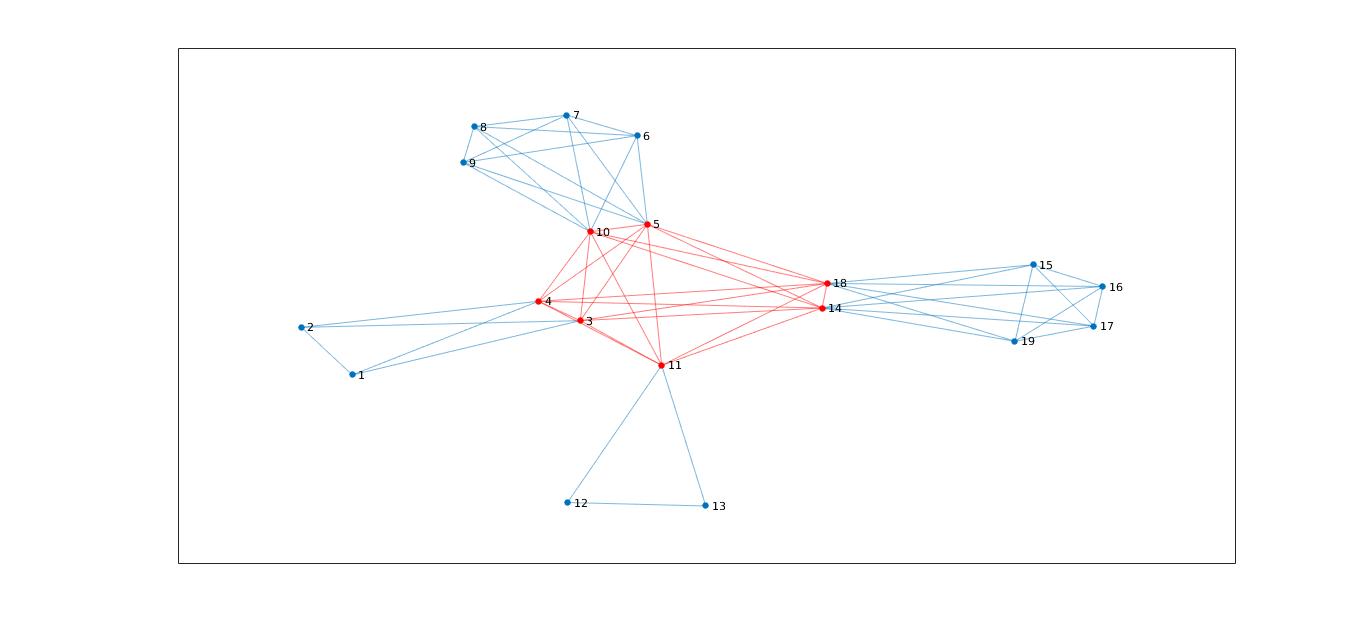}
\caption{Example of block-structured graph with a complete peripheral sub-graph. Here we have a 4-blocks-structured graph linked by a set of peripheral nodes. For the first block the set of central nodes is $C_1^c=\{1,2\}$ and the set of peripheral nodes is $C_1^p=\{3,4\}$. The set of all peripheral nodes of the graph is given by the set of nodes $C^p=\{3,4,5,10,11,14,18\}$. }
\label{comp-subgraph}
\end{figure}

\begin{rem}
\begin{enumerate}
\item Since $\mathcal{Z}$ is a finite state space, the functions $\gamma_{z,z'}^c$ and $\gamma_{z,z'}^p$ are bounded on $\mathcal{Z}$. Moreover, since $\mathcal{Z}\subset \mathbb{N}$ and that every bounded function on $\mathbb{N}$ is automatically Lipschitz, $\gamma_{z,z'}^c$ and $\gamma_{z,z'}^p$ are also Lipschitz. Denote by $\bar{\gamma}>0$ the maximum bound and by $L_{\gamma}$ the maximum Lipschitz coefficient of the functions $\gamma_{z,z'}^c$ and $\gamma_{z,z'}^p$ for all $(z,z')\in\mathcal{E}$.

\item Condition $(\ref{cond-regul})$  is satisfied for example if, for each peripheral node $n\in C_j^p$,   $M_i^n/N_i^p\rightarrow 1$ as $N\rightarrow\infty$ for all $1\leq i\leq r$. Indeed, under this condition we define
\begin{align}
\alpha_j^c&=\lim_{N\rightarrow\infty}\frac{N_j^c}{N_j^c+N_1^p+\cdots+N_r^p},\quad\forall 1\leq j\leq r,\\
 q_{j,i}&=\lim_{N\rightarrow\infty}\frac{N_{i}^p}{N_j^c+N_1^p+\cdots+N_r^p},\quad\forall 1\leq j,i\leq r,
 \end{align}
and thus, one can easily verify that,  as $N\rightarrow\infty$,
\begin{align}
\frac{M_i^n}{deg(n)+1}\rightarrow q_{ji}\quad\text{and}\quad\frac{N_j^c}{deg(n)+1}\rightarrow\alpha_{j}^c.
\end{align}

\item A special case where condition  $(\ref{cond-regul})$ is satisfied is when the peripheral sub-graph is complete, that is, when all peripheral nodes are connected to each other (see Figure $\ref{comp-subgraph}$). In such a case, the peripheral nodes of the same block are homogeneous.

\item Even though condition $(\ref{cond-regul})$ is somehow restrictive since it imposes the peripheral sub-graph to be dense, the construction of the model allows to have very different degrees in each block. Moreover, one can investigate less restrictive conditions on the central nodes by relaxing the conditions of completeness and replace it with adequate regularity conditions.
\item One might contrast Condition $(\ref{cond-regul})$ with some existing conditions in the literature. Consider for example the condition imposed in \cite{Budh+Mukh+Wu2019} for a supermarket model on sparse graphs to asymptotically behave as on a clique. The condition in \cite{Budh+Mukh+Wu2019} relies on local properties of the graph by imposing direct neighbors of any node to have asymptotically similar degrees, see \cite[Cond. 1 (ii)]{Budh+Mukh+Wu2019}. This condition is nevertheless violated by our model. Indeed, Condition $(\ref{cond-regul})$ allows central and peripheral nodes of the same block to have very different degrees, even though they are neighbors, which goes beyond \cite[Cond. 1 (ii)]{Budh+Mukh+Wu2019}. In addition, under our condition, $d_{\max} (G)/d_{\min}(G)$ should not goes to $1$ as $N\rightarrow\infty$ neither $\max_j \left|\left(d_{\min} (C_j)/d_{\max}(C_j)\right)-1\right|$ goes to zero as proposed in \cite[Rem. 1]{Budh+Mukh+Wu2019} (here $deg_{\min}(C_j)$ and $deg_{\max}(C_j)$ refer to the min and max degrees of nodes within the same block j). In that sense, the graph we are considering in the present work is sparser than the ones covered by \cite[Cond. 1 (ii)]{Budh+Mukh+Wu2019}. Another condition to contrast with is the one proposed in \cite{Delattre+al2016} under which an $n$-dimensional diffusion system converges to a limiting Fokker-Plank equation, see \cite[eqn. (1.1) and (1.3)]{Delattre+al2016}. Note that \cite[(1.5),(1.7)]{Delattre+al2016} impose global regularity conditions in the sense that the degrees of all the nodes should converge to the same limit, which is clearly not imposed by Condition $(\ref{cond-regul})$.    
\item While the current paper considers the case of static graphs, one can investigate the case where the underlying graph topology is random. For example, it is of interest for some applications to have the scenario where the connections between the peripheral nodes are allowed to be random. One then can search for the adequate conditions to impose on the edges dynamics for the propagation of chaos property to hold. This is however goes beyond the scope of the current work. 
\end{enumerate}
\label{rem-regul-cond}
\end{rem}

\subsection{Existence and uniqueness results}

The following theorem shows the existence and uniqueness of the solution of the McKean-Vlasov limiting equation $(\ref{limit-syst})$.
\begin{theo} 
 Suppose that Condition \ref{Cond-prin} holds true. Then, for a given initial condition $\big((\bar{X}_n^{c}(0),\bar{X}_m^{p}(0)),n\in C_j^c,m\in C_j^p;1\leq j\leq r\big)$, the McKean-Vlasov system $(\ref{limit-syst})$ has a unique solution over any finite interval of time $[0,T]$. In addition, the solution of the limiting system $(\ref{limit-syst})$ depends continuously on the initial condition in the following sense: if $(\bar{X}(t),t\in[0,T])$ and $(\bar{X}'(t),t\in[0,T])$ are two solutions of $(\ref{limit-syst})$ with two different initial conditions $(\bar{X}(0))$ and $(\bar{X}'(t))$, respectively, then there exist a constant $A_T$ such that
\begin{equation}
\begin{split}
\mathbb{E}[\|\bar{X}-\bar{X}'\|_T]\leq 2E[\|\bar{X}(0)-\bar{X}'(0)\|]e^{A_T}.
\end{split}
\label{cont-init}
\end{equation}
\label{theo1}
\end{theo}

\proof For $1\leq j\leq r$, let, with slight abuse of notations,
\[
    e_{j,c}: (x^c_1, x_1^p,\ldots, x_r^c,x_r^p) \in (\mathcal{Z}^{2r})\rightarrow x_j^c\in\mathcal{Z},
\]
and
\[
    e_{j,p}: (x^c_1, x_1^p,\ldots, x_r^c,x_r^p) \in (\mathcal{Z}^{2r})\rightarrow x_j^p\in\mathcal{Z}
\]
  be the $c$-th and the $p$-th component of the $j$-th projection, respectively. Moreover, for $t\leq T$, denote $p_t: f\in \mathcal{D}([0,T],\mathcal{Z}^{2r})\rightarrow f(t)\in\mathcal{Z}^{2r}$.
  
For $M\in\mathcal{M}_1\big(\mathcal{D}([0,T],\mathcal{Z}^{2r})\big)$, denote $M (t) = M \circ p_t^{-1}$. Define $\psi$ and $\phi$ the mappings that associate to $M$ respectively the solution and its corresponding law of the system starting at $\bar{X}_0=(\bar{X}_0^{1,c},\bar{X}_0^{1,p},\ldots,\bar{X}_0^{r,c},\bar{X}_0^{r,p})$ by, at each $t\in (0,T]$, 
\begin{equation}
\begin{split}
 \bar{X}^{c}_{j}(t)&= \bar{X}^{c}_{j}(0)+\int_{[0,t]\times\mathbb{R}_+}\sum_{(z,z')\in\mathcal{E}}\mathds{1}_{ \bar{X}^{c}_{j}(s-)=z}(z'-z)\mathds{1}_{\left[0,\lambda^c_{z,z'}\left(\upsilon_j(s-)\right)\right]}(y)\mathcal{N}_{j}^c(ds,dy),\\
 \bar{X}^{p}_{j}(t)&= \bar{X}^{p}_{j}(0)+\int_{[0,t]\times\mathbb{R}_+}\sum_{(z,z')\in\mathcal{E}}\mathds{1}_{ \bar{X}^{p}_{j}(s-)=z}(z'-z)\mathds{1}_{\left[0,\lambda^p_{z,z'}\left(\upsilon_{n,j}(s-)\right)\right]}(y)\mathcal{N}_{j}^p(ds,dy),
\label{limit-syst-test}
\end{split}
\end{equation}
 for $1\leq j\leq r$ where $\mu_j^c(t)=M(t)\circ e_{j,c}^{-1}$, $\mu_j^p(t)=M(t)\circ e_{j,p}^{-1}$ and the vectors $\upsilon_j(t)$ and $\upsilon_{n,j}(t)$ are given by $(\ref{ups-funct})$. Thus, $\psi (M)=\big((\bar{X}^c_{j}(t), \bar{X}^p_{j}(t),t\geq 0),1\leq j\leq r\big)$ and $\phi (M)=\mathcal{L}\big((\bar{X}^c_{j}(t), \bar{X}^p_{j}(t),t\geq 0),1\leq j\leq r\big)$. Observe that if $\bar{X}$ is a solution of $(\ref{limit-syst})$, then its law is a fixed point of $\phi$. Conversely, if $M$  is a fixed point of $\phi$ for the system $(\ref{limit-syst-test})$, then  the corresponding solution $\psi (M)$ defines a solution of the limiting system $(\ref{limit-syst})$. The idea is then to prove the existence of a fixed point of $\phi$.

Take  $M_1,M_2\in\mathcal{M}_1\big(\mathcal{D}([0,T],\mathcal{Z}^{2r})\big)$. Set $\bar{X}_1=(\bar{X}^{1,c}_1,\bar{X}^{1,p}_1\ldots,\bar{X}^{1,c}_r,\bar{X}^{1,p}_r)=\psi (M_1)$ and $\bar{X}_2=(\bar{X}^{2,c}_1,\bar{X}^{2,p}_1\ldots,\bar{X}^{2,c}_r,\bar{X}^{2,p}_r)=\psi (M_2)$. Thus, $\mathcal{L}(\bar{X}_1)=\phi (M_1)$ and $\mathcal{L}(\bar{X}_2)=\phi (M_2)$. Moreover, for all $t\in [0,T]$ denote  $\mu_1(t)=(\mu_1^{1,c}(t),\mu_1^{1,p}(t),\ldots,\mu_r^{1,c}(t),\mu_r^{1,p}(t))$ and $\mu_2(t)=(\mu_1^{2,c}(t),\mu_1^{2,p}(t),\ldots,\mu_r^{2,c}(t),\mu_r^{2,p}(t))$ with $\mu_j^{1,c}(t)=M_1(t)\circ e_{j,c}^{-1}$, $\mu_j^{1,p}(t)=M_1(t)\circ e_{j,p}^{-1}$, $\mu_j^{2,c}(t)=M_2(t)\circ e_{j,c}^{-1}$ and $\mu_j^{2,p}(t)=M_2(t)\circ e_{j,p}^{-1}$ for $1\leq j\leq r$. According to $(\ref{ups-funct})$ we introduce the following notations:
\begin{equation}
\begin{split}
\upsilon^{1,c}_j(t)&=(\mu^{1,c}_j(t),\mu^{1,p}_j(t),p_j^{c},p_j^{p}), \\
\upsilon^{2,c}_j(t)&=(\mu^{2,c}_j(t),\mu^{2,p}_j(t),p_j^{c},p_j^{p}), \\
\upsilon^{1,p}_{n,j}(t)&=(\mu^{1,c}_j(t),\mu^{1,p}_1(t),\ldots,\mu^{1,p}_r(t),\alpha_j^c,q_{j,1},\ldots,q_{j,r}),\\
\upsilon^{2,p}_{n,j}(t)&=(\mu^{2,c}_j(t),\mu^{2,p}_1(t),\ldots,\mu^{2,p}_r(t),\alpha_j^c,q_{j,1},\ldots,q_{j,r}).
\label{ups-not}
\end{split}
\end{equation}

We first prove that $\phi$ is a contraction mapping on $\mathcal{M}_1\big(\mathcal{D}([0,T],\mathcal{Z}^{2r})\big)$, that is, for any $t\in[0,T]$,
\begin{align}
\mathcal{W}_{1,t}\bigg(\phi(M_1),\phi(M_2)\bigg)\leq C(t)\mathbb{E}\bigg[\int_{[0,t]} \mathcal{W}_{1,s}(M_1,M_2) ds\bigg].
\label{recur-ineq}
\end{align}

Indeed, for any $1\leq j\leq r$ we have that
\begin{equation}
\begin{split}
\|\bar{X}^{1,c}_{j}-\bar{X}^{2,c}_{j}\|_t&\leq \int_{[0,t]\times\mathbb{R}_+}\bigg|\sum_{(z,z')\in\mathcal{E}}(z'-z)\bigg\{\mathds{1}_{\bar{X}^{1,c}_{j}(s-)=z}\mathds{1}_{\left[0,\lambda^c_{z,z'}(\upsilon^{1,c}_j(s))\right]}(y) \\
                                             &\qquad\qquad\qquad\qquad -\mathds{1}_{\bar{X}^{2,c}_{j}(s-)=z}\mathds{1}_{\left[0,\lambda^c_{z,z'}(\upsilon^{2,c}_{j}(s))\right]}(y)\bigg\}\bigg| \mathcal{N}_{j}^c(ds,dy).
\end{split}
\end{equation}

Using martingale argument (see $(\ref{ineq-mart1})$) and taking the expectation one gets, by adding and subtracting terms (see $(\ref{ineq-2})$), for any $t\in [0,T]$,
\begin{equation}
\begin{split}
\mathbb{E}\left[\|\bar{X}^{1,c}_{j}-\bar{X}^{2,c}_{j}\|_t\right]&\leq K\mathbb{E}\bigg[\int_{[0,t]}\sum_{(z,z')\in\mathcal{E}}\bigg|\big(\mathds{1}_{\bar{X}^{1,c}_{j}(s)=z}-\mathds{1}_{\bar{X}^{2,c}_{j}(s)=z}\big)\lambda^c_{z,z'}(\upsilon^{1,c}_j(s))
                                             \\
                                             &\qquad\qquad\qquad\qquad\qquad+\bigg(\lambda^c_{z,z'}(\upsilon^{1,c}_j(s))-\lambda^c_{z,z'}(\upsilon^{2,c}_j(s))\bigg) \bigg|ds\bigg].
\end{split}
\end{equation}

Recall the definition of the functions $\lambda^c_{z,z'}$ in $(\ref{lambda-c-func})$. Given that $\mu_j^c(t)$ and $\mu_j^p(t)$ are probability measures and using the boundedness of the functions $\gamma^c_{z,z'}$ and $\gamma^p_{z,z'}$, one easily get that
\begin{align}
\lambda^c_{z,z'}(\upsilon^{1,c}_j(s))\leq \bar{\gamma},
\end{align}
and
 \begin{equation}
 \begin{split}
 \bigg|\lambda^c_{z,z'}(\upsilon_j^{1,c}(s))-\lambda^c_{z,z'}(\upsilon_j^{2,c}(s))\bigg|&\leq p_j^c\bar{\gamma}\bigg|\langle 1, \mu_j^{1,c}(s)-\mu_j^{2,c}(s)\rangle\bigg|+p_j^p\bar{\gamma}\bigg|\langle 1, \mu_j^{1,p}(s)-\mu_j^{2,p}(s)\rangle\bigg|\\
 &\leq p_j^c\bar{\gamma} d_{BL} \bigg(\mu_j^{1,c}(s),\mu_j^{2,c}(s)\bigg)+p_j^p\bar{\gamma}d_{BL} \bigg(\mu_j^{1,p}(s),\mu_j^{2,p}(s)\bigg),
 \end{split}
 \end{equation}
 Therefore we obtain,
\begin{equation}
\begin{split}
\mathbb{E}\left[\|\bar{X}^{1,c}_{j}-\bar{X}^{2,c}_{j}\|_t\right]&\leq K|\mathcal{E}|\bar{\gamma}\mathbb{E}^0\bigg[\int_{[0,t]}\bigg(\big|\bar{X}^{1,c}_{n,j}(s)-\bar{X}^{2,c}_{n,j}(s)\big|
                                             \\
                                             &\qquad\qquad\qquad\qquad\qquad+p_j^c d_{BL} \big(\mu_j^{1,c}(s),\mu_j^{2,c}(s)\big)+p_j^pd_{BL} \big(\mu_j^{1,p}(s),\mu_j^{2,p}(s)\big)\bigg) ds\bigg].
\label{BL-ineq1}
\end{split}
\end{equation}

Using $(\ref{lambda-p-func})$ and the same previous steps we find, for any $1\leq j\leq r$,
\begin{equation}
\begin{split}
\mathbb{E}\left[\|\bar{X}_{j}^{1,p}-\bar{X}_{j}^{2,p} \|_t\right]&\leq K|\mathcal{E}|\bar{\gamma}\mathbb{E}\bigg[\int_{[0,t]}\bigg(\big|\bar{X}^{1,p}_{j}(s)-\bar{X}^{2,p}_{j}(s)\big|+\alpha_j^c d_{BL} \big(\mu_j^{1,c}(s),\mu_j^{2,c}(s)\big) \\
                                             &\qquad\qquad\qquad\qquad+q_{j1}d_{BL} \big(\mu_1^{1,p}(s),\mu_1^{2,p}(s)\big)+\cdots+q_{jr}d_{BL} \big(\mu_r^{1,p}(s),\mu_r^{2,p}(s)\big)\bigg)ds\bigg].
                                             \label{BL-ineq2}
\end{split}
\end{equation}

On one hand, we have, from the Kantorovich-Rubinstein theorem, that for $1\leq j\leq r$ and $\alpha\in\{c,p\}$,
\begin{align}
d_{BL} \big(\mu_j^{1,\alpha}(s),\mu_j^{2,\alpha}(s)\big)=\mathcal{W}_1\big(\mu_j^{1,\alpha}(s),\mu_j^{2,\alpha}(s)\big).
\label{Kant-Rub}
\end{align}

On the other hand, we can easily verify that
\begin{align}
\mathcal{W}_1\big(\mu_j^{1,\alpha}(s),\mu_j^{2,\alpha}(s)\big)\leq\mathcal{W}_{1,s}\big(M_1,M_2\big).
\label{W-ineq}
\end{align}

Therefore, using $(\ref{Kant-Rub})$ and $(\ref{W-ineq})$, and taking the supremum over $1\leq j\leq r$ in $(\ref{BL-ineq1})$ and $(\ref{BL-ineq2})$ we obtain
\begin{equation}
\begin{split}
\mathbb{E}\bigg[\sup_{1\leq j\leq r}\|\bar{X}^{1,c}_{j}-\bar{X}^{2,c}_{j}\|_t\bigg]&\leq K|\mathcal{E}|\bar{\gamma}\mathbb{E}\bigg[\int_{[0,t]}\bigg(\sup_{1\leq j\leq r}\big\|\bar{X}^{1,c}_{n,j}-\bar{X}^{2,c}_{n,j}\big\|_s + \mathcal{W}_{1,s}(M_1,M_2)\bigg) ds\bigg],
\label{BL-ineq11}
\end{split}
\end{equation}

\begin{equation}
\begin{split}
\mathbb{E}\bigg[\sup_{1\leq j\leq r}\|\bar{X}^{1,p}_{j}-\bar{X}^{2,p}_{j}\|_t\bigg]&\leq K|\mathcal{E}|\bar{\gamma}\mathbb{E}\bigg[\int_{[0,t]}\bigg(\sup_{1\leq j\leq r}\big\|\bar{X}^{1,p}_{n,j}-\bar{X}^{2,p}_{n,j}\big\|_s + \mathcal{W}_{1,s}(M_1,M_2)\bigg) ds\bigg].
\label{BL-ineq22}
\end{split}
\end{equation}

Adding side by side the two last inequalities and applying the Gr\"{o}nwall's lemma we obtain

\begin{equation}
\begin{split}
\mathbb{E}\bigg[\sup_{1\leq j\leq r}\|\bar{X}^{1,c}_{j}-\bar{X}^{2,c}_{j}\|_t+\sup_{1\leq j\leq r}\|\bar{X}^{1,p}_{j}-\bar{X}^{2,p}_{j}\|_t\bigg]&\leq K|\mathcal{E}|\bar{\gamma}\mathbb{E}\bigg[\int_{[0,t]} \mathcal{W}_{1,s}(M_1,M_2) ds\bigg]e^{K|\mathcal{E}|\bar{\gamma}t}.
\end{split}
\end{equation}
Hence,
\begin{equation}
\begin{split}
\mathbb{E}\bigg[\|\bar{X}^{1}-\bar{X}^{2}\|_t\bigg]&\leq C(t)\mathbb{E}\bigg[\int_{[0,t]} \mathcal{W}_{1,s}(M_1,M_2) ds\bigg],
\end{split}
\end{equation}
with $C(t)=K|\mathcal{E}|\bar{\gamma}e^{K|\mathcal{E}|\bar{\gamma}t}$. Observe that from the definition of the Wasserstein distance we have that
\begin{align*}
\mathcal{W}_{1,t}\bigg(\phi(M_1),\phi(M_2)\bigg)\leq \mathbb{E}\bigg[\|\bar{X}^{1}-\bar{X}^{2}\|_t\bigg],
\end{align*}
from which we deduce $(\ref{recur-ineq})$.

We now consider the following recursive scheme:
\begin{itemize}
\item $M_0\in\mathcal{M}_1\big(\mathcal{D}([0,T],\mathcal{Z}^{2r})\big)$;
\item $M_{k+1}=\phi (M_k),\quad k\geq 0$.
\end{itemize}

By iterating the formula in $(\ref{recur-ineq})$ and using the fact that $\mathcal{W}_{1,t}(M_1,M_0)$ is increasing in $t$ we find that
\begin{align*}
\mathcal{W}_{1,t}(M_{k+2},M_{k+1})\leq \frac{(tC(t))^k}{k!}\mathcal{W}_{1,t}(M_1,M_0)
\end{align*}
for all $k\geq 0$. It is easy to verify that  $\mathcal{W}_{1,t}(M_1,M_0)<\infty$ and thus the sequence $\{M_k\}_{k\geq 0}$ is a Cauchy sequence. Note that the space $\mathcal{M}_1\big(\mathcal{D}([0,T],\mathcal{Z}^{2r})\big)$ endowed with the Wasserstein distance $\mathcal{W}_{p,T}$ is complete (see \cite{Bolley2008}). Hence the sequence $\{M_k\}_{k\geq 0}$ converges to some measure $M$ in $\mathcal{M}_1\big(\mathcal{D}([0,T],\mathcal{Z}^{2r})\big)$ which is a fixed point of $\phi$ on $\mathcal{M}_1\big(\mathcal{D}([0,T],\mathcal{Z}^{2r})\big)$. This proves existence of solution of equation $(\ref{limit-syst-test})$ and thus for the equation in $(\ref{limit-syst})$. Uniqueness follows using again $(\ref{recur-ineq})$ and Gr\"{o}nwall's lemma.

Define by $(\bar{X}^1(t))=(\bar{X}^{1,c}_{n,j}(t),\bar{X}^{1,p}_{m,j}(t),n\in C_j^c,m\in C_j^p, 1\leq j\leq r)$ and $(\bar{X}^2(t))=(\bar{X}^{2,c}_{n,j}(t),\bar{X}^{2,p}_{m,j}(t),n\in C_j^c,m\in C_j^p, 1\leq j\leq r)$ the two solutions of $(\ref{limit-syst})$ with respective initial conditions $(\bar{X}^1(0))$ and $(\bar{X}^2(0))$. Denote by $\mu_j^{1,c}(t)=\mathcal{L}(\bar{X}^{1,c}_{n,j}(t))$ and $\mu_j^{1,p}(t)=\mathcal{L}(\bar{X}^{1,p}_{m,j}(t))$ the probability measures corresponding to the first solution. Similarly, denote by  $\mu_j^{2,c}(t)=\mathcal{L}(\bar{X}^{2,c}_{n,j}(t))$ and $\mu_j^{2,p}(t)=\mathcal{L}(\bar{X}^{2,p}_{m,j}(t))$ the the probability measures corresponding to the second solution. Using again the notation in $(\ref{ups-not})$ we find that, for any $1\leq j\leq r$ and $t\in [0,T]$,
\begin{equation}
\begin{split}
\|\bar{X}^{1,c}_{n,j}-\bar{X}^{2,c}_{n,j}\|_t&\leq |\bar{X}^{1,c}_{n,j}(0)-\bar{X}^{2,c}_{n,j}(0)|\\
                                        &\qquad+\int_{[0,t]\times\mathbb{R}_+}\bigg|\sum_{(z,z')\in\mathcal{E}}(z'-z)\bigg\{\mathds{1}_{\bar{X}^{1,c}_{n,j}(s-)=z}\mathds{1}_{\left[0,\lambda^c_{z,z'}(\upsilon^{1,c}(s))\right]}(y) \\
                                             &\qquad\qquad\qquad\qquad\qquad\qquad\qquad -\mathds{1}_{\bar{X}^{2,c}_{n,j}(s-)=z}\mathds{1}_{\left[0,\lambda^c_{z,z'}(\upsilon^{2,c}(s))\right]}(y)\bigg\}\bigg| \mathcal{N}_{n,j}^c(ds,dy).
\end{split}
\end{equation}
Using martingale argument (see $(\ref{ineq-mart1})$), taking the conditional expectation $E^0$ given $(\bar{X}^1(0),\bar{X}^2(0))$  and finally adding and subtracting terms (see $(\ref{ineq-2})$) we find that, for $t\in [0,T]$,
\begin{equation}
\begin{split}
\mathbb{E}^0\left[\|\bar{X}^{1,c}_{n,j}-\bar{X}^{2,c}_{n,j}\|_t\right]&\leq |\bar{X}^{1,c}_{n,j}(0)-\bar{X}^{2,c}_{n,j}(0)|\\
&+K\mathbb{E}^0\bigg[\int_{[0,t]}\sum_{(z,z')\in\mathcal{E}}\bigg|\bigg(\mathds{1}_{\bar{X}^{1,c}_{n,j}(s)=z}-\mathds{1}_{\bar{X}^{2,c}_{n,j}(s)=z}\bigg)\lambda^c_{z,z'}(\upsilon^{1,c}(s))
                                             \\
                                             &\qquad\qquad\qquad\qquad\qquad+\bigg(\lambda^c_{z,z'}(\upsilon^{1,c}(s)(s))-\lambda^c_{z,z'}(\upsilon^{2,c}(s))\bigg) \bigg|ds\bigg].
\end{split}
\end{equation}
Given that $\mu_j^{1,c}(t)$ and $\mu_j^{1,p}(t)$ are probability measures then, by using $(\ref{lambda-c-func})$ and the boundedness of the function $\gamma^c_{z,z'}$ and $\gamma^p_{z,z'}$ we find that
\begin{align}
\lambda^c_{z,z'}(\upsilon^{1,c}(s))\leq \bar{\gamma},
\end{align}
and by the Lipschitz property of the functions $\gamma^c_{z,z'}$ we find that
 \begin{align}
 \bigg|\lambda^c_{z,z'}(\upsilon^{1,c}(s))-\lambda^c_{z,z'}(\upsilon^{2,c}(s))\bigg|\leq p_j^c\bar{\gamma}\mathbb{E}[\|\bar{X}^{1,c}_{n,j}(s)-\bar{X}^{2,c}_{n,j}(s)\|]+p_j^p\bar{\gamma}\mathbb{E}[\|\bar{X}^{1,p}_{n,j}(s)-\bar{X}^{2,p}_{n,j}(s)\|].
 \end{align}
 Therefore we obtain that
 \begin{equation}
\begin{split}
\mathbb{E}^0\left[\|\bar{X}^{1,c}_{n,j}-\bar{X}^{2,c}_{n,j}\|_t\right]&\leq |\bar{X}^{1,c}_{n,j}(0)-\bar{X}^{2,c}_{n,j}(0)|+K\bar{\gamma}|\mathcal{E}|\int_{[0,t]}\bigg(\mathbb{E}^0[\|\bar{X}^{1,c}_{n,j}(s)-\bar{X}^{2,c}_{n,j}(s)\|]+p_j^c\mathbb{E}[\|\bar{X}^{1,c}_{n,j}(s)-\bar{X}^{2,c}_{n,j}(s)\|]
                                             \\
                                             &\qquad\qquad\qquad\qquad\qquad\qquad\qquad\qquad+p_j^p\mathbb{E}[\|\bar{X}^{1,p}_{n,j}(s)-\bar{X}^{2,p}_{n,j}(s)\|]\bigg) ds.
\end{split}
\end{equation}
Taking the expectation in the two sides of the last inequality, and recalling that $p_j^c+p_j^p=1$ we obtain
 \begin{equation}
\begin{split}
\mathbb{E}\left[\|\bar{X}^{1,c}_{n,j}-\bar{X}^{2,c}_{n,j}\|_t\right]&\leq \mathbb{E}[|\bar{X}^{1,c}_{n,j}(0)-\bar{X}^{2,c}_{n,j}(0)|]+4K\bar{\gamma}|\mathcal{E}|\int_{[0,t]}\bigg(\mathbb{E}[\|\bar{X}^{1,c}_{n,j}(s)-\bar{X}^{2,c}_{n,j}(s)\|]+\mathbb{E}[\|\bar{X}^{1,p}_{n,j}(s)-\bar{X}^{2,p}_{n,j}(s)\|]\bigg) ds.
\end{split}
\end{equation}
Taking the maximum over $n\in C_j^c$ and over $1\leq j\leq r$ gives to us
 \begin{equation}
\begin{split}
\mathbb{E}[\max_{\substack{n\in C_j^c\\1\leq j\leq r}}\|\bar{X}^{1,c}_{n,j}-\bar{X}^{2,c}_{n,j}\|_t]&\leq \mathbb{E}[\max_{\substack{n\in C_j^c\\1\leq j\leq r}}|\bar{X}^{1,c}_{n,j}(0)-\bar{X}^{2,c}_{n,j}(0)|]\\
&+4K\bar{\gamma}|\mathcal{E}|\int_{[0,t]}\bigg(\mathbb{E}[\max_{\substack{n\in C_j^c\\1\leq j\leq r}}\|\bar{X}^{1,c}_{n,j}(s)-\bar{X}^{2,c}_{n,j}(s)\|]+\mathbb{E}[\max_{\substack{n\in C_j^p\\1\leq j\leq r}}\|\bar{X}^{1,p}_{n,j}(s)-\bar{X}^{2,p}_{n,j}(s)\|]\bigg) ds.
\label{cont-init1}
\end{split}
\end{equation}
Using similar arguments we find that, for any $n\in C_j^p$ with $1\leq j\leq r$
\begin{equation}
\begin{split}
\mathbb{E}^0\left[\|\bar{X}_{n,j}^{1,p}-\bar{X}_{n,j}^{2,p} \|_t\right]&\leq |\bar{X}_{n,j}^{1,p}(0)-\bar{X}_{n,j}^{2,p}(0)|+K\mathbb{E}^0\bigg[\int_{[0,t]}\sum_{(z,z')\in\mathcal{E}}\bigg|\bigg(\bar{X}^{1,p}_{n,j}(s)-\bar{X}^{2,p}_{n,j}(s)\bigg)\lambda^p_{z,z'}\left(\upsilon^{1,p}(s)\right) \\
                                             &\qquad\qquad\qquad\qquad\qquad\qquad\qquad\qquad\qquad\qquad+\bigg(\lambda^p_{z,z'}\left(\upsilon^{1,p}_{n,j}(s)\right)-\lambda^p_{z,z'}\left(\upsilon^{2,p}(s)\right)\bigg) \bigg|\bigg]ds.
\end{split}
\end{equation}
By $(\ref{lambda-p-func})$ and the Lipschitz boundedness property of the functions $\gamma^c_{z,z'}$ and $\gamma^p_{z,z'}$ we find that
\begin{equation}
\begin{split}
\mathbb{E}^0\left[\|\bar{X}_{n,j}^{1,p}-\bar{X}_{n,j}^{2,p} \|_t\right]&\leq |\bar{X}_{n,j}^{1,p}(0)-\bar{X}_{n,j}^{2,p}(0)|+K\bar{\gamma}|\mathcal{E}|\int_{[0,t]}\bigg(\mathbb{E}^0[\|\bar{X}^{1,p}_{n,j}(s)-\bar{X}^{2,p}_{n,j}(s)\|] \\
                                             &+ \alpha_j^c\mathbb{E}[\|\bar{X}^{1,c}_{n,j}(s)-\bar{X}^{2,c}_{n,j}(s)\|]+q_{j1}\mathbb{E}[\|\bar{X}^{1,p}_{n,1}(s)-\bar{X}^{2,p}_{n,1}(s)\|]+\cdots+q_{jr}\mathbb{E}[\|\bar{X}^{1,p}_{n,r}(s)-\bar{X}^{2,p}_{n,r}(s)\|]\bigg)ds.
\end{split}
\end{equation}
Recall that $\alpha_j^c+q_{j1}+\cdots+q_{jr}=1$. Then by taking the expectation in the two sides of the last inequality we obtain
\begin{equation}
\begin{split}
\mathbb{E}\left[\|\bar{X}_{n,j}^{1,p}-\bar{X}_{n,j}^{2,p} \|_t\right]&\leq \mathbb{E}[|\bar{X}_{n,j}^{1,p}(0)-\bar{X}_{n,j}^{2,p}(0)|]+4K\bar{\gamma}|\mathcal{E}|\int_{[0,t]}\bigg(\mathbb{E}[\|\bar{X}^{1,c}_{n,j}(s)-\bar{X}^{2,c}_{n,j}(s)\|]\\
&\qquad\qquad\qquad\qquad\qquad+\mathbb{E}[\|\bar{X}^{1,p}_{n,1}(s)-\bar{X}^{2,p}_{n,1}(s)\|]+\cdots+\mathbb{E}[\|\bar{X}^{1,p}_{n,r}(s)-\bar{X}^{2,p}_{n,r}(s)\|]\bigg)ds.
\end{split}
\end{equation}
Taking the maximum over $n\in C_j^p$ and then over $1\leq j\leq r$ we find that
\begin{equation}
\begin{split}
\mathbb{E}[\max_{\substack{n\in C_j^p\\1\leq j\leq r}}\|\bar{X}_{n,j}^{1,p}-\bar{X}_{n,j}^{2,p} \|_t]&\leq E[\max_{\substack{n\in C_j^p\\1\leq j\leq r}}|\bar{X}_{n,j}^{1,p}(0)-\bar{X}_{n,j}^{2,p}(0)|]\\
  &\quad+4K\bar{\gamma}|\mathcal{E}|\int_{[0,t]}\bigg(r\mathbb{E}[\max_{\substack{n\in C_j^p\\1\leq j\leq r}}\|\bar{X}^{1,p}_{n,j}(s)-\bar{X}^{2,p}_{n,j}(s)\|]+\mathbb{E}[\max_{\substack{n\in C_j^c\\1\leq j\leq r}}\|\bar{X}^{1,c}_{n,j}(s)-\bar{X}^{2,c}_{n,j}(s)\|]\bigg)ds  \\
  &\leq E[\max_{\substack{n\in C_j^p\\1\leq j\leq r}}|\bar{X}_{n,j}^{1,p}(0)-\bar{X}_{n,j}^{2,p}(0)|]\\
  &\quad+4K\bar{\gamma}|\mathcal{E}|(1+r)\int_{[0,t]}\bigg(\mathbb{E}[\max_{\substack{n\in C_j^p\\1\leq j\leq r}}\|\bar{X}^{1,p}_{n,j}(s)-\bar{X}^{2,p}_{n,j}(s)\|]+\mathbb{E}[\max_{\substack{n\in C_j^c\\1\leq j\leq r}}\|\bar{X}^{1,c}_{n,j}(s)-\bar{X}^{2,c}_{n,j}(s)\|]\bigg)ds.
  \label{cont-init2}
\end{split}
\end{equation}
Now $(\ref{cont-init1})$ and $(\ref{cont-init2})$ together leads to
\begin{equation}
\begin{split}
\mathbb{E}[\|\bar{X}^1-\bar{X}^2\|_t]&=\mathbb{E}\bigg[\max\bigg(\max_{\substack{n\in C_j^c\\1\leq j\leq r}}\|\bar{X}^{1,c}_{n,j}-\bar{X}^{2,c}_{n,j}\|_t,\max_{\substack{n\in C_j^p\\1\leq j\leq r}}\|\bar{X}_{n,j}^{1,p}-\bar{X}_{n,j}^{2,p} \|_t\bigg)\bigg]\\
&\leq \mathbb{E}\bigg[\max_{\substack{n\in C_j^c\\1\leq j\leq r}}\|\bar{X}^{1,c}_{n,j}-\bar{X}^{2,c}_{n,j}\|_t+\max_{\substack{n\in C_j^p\\1\leq j\leq r}}\|\bar{X}_{n,j}^{1,p}-\bar{X}_{n,j}^{2,p} \|_t\bigg]\\
 &\leq E[\max_{\substack{n\in C_j^p\\1\leq j\leq r}}|\bar{X}_{n,j}^{1,c}(0)-\bar{X}_{n,j}^{2,c}(0)|]+ E[\max_{\substack{n\in C_j^p\\1\leq j\leq r}}|\bar{X}_{n,j}^{1,p}(0)-\bar{X}_{n,j}^{2,p}(0)|]\\
  &\quad+4K\bar{\gamma}|\mathcal{E}|(2+r)\int_{[0,t]}\bigg(\mathbb{E}[\max_{\substack{n\in C_j^p\\1\leq j\leq r}}\|\bar{X}^{1,p}_{n,j}(s)-\bar{X}^{2,p}_{n,j}(s)\|]+\mathbb{E}[\max_{\substack{n\in C_j^c\\1\leq j\leq r}}\|\bar{X}^{1,c}_{n,j}(s)-\bar{X}^{2,c}_{n,j}(s)\|]\bigg)ds\\
   &\leq 2E[\|\bar{X}^1(0)-\bar{X}^2(0)\|]+8K\bar{\gamma}|\mathcal{E}|(2+r)\int_{[0,t]}\bigg(\mathbb{E}[\|\bar{X}^1(s)-\bar{X}^2(s)\|]\bigg)ds.
\end{split}
\end{equation}
Finally, by the Gr\"{o}nwall's we conclude that
\begin{equation}
\begin{split}
\mathbb{E}[\|\bar{X}^1-\bar{X}^2\|_t]\leq 2E[\|\bar{X}^1(0)-\bar{X}^2(0)\|]e^{8K\bar{\gamma}|\mathcal{E}|(2+r)t}.
\end{split}
\end{equation}
Defining $A_t=8K\bar{\gamma}|\mathcal{E}|(2+r)t$ leads to $(\ref{cont-init})$. The theorem is proved. 
\carre

\subsection{Weak convergence for converging initial condition}

In this section we establish the weak convergence of $(\ref{SDE-rep})$ towards the limiting McKean-Vlasov system $(\ref{limit-syst})$. First, recall the notions of multi-exchangeability and multi-chaoticity introduced in \cite{Graha2008}.

 \begin{defi}
A sequence of random variables $(X_{n,k}, 1 \leq n \leq N_k, 1 \leq k \leq K)$ indexed by $N=(N_k)\in \mathbb{N}^K$ is said to be multi-exchangeable if its law is invariant under permutation of the indexes within the classes, that is, for $1 \leq k \leq K$ and any permutations $\sigma_k$ of $\{1,\ldots, N_k\}$, the following equality holds in distribution
\begin{align*}
(X_{\sigma_k(n),k} , 1 \leq n \leq N_k, 1 \leq k \leq K) \overset{dist}{=} (X_{n,k}, 1 \leq n \leq N_k, 1 \leq k \leq K).
\end{align*}
A sequence of random variables $(X_{n,k}, 1\leq n\leq N_k, 1\leq k\leq K)$ indexed by $N=(N_k)\in \mathbb{N}^K$ is $P_1\otimes\cdots\otimes P_K$-multi-chaotic if, for any $m\geq 1$, the convergence in distribution
\begin{align*}
\lim_{N\rightarrow\infty}(X^N_{n,k}, 1\leq n \leq m, 1 \leq k \leq K) \overset{dist}{=} P_1^{\otimes m}\otimes\cdots\otimes P_K^{\otimes m}
\end{align*}
holds for the topology of the uniform convergence on compact sets, where $P_k$, for $1\leq k \leq K$, is a probability distribution on $\mathbb{R}_+$,  and with the convention that $N$ goes to infinity when $\min N_k$ goes to infinity.
\label{defi-multi-chaos}
 \end{defi}

The following theorem is the main result for the weak convergence.
\begin{theo}
Suppose that Condition \ref{Cond-prin} holds true. Moreover, suppose that the initial conditions $(X_n^{c,N}(0),X_m^{p,N}(0),n\in C_j^c,m\in C_j^p;1\leq j\leq r)$ are multi-exchangeable and $\nu^{1,c}\otimes\nu^{1,p}\cdots \nu^{r,c}\otimes\nu^{r,p}$-multi-chaotic. Then, for any $t\in [0,T]$, as $N\rightarrow\infty$,
\begin{align}
\max_{1\leq j\leq r}\max_{n\in C_j^c}\mathbb{E}\left[\|X^{c}_{n,j}-\bar{X}^{c}_{n,j}\|_t \right]+\max_{1\leq j\leq r}\max_{n\in C_j^p}\mathbb{E}\left[\|X^{p}_{n,j}-\bar{X}^{p}_{n,j}\|_t \right]\rightarrow 0,
\label{L_1-norm}
\end{align}
 and the sequence of processes $\big((X_n^{c,N}(t),X_m^{p,N}(t),t\geq 0),n\in C_j^c,m\in C_j^p;1\leq j\leq r\big)$, the solutions of the SDE $(\ref{SDE-rep})$ with initial conditions $(X_n^{c,N}(0),X_m^{p,N}(0),n\in C_j^c,m\in C_j^p;1\leq j\leq r)$, is  $P_{\bar{X}}$-multi-chaotic, where $P_{\bar{X}}=\mu_{1}^c\otimes\mu_{1}^p\cdots \mu^c_{r}\otimes \mu_{r}^p$ is the distribution of the process $\big((\bar{X}_n^{c}(t),\bar{X}_m^{p}(t),t\geq 0),n\in C_j^c,m\in C_j^p;1\leq j\leq r\big)$, the solution of the limiting SDE $(\ref{limit-syst})$ with initial distribution $\nu^{1,c}\otimes\nu^{1,p}\cdots \nu^{r,c}\otimes\nu^{r,p}$.
\label{theo2}
\end{theo}

Before proceeding to the proof, we recall, without a proof, an elementary result on (conditionally) i.i.d. random variables.
\begin{lem}
 Let $\{S_i:i=1,\ldots,n\}$ be a collection of $\mathbb{S}$-valued random variables defined on some probability space $(\Omega,\mathcal{F},\mathbb{P})$, where $\mathbb{S}$ is some Polish space. Suppose that $S_i$ for $i=1,\ldots,n$ are conditionally i.i.d. given some $\sigma$-algebra $\mathcal{G}\subset\mathcal{F}$. Then, for any $k\in\mathbb{N}$, there exists a positive finite constant $0<a_k<\infty$ such that,
 \begin{align}
 \sup_{\|f\|_{\infty}\leq 1}\mathbb{E}\bigg|\frac{1}{n}\sum_{i=1}^n(S_i-\mathbb{E}[S_i|\mathcal{G}])\bigg|^k\leq \frac{a_k}{N^{k/2}}.
 \label{moment-bound}
\end{align}
\end{lem}

\paragraph{Proof of Theorem \ref{theo2}.}
 We use a coupling method. Let $Y(t)=\big((Y_n^{c}(t),Y_m^{p}(t)),n\in C_j^c,m\in C_j^p;1\leq j\leq r\big)$ be the solution of the limiting SDE $(\ref{limit-syst})$ with the initial conditions of the process $X(t)=\big((X_n^{c}(t),X_m^{p}(t)),n\in C_j^c,m\in C_j^p;1\leq j\leq r\big)$, the solution of the SDE $(\ref{SDE-rep})$, given by
\begin{align*}
Y(0)=\big((Y_n^{c}(0),Y_m^{p}(0)),n\in C_j^c,m\in C_j^p;1\leq j\leq r\big)=\big((X_n^{c}(0),X_m^{p}(0)),n\in C_j^c,m\in C_j^p;1\leq j\leq r\big).
\end{align*}
Also, we define the processes $Y(t)$   and  $X(t)$ on the same probability space by taking the same sequences of Poisson random measures $\{\mathcal{N}^{c}_{n,j}\}$ (resp. $\{\mathcal{N}^{p}_{n,j}\}$) in both cases. We first prove that the two processes are asymptotically close, that is, for any $t\in [0,T]$,

\begin{align}
\max_{1\leq j\leq r}\max_{n\in C_j^c}\mathbb{E}\left[\|X^{c}_{n,j}-Y^{c}_{n,j}\|_t \right]+\max_{1\leq j\leq r}\max_{n\in C_j^p}\mathbb{E}\left[\|X^{p}_{n,j}-Y^{p}_{n,j}\|_t \right]\rightarrow 0.
\label{Asymp-close}
\end{align}

 We treat the central and peripheral nodes in two separate steps.

\subparagraph*{Step 1.} Fix $1\leq j\leq r$. For each central node $n\in C_j^c$ and any $t\in [0,T]$,

\begin{equation}
\begin{split}
\mathbb{E}\left[\|X_{n,j}^{c}-Y_{n,j}^{c}\|_t\right]&=\mathbb{E}\left[\sup_{0\leq s\leq t}\big|X_{n,j}^{c}(s)-Y_{n,j}^{c}(s)\big|\right]\\
                                             &\leq\mathbb{E}\bigg[\bigg|\int_{[0,t]\times\mathbb{R}_+}\sum_{(z,z')\in\mathcal{E}}\mathds{1}_{X^{c}_{n,j}(s-)=z}(z'-z)\mathds{1}_{\left[0,\lambda^c_{z,z'}(\upsilon_j^{N}(s))\right]}(y)\mathcal{N}_{n,j}^c(ds,dy) \\
                                             &\qquad\qquad -\int_{[0,t]\times\mathbb{R}_+}\sum_{(z,z')\in\mathcal{E}}\mathds{1}_{Y^{c}_{n,j}(s-)=z}(z'-z)\mathds{1}_{\left[0,\lambda^c_{z,z'}(\upsilon_j(s))\right]}(y)\mathcal{N}_{n,j}^c(ds,dy) \bigg| \bigg]\\
                                             &\leq\mathbb{E}\bigg[\int_{[0,t]\times\mathbb{R}_+}\bigg|\sum_{(z,z')\in\mathcal{E}}(z'-z)\bigg\{\mathds{1}_{X^{c}_{n,j}(s-)=z}\mathds{1}_{\left[0,\lambda^c_{z,z'}(\upsilon_j^{N}(s))\right]}(y) \\
                                             &\qquad\qquad -\mathds{1}_{Y^{c}_{n,j}(s-)=z}\mathds{1}_{\left[0,\lambda^c_{z,z'}(\upsilon_j(s))\right]}(y)\bigg\}\mathcal{N}_{n,j}^c(ds,dy) \bigg| \bigg].
\end{split}
\label{ineq-1}
\end{equation}
Denote by $\mathcal{F}_t$ the filtration generated by the Poisson random measures and defined by,
\begin{align*}
\mathcal{F}_t=\sigma\big\langle \mathcal{N}_{n,j}^c(A\times B):n\in C_j^c \cup \mathcal{N}_{m,j}^p(A\times B):m\in C_j^p, A\in\mathcal{B}(\mathbb{R}_+),B\in\mathcal{B}([0,T]))     \big\rangle.
\end{align*}
Then,  $X^{c}_{n,j}(t)$ and $Y^{c}_{n,j}(t)$ are adapted to the filtration $\mathcal{F}_t$. Therefore, the two processes

\begin{equation}
\begin{split}
&\int_{[0,t]\times\mathbb{R}_+}\sum_{(z,z')\in\mathcal{E}}\mathds{1}_{X^{c}_{n,j}(s-)=z}(z'-z)\mathds{1}_{\left[0,\lambda^c_{z,z'}(\upsilon_j^{N}(s))\right]}(y)[\mathcal{N}_{n,j}^c(ds,dy)-dsdy],\\
&\int_{[0,t]\times\mathbb{R}_+}\sum_{(z,z')\in\mathcal{E}}\mathds{1}_{Y^{c}_{n,j}(s-)=z}(z'-z)\mathds{1}_{\left[0,\lambda^c_{z,z'}(\upsilon_j(s))\right]}(y)[\mathcal{N}_{n,j}^c(ds,dy)-dsdy],
\label{ineq-mart1}
\end{split}
\end{equation}
are $\mathcal{F}_t$-martingales. Furthermore, $(\ref{ineq-1})$ reduces to,
\begin{equation}
\begin{split}
\mathbb{E}\left[\|X_{n,j}^{c}-Y_{n,j}^{c} \|_t\right]&\leq \mathbb{E}\bigg[\int_{[0,t]}\bigg|\sum_{(z,z')\in\mathcal{E}}(z'-z)\bigg\{\mathds{1}_{X^{c}_{n,j}(s)=z}\lambda^c_{z,z'}(\upsilon_j^{N}(s)) -\mathds{1}_{Y^{c}_{n,j}(s)=z}\lambda^c_{z,z'}(\upsilon_j(s))\bigg\} \bigg|ds \bigg].
\end{split}
\end{equation}
Recall that $K=|\mathcal{Z}|$ is the number of colors. By adding and subtracting terms we obtain
   \begin{equation}
\begin{split}
\mathbb{E}&\left[\|X_{n,j}^{c}-Y_{n,j}^{c} \|_t\right]\\
&\leq \mathbb{E}\bigg[\int_{[0,t]}\bigg|\sum_{(z,z')\in\mathcal{E}}(z'-z)\bigg\{\mathds{1}_{X^{c}_{n,j}(s)=z}\lambda^c_{z,z'}(\upsilon_j^{N}(s)) -\mathds{1}_{Y^{c}_{n,j}(s)=z}\lambda^c_{z,z'}(\upsilon_j(s))\bigg\} \bigg|ds \bigg] \\
                                             &\leq   K\mathbb{E}\bigg[\int_{[0,t]}\bigg|\sum_{(z,z')\in\mathcal{E}}(z'-z)\bigg\{\mathds{1}_{X^{c}_{n,j}(s)=z}\lambda^c_{z,z'}(\upsilon_j^{N}(s))-\mathds{1}_{Y^{c}_{n,j}(s)=z}\lambda^c_{z,z'}(\upsilon_j(s))\bigg\} \bigg|ds \bigg] \\
                                             &\leq K\mathbb{E}\bigg[\int_{[0,t]}\bigg|\sum_{(z,z')\in\mathcal{E}}\bigg\{\mathds{1}_{X^{c}_{n,j}(s)=z}\lambda^c_{z,z'}(\upsilon_j^{N}(s))-\mathds{1}_{Y^{c}_{n,j}(s)=z}\lambda^c_{z,z'}(\upsilon_j^{N}(s))\\
                                             &\qquad\qquad\qquad\qquad\qquad+\mathds{1}_{Y^{c}_{n,j}(s)=z}\lambda^c_{z,z'}(\upsilon_j^{N}(s))-\mathds{1}_{Y^{c}_{n,j}(s)=z}\lambda^c_{z,z'}(\upsilon_j(s))\bigg\} \bigg|ds \bigg]  \\
                                             &\leq  K\int_{[0,t]}\sum_{(z,z')\in\mathcal{E}}\mathbb{E}\bigg[\bigg|\bigg(\mathds{1}_{X^{c}_{n,j}(s)=z}-\mathds{1}_{Y^{c}_{n,j}(s)=z}\bigg)\lambda^c_{z,z'}(\upsilon_j^{N}(s))
                                             +\bigg(\lambda^c_{z,z'}(\upsilon_j^{N}(s))-\lambda^c_{z,z'}(\upsilon_j(s))\bigg) \bigg|\bigg]ds.
\end{split}
\label{ineq-2}
\end{equation}

The idea now is to find a bound for the right-hand side of $(\ref{ineq-2})$. Let us start by the second term. Again by adding and subtracting terms we get

\begin{equation}
\begin{split}
\mathbb{E}\bigg[\bigg|\lambda^c_{z,z'}(\upsilon_j^{N}(s))-\lambda^c_{z,z'}(\upsilon_j(s)) \bigg|\bigg]&= \mathbb{E}\bigg[\bigg| \bigg(\frac{1}{N_j}\sum_{n\in  C_j^c}\gamma_{z,z'}^c(X_{n,j}^{c}(s))+\frac{1}{N_j}\sum_{n\in  C_j^p}\gamma_{z,z'}^p(X_{n,j}^{p}(s))\bigg)\\
&\qquad-\bigg(\frac{1}{N_j}\sum_{n\in  C_j^c}\gamma_{z,z'}^c(Y_{n,j}^{c}(s))+\frac{1}{N_j}\sum_{n\in  C_j^p}\gamma_{z,z'}^p(Y_{n,j}^{p}(s))\bigg)\\
&\qquad+\bigg(\frac{1}{N_j}\sum_{n\in  C_j^c}\gamma_{z,z'}^c(Y_{n,j}^{c}(s))+\frac{1}{N_j}\sum_{n\in  C_j^p}\gamma_{z,z'}^p(Y_{n,j}^{p}(s))\bigg)\\
&\qquad-\bigg(p_j^c\int_{\mathcal{Z}}\gamma^c_{z,z'}(x)\mu_j^c(s)ds+p_j^p\int_{\mathcal{Z}}\gamma^p_{z,z'}(x)\mu_j^p(s)ds\bigg)\bigg|\bigg]\\
&\leq \mathbb{E}\bigg[\bigg| \bigg(\frac{1}{N_j}\sum_{n\in  C_j^c}\gamma_{z,z'}^c(X_{n,j}^{c}(s))+\frac{1}{N_j}\sum_{n\in  C_j^p}\gamma_{z,z'}^p(X_{n,j}^{p}(s))\bigg)\\
&\qquad-\bigg(\frac{1}{N_j}\sum_{n\in  C_j^c}\gamma_{z,z'}^c(Y_{n,j}^{c}(s))+\frac{1}{N_j}\sum_{n\in  C_j^p}\gamma_{z,z'}^p(Y_{n,j}^{p}(s))\bigg)\bigg|\bigg]\\
&+\mathbb{E}\bigg[\bigg|\bigg(\frac{1}{N_j}\sum_{n\in  C_j^c}\gamma_{z,z'}^c(Y_{n,j}^{c}(s))+\frac{1}{N_j}\sum_{n\in  C_j^p}\gamma_{z,z'}^p(Y_{n,j}^{p}(s))\bigg)\\
&\qquad-\bigg(p_j^c\int_{\mathcal{Z}}\gamma^c_{z,z'}(x)\mu_j^c(s)ds+p_j^p\int_{\mathcal{Z}}\gamma^p_{z,z'}(x)\mu_j^p(s)ds\bigg)\bigg|\bigg].
\end{split}
\label{ineq-3}
\end{equation}
From the Lipschitz property of the functions $\gamma^c_{z,z'}$ and $\gamma^p_{z,z'}$, the first expectation in the right-hand side of $(\ref{ineq-3})$ is bounded as follows:
\begin{equation}
  \begin{split}
  \mathbb{E}\bigg[&\bigg| \bigg(\frac{1}{N_j}\sum_{n\in  C_j^c}\gamma_{z,z'}^c(X_{n,j}^{c}(s))+\frac{1}{N_j}\sum_{n\in  C_j^p}\gamma_{z,z'}^p(X_{n,j}^{p}(s))\bigg)\\
                                         &\qquad-\bigg(\frac{1}{N_j}\sum_{n\in  C_j^c}\gamma_{z,z'}^c(Y_{n,j}^{c}(s))+\frac{1}{N_j}\sum_{n\in  C_j^p}\gamma_{z,z'}^p(Y_{n,j}^{p}(s))\bigg)\bigg|\bigg]\\
  &\leq \frac{N_j^c}{N_j}\frac{1}{N_j^c}\sum_{n\in  C_j^c}L_{\gamma}\mathbb{E}\bigg[\bigg|X_{n,j}^{c}(s))-Y_{n,j}^{c}(s))\bigg|\bigg]+\frac{N_j^p}{N_j}\frac{1}{N_j^p}\sum_{n\in  C_j^p}L_{\gamma}\mathbb{E}\bigg[\bigg|X_{n,j}^{p}(s))-Y_{n,j}^{p}(s))\bigg|\bigg]\\
  &\leq \frac{N_j^c}{N_j} L_{\gamma}\max_{m\in C_j^c}\mathbb{E}\big\|X_{m,j}^{c}-Y_{m,j}^{c}\big\|_s+\frac{N_j^p}{N_j} L_{\gamma}\max_{m\in C_j^p}\mathbb{E}\big\|X_{m,j}^{p}-Y_{m,j}^{p}\big\|_s,
  \end{split}
  \label{ineq-4}
  \end{equation}
 where $L_{\gamma}$ is the maximum Lipschitz constant of the functions $\gamma_{z,z'}^c$ and $\gamma_{z,z'}^p$ for all $(z,z')\in\mathcal{Z}^2$. Moreover, by adding and subtracting terms and using the fact that both $\{Y_{n,j}^{c}(s)\}$ and $\{Y_{n,j}^p(s)\}$ are sequences of i.i.d. random variables, the second expectation of the right-hand side of $(\ref{ineq-3})$ can be bounded as follows:
\begin{equation}
\begin{split}
\mathbb{E}&\bigg[\bigg|\bigg(\frac{1}{N_j}\sum_{n\in  C_j^c}\gamma_{z,z'}^c(Y_{n,j}^{c}(s))+\frac{1}{N_j}\sum_{n\in  C_j^p}\gamma_{z,z'}^p(Y_{n,j}^{p}(s))\bigg)-\bigg(p_j^c\int_{\mathcal{Z}}\gamma^c_{z,z'}(x)\mu_j^c(s)ds+p_j^p\int_{\mathcal{Z}}\gamma^p_{z,z'}(x)\mu_j^p(s)ds\bigg)\bigg|\bigg]\\
&\leq \mathbb{E}\bigg[\bigg|\frac{1}{N_j}\sum_{n\in  C_j^c}\gamma_{z,z'}^c(Y_{n,j}^{c}(s))-p_j^c\int_{\mathcal{Z}}\gamma^c_{z,z'}(x)\mu_j^c(s)ds\bigg|\bigg]+\mathbb{E}\bigg[\bigg|\frac{1}{N_j}\sum_{n\in  C_j^p}\gamma_{z,z'}^p(Y_{n,j}^{p}(s))-p_j^p\int_{\mathcal{Z}}\gamma^p_{z,z'}(x)\mu_j^p(s)ds\bigg|\bigg]\\
&\leq \mathbb{E}\bigg[\bigg|p_j^c\frac{1}{N_j^c}\sum_{n\in  C_j^c}\bigg(\gamma_{z,z'}^c(Y_{n,j}^{c}(s))-\mathbb{E}\big[\gamma_{z,z'}^c(Y_{n,j}^{c}(s))\big]\bigg)\bigg|\bigg]+\bigg|\frac{1}{N_j}-\frac{p_j^c}{N_j^c}\bigg|\mathbb{E}\bigg[\bigg|\sum_{n\in C_j^c}\gamma_{z,z'}^c(Y_{n,j}^{c}(s))\bigg|\bigg]\\
&+\mathbb{E}\bigg[\bigg|p_j^p\frac{1}{N_j^p}\sum_{n\in  C_j^p}\bigg(\gamma_{z,z'}^p(Y_{n,j}^{p}(s))-\mathbb{E}\big[\gamma_{z,z'}^p(Y_{n,j}^{c}(s))\big]\bigg)\bigg|\bigg]+\bigg|\frac{1}{N_j}-\frac{p_j^p}{N_j^p}\bigg|\mathbb{E}\bigg[\bigg|\sum_{n\in C_j^p}\gamma_{z,z'}^p(Y_{n,j}^{p}(s))\bigg|\bigg].
\end{split}
\label{ineq-5}
\end{equation}

Note that, by the exchangeability of $\{Y_{n,j}^{c}(s), n\in C_j^c\}$ and the boundedness of the functions $\gamma^c_{z,z'}$, we obtain
\begin{equation}
\bigg|\frac{1}{N_j}-\frac{p_j^c}{N_j^c}\bigg|\mathbb{E}\bigg[\sum_{n\in C_j^c}\gamma_{z,z'}^c(Y_{n,j}^{c}(s))\bigg]=\bigg|\frac{N_j^c}{N_j}-p_j^c\bigg|\mathbb{E}\bigg[\gamma_{z,z'}^c(Y_{1,j}^{c}(s))\bigg]\leq\bar{\gamma} \bigg|\frac{N_j^c}{N_j}-p_j^c\bigg|,
\label{ineq-6}
\end{equation}
which vanishes as $N\rightarrow\infty$ by $(\ref{p-regul})$. In the same manner, the fourth term in the right member of $(\ref{ineq-5})$ is also bounded as follows:
\begin{equation}
\bigg|\frac{1}{N_j}-\frac{p_j^p}{N_j^p}\bigg|\mathbb{E}\bigg[\sum_{n\in C_j^p}\gamma_{z,z'}^p(Y_{n,j}^{p}(s))\bigg]\leq\bar{\gamma} \bigg|\frac{N_j^p}{N_j}-p_j^p\bigg|,
\label{ineq-6-2}
\end{equation}
which also goes to zero as $N\rightarrow\infty$ again by $(\ref{p-regul})$. Furthermore, using $(\ref{moment-bound})$, the first and the third expectations in $(\ref{ineq-5})$ are bounded by $\frac{\kappa_1 p_j^c}{\sqrt{N_j^c}}$ and $\frac{\kappa_2 p_j^p}{\sqrt{N_j^p}}$, respectively, where $\kappa_1$ and $\kappa_2$ are positives constants.
  \\
\indent Now, take a look at the first term of the right-hand side of $(\ref{ineq-2})$. Since $X^c_{n,j}$ and $Y^c_{n,j}$ are $\mathcal{Z}$-valued, which is a subset of $\mathbb{N}$, one can easily find that
 \begin{equation}
 \begin{split}
 \mathbb{E}\bigg[\lambda^c_{z,z'}(\upsilon_j^{N}(s))\bigg|(\mathds{1}_{X^{c}_{n,j}(s)=z}-\mathds{1}_{Y^{c}_{n,j}(s)=z})\bigg|\bigg]&\leq  \mathbb{E}\bigg[\lambda^c_{z,z'}(\upsilon_j^{N}(s))\big|X^{c}_{n,j}(s)-Y^{c}_{n,j}(s)\big|\bigg]\\
 &\leq\bar{\gamma}\mathbb{E}\big|X^{c}_{n,j}(s)-Y^{c}_{n,j}(s)\big|\\
  &\leq\bar{\gamma}\mathbb{E}\big\|X^{c}_{n,j}-Y^{c}_{n,j}\big\|_s,
 \end{split}
 \label{ineq-7}
 \end{equation}
 where $\bar{\gamma}$ is the upper bound of the function  $\gamma_{z,z'}^c$ and $\gamma_{z,z'}^p$ for all $(z,z')\in\mathcal{Z}^2$. Finally, by combining $(\ref{ineq-2}), (\ref{ineq-3}), (\ref{ineq-4}), (\ref{ineq-5}), (\ref{ineq-6})$, $(\ref{ineq-6-2})$ and $(\ref{ineq-7})$ we obtain
\begin{equation}
\begin{split}
\mathbb{E}\left[\|X_{n,j}^{c}-Y_{n,j}^{c}(t) \|_t\right]\leq K |\mathcal{E}|\int_{0}^t \bigg[&\bar{\gamma}\mathbb{E}\big\|X^{c}_{n,j}-Y^{c}_{n,j}\big\|_s+\frac{N_j^c}{N_j} L_{\gamma}\max_{m\in C_j^c}\mathbb{E}\big\|X_{m,j}^{c}-Y_{m,j}^{c}\big\|_s \\
                                                   &+\frac{N_j^p}{N_j} L_{\gamma}\max_{m\in C_j^p}\mathbb{E}\big\|X_{m,j}^{p}-Y_{m,j}^{p}\big\|_s+\bar{\gamma} \bigg|\frac{N_j^c}{N_j}-p_j^c\bigg|+\bar{\gamma} \bigg|\frac{N_j^p}{N_j}-p_j^p\bigg|\\
                                                   &+\frac{\kappa_1 p_j^c}{\sqrt{N_j^c}}+\frac{\kappa_2 p_j^p}{\sqrt{N_j^p}}\bigg]ds,
\end{split}
\label{ineq-fin-cent}
\end{equation}
 where $|\mathcal{E}|$ stands for the cardinal of the set of edges $\mathcal{E}$ contained in the graph $(\mathcal{Z},\mathcal{E})$. Taking the maximum over $n\in C_j^c$ and then over $1\leq j\leq r$ in $(\ref{ineq-fin-cent})$ we obtain
\begin{equation}
\begin{split}
\max_{1\leq j\leq r}\max_{n\in C_j^c}\mathbb{E}\left[\|X_{n,j}^{c}-Y_{n,j}^{c}\|_t\right]\leq K |\mathcal{E}|&\int_{0}^t \bigg[ (\bar{\gamma}+\frac{N_{j^*}^c}{N_{j^*}} L_{\gamma})\max_{1\leq j\leq r}\max_{n\in C_{j}^c}\mathbb{E}\big\|X_{n,j}^{c}-Y_{n,j}^{c}\big\|_s \\
                                                   &+\frac{N_{j^*}^p}{N_{j^*}} L_{\gamma}\max_{1\leq j\leq r}\max_{n\in C_{j}^p}\mathbb{E}\big\|X_{n,j}^{p}-Y_{n,j}^{p}\big\|_s+\bar{\gamma} \bigg|\frac{N_{j^*}^c}{N_{j^*}}-p_{j^*}^c\bigg|+\bar{\gamma} \bigg|\frac{N_{j^*}^p}{N_{j^*}}-p_{j^*}^p\bigg|\\
                                                   &+\frac{\kappa_1 p_{j^*}^c}{\sqrt{N_{j^*}^c}}+\frac{\kappa_2 p_{j^*}^p}{\sqrt{N_{j^*}^p}}\bigg]ds,
\end{split}
\label{ineq-fin-cent-max}
\end{equation}
 where $(n^*,j^*)=\argmax\limits_{\substack{1\leq j\leq r, n\in C_j^c}}\mathbb{E}\|X_{n,j}^{c}-Y_{n,j}^{c} \|_t$.


\subparagraph*{Step 2.}  Fix a block $1\leq j\leq r$ and a peripheral node $n\in C_j^p$. For any $t\in [0,T]$, we have
\begin{equation}
\begin{split}
\mathbb{E}\left[\|X_{n,j}^{p}-Y_{n,j}^p \|_t\right]&=\mathbb{E}\left[\sup_{0\leq s\leq t}\big|X_{n,j}^{p}(s)-Y_{n,j}^p(s)\big|\right]\\
                                             &\leq\mathbb{E}\bigg[\bigg|\int_{[0,t]\times\mathbb{R}_+}\sum_{(z,z')\in\mathcal{E}}\mathds{1}_{X^{p}_{n,j}(s)=z}(z'-z)\mathds{1}_{\left[0,\lambda^p_{z,z'}(\upsilon_{n,j}^N(s))\right]}(y)\mathcal{N}_{n,j}^p(ds,dy)\\
                                             &\qquad\qquad -\int_{[0,t]\times\mathbb{R}_+}\sum_{(z,z')\in\mathcal{E}}\mathds{1}_{Y^{p}_{n,j}(s)=z}(z'-z)\mathds{1}_{\left[0,\lambda^p_{z,z'}(\upsilon_{n,j}(s))\right]}(y)\mathcal{N}_{n,j}^p(ds,dy) \bigg| \bigg]\\
                                             &\leq  K\int_{[0,t]}\sum_{(z,z')\in\mathcal{E}}\mathbb{E}\bigg[\bigg|\bigg(\mathds{1}_{X^{p}_{n,j}(s)=z}-\mathds{1}_{Y^{p}_{n,j}(s)=z}\bigg)\lambda^p_{z,z'}\left(\upsilon^N_{n,j}(s)\right) \\
                                             &\qquad\qquad+\bigg(\lambda^p_{z,z'}\left(\upsilon^N_{n,j}(s)\right)-\lambda^p_{z,z'}\left(\upsilon_{n,j}(s)\right)\bigg) \bigg|\bigg]ds,
\end{split}
\label{ineq-p1}
\end{equation}
where the last inequality is obtained by following the same steps as in $(\ref{ineq-mart1})$ and $(\ref{ineq-2})$. Again, given that $X_{n,j}^p$ and $Y_{n,j}^p$ are $\mathcal{Z}$-valued and that $\mathcal{Z}\subset\mathbb{N}$, the first expectation in the right-hand side of $(\ref{ineq-p1})$ can be bounded as follows:
\begin{equation}
\begin{split}
\mathbb{E}\bigg[\bigg|\bigg(\mathds{1}_{X^{p}_{n,j}(s)=z}-\mathds{1}_{Y^{p}_{n,j}(s)=z}\bigg)\lambda^p_{z,z'}\left(\upsilon^N_{n,j}(s)\right)\bigg|\bigg]\leq \bar{\gamma}\mathbb{E}\bigg[\big\|X^{p}_{n,j}-Y^{p}_{n,j}\big\|_s\bigg].
\end{split}
\label{ineq-p2}
\end{equation}

It remains to bound the second term in the right-hand side of $(\ref{ineq-p1})$. Under Condition $\ref{Cond-prin}$ we get
\begin{equation}
\begin{split}
\mathbb{E}\bigg[\bigg|\bigg(\lambda^p_{z,z'}\left(\upsilon^N_{n,j}(s)\right)-\lambda^p_{z,z'}&\big(\upsilon_{n,j}(s)\big)\bigg) \bigg|\bigg]=\mathbb{E}\bigg[\bigg|\frac{1}{deg(n)+1}\bigg(\sum_{n\in C_j^c}\gamma^c_{z,z'}(X_{n,j}^{c}(s))+\sum_{\substack{m\in C_1^p\\(m,n)\in\Xi}}\gamma_{z,z'}^p(X_{m,1}^{p}(s))+\\
&\ldots+\sum_{\substack{n\in C_j^p}}\gamma_{z,z'}^p(X_{n,j}^{p}(s))+\ldots+\sum_{\substack{m\in C_r^p\\(m,n)\in\Xi}}\gamma_{z,z'}^p(X_{m,r}^{p}(s))\bigg)\\
&-\bigg(\alpha_j^c\int_{\mathcal{Z}}\gamma_{z,z'}^c(x)\mu_j^c(s)ds+q_{j1}\int_{\mathcal{Z}}\gamma_{z,z'}^p(x)\mu_1^p(s)ds+\cdots+q_{jr}\int_{\mathcal{Z}}\gamma_{z,z'}^p(x)\mu_r^p(s)ds\bigg)\bigg|\bigg].
\end{split}
\end{equation}
By rearranging the terms and using the triangle inequality we obtain
\begin{equation}
\begin{split}
 \mathbb{E}\bigg[\bigg|\bigg(\lambda^p_{z,z'}\left(\upsilon^N_{n,j}(s)\right)-\lambda^p_{z,z'}\left(\upsilon_{n,j}(s)\right)\bigg) \bigg|\bigg]  &\leq\mathbb{E}\bigg[\bigg| \frac{1}{deg(n)+1}\sum_{n\in C_j^c}\gamma^c_{z,z'}(X_{n,j}^{c}(s))-\alpha_j^c\int_{\mathcal{Z}}\gamma_{z,z'}^c(x)\mu_j^c(s)ds  \bigg|\bigg]\\
  &+\mathbb{E}\bigg[\bigg|\frac{1}{deg(n)+1}\sum_{\substack{m\in C_1^p\\(m,n)\in\Xi}}\gamma_{z,z'}^p(X_{m,1}^{p}(s))-q_{j1}\int_{\mathcal{Z}}\gamma_{z,z'}^p(x)\mu_1^p(s)ds\\
                                             &\qquad\vdots\\
                                             &\qquad+\frac{1}{deg(n)+1}\sum_{\substack{n\in C_j^p}}\gamma_{z,z'}^p(X_{n,j}^{p}(s))-q_{jj}\int_{\mathcal{Z}}\gamma_{z,z'}^p(x)\mu_j^p(s)ds\\
                                             &\qquad\vdots\\
                                             &\qquad+\frac{1}{deg(n)+1}\sum_{\substack{m\in C_r^p\\(m,n)\in\Xi}}\gamma_{z,z'}^p(X_{m,r}^{p}(s))-q_{jr}\int_{\mathcal{Z}}\gamma_{z,z'}^p(x)\mu_r^p(s)ds \bigg|\bigg].
\end{split}
\end{equation}

Note that $\int_{\mathcal{Z}}\gamma_{z,z'}^p(x)\mu_i^p(s)ds= \mathbb{E}[\gamma^p_{z,z'}(Y_{m,i}^{p}(s))]$ for $m\in C_i^p$ and $\int_{\mathcal{Z}}\gamma_{z,z'}^c(x)\mu_j^c(s)ds=\mathbb{E}[\gamma^c_{z,z'}(Y_{n,j}^{c}(s))]$ for $n\in C_j^c$. Then, by using the exchangeability of $Y_{n,j}^{c}(t)$ for $n\in C_j^c$ we obtain
\begin{equation}
\begin{split}
 \mathbb{E}\bigg[\bigg|\bigg(\lambda^p_{z,z'}\left(\upsilon^N_{n,j}(s)\right)-\lambda^p_{z,z'}\left(\upsilon_{n,j}(s)\right)\bigg) \bigg|\bigg]  &\leq\mathbb{E}\bigg[\bigg| \frac{1}{deg(n)+1}\sum_{n\in C_j^c}\gamma^c_{z,z'}(X_{n,j}^{c}(s))-\alpha_j^c\mathbb{E}[\gamma^c_{z,z'}(Y_{n,j}^{c}(s))]  \bigg|\bigg] \\
                                             &+\mathbb{E}\bigg[\bigg|\frac{1}{deg(n)+1}\sum_{\substack{m\in C_1^p\\ (m,n)\in\Xi}}\gamma_{z,z'}^p(X_{m,1}^{p}(s))-q_{j1}\mathbb{E}[\gamma^p_{z,z'}(Y_{m,1}^{p}(s))] \\
                                             &\quad\vdots\\
                                             &\quad+\frac{1}{deg(n)+1}\sum_{\substack{n\in C_j^p}}\gamma_{z,z'}^p(X_{n,j}^{p}(s))-q_{jj}\mathbb{E}[\gamma^p_{z,z'}(Y_{m,j}^{p}(s))] \\
                                             &\quad\vdots\\
                                             &\quad+\frac{1}{deg(n)+1}\sum_{\substack{m\in C_r^p\\(m,n)\in\Xi}}\gamma_{z,z'}^p(X_{m,r}^{p}(s))-q_{jr}\mathbb{E}[\gamma^p_{z,z'}(Y_{m,r}^{p}(s))]\bigg|\bigg].
\end{split}
\label{ineq-p3}
\end{equation}

 Observe that there are $r+1$ terms  on the right-hand side of the last inequality. Let us start with the first expectation. By adding and subtracting terms we obtain
   \begin{equation}
   \begin{split}
&\mathbb{E}\bigg[\bigg| \frac{1}{deg(n)+1}\sum_{n\in C_j^c}\gamma^c_{z,z'}(X_{n,j}^{c}(s))-\alpha_j^c\mathbb{E}[\gamma^c_{z,z'}(Y_{n,j}^{c}(s))]  \bigg|\bigg]\\
 &\leq \mathbb{E}\bigg[\bigg| \frac{1}{deg(n)+1}\sum_{n\in C_j^c}\bigg(\gamma^c_{z,z'}(X_{n,j}^{c}(s))-\gamma^c_{z,z'}(Y_{n,j}^{c}(s))\bigg)\bigg|\bigg]\\
&\qquad+\mathbb{E}\bigg[\bigg|\frac{1}{deg(n)+1}\sum_{n\in C_j^c}\gamma^c_{z,z'}(Y_{n,j}^{c}(s))-\alpha_j^c\mathbb{E}[\gamma^c_{z,z'}(Y_{n,j}^{c}(s))]  \bigg|\bigg].
   \end{split}
 \label{ineq-p4}
   \end{equation}

Note that, by the Lipschitz property of the functions $\gamma^c_{z,z'}$
\begin{equation}
\begin{split}
\mathbb{E}\bigg[\bigg| \frac{1}{deg(n)+1}\sum_{n\in C_j^c}\bigg(\gamma^c_{z,z'}(X_{n,j}^{c}(s))-\gamma^c_{z,z'}(Y_{n,j}^{c}(s))\bigg)\bigg|\bigg]\leq \frac{N_j^c}{deg(n)+1}L_{\gamma}\max_{m\in C_j^c}\mathbb{E}\left\|X_{m,j}^{c}-Y_{m,j}^{c}\right\|_s                .
\end{split}
 \label{ineq-p4-1}
\end{equation}

Moreover, using $(\ref{moment-bound})$ together with the exchangeability of $\{Y_{n,j}^{c}(s),n\in C_j^c\}$ leads to
\begin{equation}
\begin{split}
&\mathbb{E}\bigg[\bigg|\frac{1}{deg(n)+1}\sum_{n\in C_j^c}\gamma^c_{z,z'}(Y_{n,j}^{c}(s))-\alpha_j^c\mathbb{E}[\gamma^c_{z,z'}(Y_{n,j}^{c}(s))]  \bigg|\bigg]\\
&\qquad\leq \mathbb{E}\bigg[\bigg|\frac{\alpha_j^c}{N_j^c}\sum_{n\in C_j^c}\bigg(\gamma^c_{z,z'}(Y_{n,j}^{c}(s))-\mathbb{E}[\gamma^c_{z,z'}(Y_{n,j}^{c}(s))]\bigg)  \bigg|\bigg]+\mathbb{E}\bigg[\bigg|\left(\frac{1}{deg(n)+1}-\frac{\alpha_j^c}{N_j^c}\right)\sum_{n\in C_j^c}\gamma^c_{z,z'}(Y_{n,j}^{c}(s))\bigg|\bigg]\\
&\qquad\leq \frac{\alpha_j^c\kappa_3}{\sqrt{N_j^c}}+\left|\frac{N_j^c}{deg(n)+1}-\alpha_j^c\right|\bar{\gamma}.
\end{split}
 \label{ineq-p4-2}
\end{equation}

Now, let us examine the remaining $r$ terms on the right-hand side of $(\ref{ineq-p3})$. For simplicity of notation, denote the left side of $(\ref{ineq-p3})$ by $\mathcal{I}$. Therefore,  by adding and subtracting terms we find
\begin{equation}
\begin{split}
&\mathcal{I}\leq \mathbb{E}\bigg[\bigg|\frac{1}{deg(n)+1}\sum_{\substack{m\in C_1^p\\(m,n)\in\Xi}}\bigg(\gamma_{z,z'}^p(X_{m,1}^{p}(s))-\gamma_{z,z'}^p(Y_{m,1}^{p}(s))\bigg)+\cdots\\
&\quad\qquad\cdots+\frac{1}{deg(n)+1}\sum_{\substack{n\in C_j^p}}\bigg(\gamma_{z,z'}^p(X_{n,j}^{p}(s))-\gamma_{z,z'}^p(Y_{n,j}^{p}(s))\bigg)+\cdots\\
&\quad\qquad\cdots+\frac{1}{deg(n)+1}\sum_{\substack{m\in C_r^p\\(m,n)\in\Xi}}\bigg(\gamma_{z,z'}^p(X_{m,r}^{p}(s))-\gamma_{z,z'}^p(Y_{m,r}^{p}(s))\bigg)\bigg|\bigg]\\
&\quad+\mathbb{E}\bigg[\bigg|\frac{1}{deg(n)+1}\sum_{\substack{m\in C_1^p\\(m,n)\in\Xi}}\gamma_{z,z'}^p(Y_{m,1}^{p}(s))-q_{j1} \mathbb{E}[\gamma^p_{z,z'}(Y_{m,1}^{p}(s))]+\cdots\\
&\quad\qquad\cdots+\frac{1}{deg(n)+1}\sum_{\substack{n\in C_j^p}}\gamma_{z,z'}^p(Y_{n,j}^{p}(s))-q_{jj} \mathbb{E}[\gamma^p_{z,z'}(Y_{n,j}^{p}(s))]+\cdots\\
&\quad\qquad\cdots+\frac{1}{deg(n)+1}\sum_{\substack{m\in C_r^p\\(m,n)\in\Xi}}\gamma_{z,z'}^p(Y_{m,r}^{p}(s))-q_{jr} \mathbb{E}[\gamma^p_{z,z'}(Y_{m,r}^{p}(s))]  \bigg|\bigg]\\
&=\mathcal{I}_1+\mathcal{I}_2.
\end{split}
\label{ineq-p5}
\end{equation}

Using the Lipschitz property of the functions $\gamma^p_{z,z'}$, and recalling that, for $i\neq j$, $M_i^n$ represents the number of peripheral  nodes of the $i-$th block connecting node $n$ (in particular $M_j^n+1=N_j^p$), we straightforwardly bound the first expectation $\mathcal{I}_1$ as follows:
\begin{equation}
\begin{split}
\mathcal{I}_1\leq L_{\gamma}\bigg(&\frac{M^n_1}{deg(n)+1}\max_{m\in C_1^p}\mathbb{E}\left\|X_{m,1}^{p}-Y_{m,1}^{p}\right\|_s+\cdots+\frac{N^p_j}{deg(n)+1}\max_{m\in C_j^p}\mathbb{E}\left\|X_{m,j}^{p}-Y_{m,j}^{p}\right\|_s+\cdots\\
&\qquad\qquad\cdots+\frac{M^n_r}{deg(n)+1}\max_{m\in C_r^p}\mathbb{E}\left\|X_{m,r}^{p}-Y_{m,r}^{p}\right\|_s\bigg).
\label{ineq-p6}
\end{split}
\end{equation}

Moreover, by adding and subtracting terms in $\mathcal{I}_2$ we obtain
\begin{equation}
\begin{split}
\mathcal{I}_2 \leq \mathcal{I}_3+\mathcal{I}_4,
\end{split}
\end{equation}
where
\begin{equation}
\begin{split}
\mathcal{I}_3&=\mathbb{E}\bigg[\bigg|\frac{q_{j1}}{M^n_1}\sum_{\substack{m\in C_1^p\\(m,n)\in\Xi}}\bigg(\gamma_{z,z'}^p(Y_{m,1}^{p}(s))-\mathbb{E}\left[\gamma^p_{z,z'}(Y_{m,1}^{p}(s))\right]\bigg)\\
&\qquad\vdots\\
&\qquad+\frac{q_{jj}}{N^p_j}\sum_{\substack{m\in C_j^p}}\bigg(\gamma_{z,z'}^p(Y_{m,r}^{p}(s))-\mathbb{E}\left[\gamma^p_{z,z'}(Y_{m,r}^{p}(s))\right]\bigg)\\
&\qquad\vdots\\
&\qquad+\frac{q_{jr}}{M^n_r}\sum_{\substack{m\in C_r^p\\ (m,n)\in\Xi}}\bigg(\gamma_{z,z'}^p(Y_{m,r}^{p}(s))-\mathbb{E}\left[\gamma^p_{z,z'}(Y_{m,r}^{p}(s))\right]\bigg)\bigg|\bigg],
\end{split}
\label{I_3}
\end{equation}
and
\begin{equation}
\begin{split}
\mathcal{I}_4=\mathbb{E}\bigg[\bigg|\left(\frac{1}{deg(n)+1}-\frac{q_{j1}}{M^n_1}\right)&\sum_{\substack{m\in C_1^p\\(m,n)\in\Xi}}\gamma_{z,z'}^p(Y_{m,1}^{p}(s))+\cdots+\left(\frac{1}{deg(n)+1}-\frac{q_{jj}}{N^p_j}\right)\sum_{\substack{m\in C_j^p}}\gamma_{z,z'}^p(Y_{m,j}^{p}(s))+\cdots\\
&\qquad\qquad\cdots+\left(\frac{1}{deg(n)+1}-\frac{q_{jr}}{M^n_r}\right)\sum_{\substack{m\in C_r^p\\(m,n)\in\Xi}}\gamma_{z,z'}^p(Y_{m,r}^{p}(s))\bigg|\bigg].
\end{split}
\label{I_4}
\end{equation}

We first prove that under Condition $\ref{Cond-prin}$, the term $\mathcal{I}_4$ goes to zero. Indeed, the triangle inequality gives to us
\begin{equation}
\begin{split}
\mathcal{I}_4\leq \mathbb{E}&\bigg[\bigg|\left(\frac{1}{deg(n)+1}-\frac{q_{j1}}{M^n_1}\right)\sum_{\substack{m\in C_1^p\\(m,n)\in\Xi}}\gamma_{z,z'}^p(Y_{m,1}^{p}(s))\bigg|\bigg]+\cdots+\mathbb{E}\bigg[\bigg|\left(\frac{1}{deg(n)+1}-\frac{q_{jj}}{N^p_j}\right)\sum_{\substack{m\in C_j^p}}\gamma_{z,z'}^p(Y_{m,j}^{p}(s))\bigg|\bigg]+\\
&\cdots+\mathbb{E}\bigg[\bigg|\left(\frac{1}{deg(n)+1}-\frac{q_{jr}}{M^n_r}\right)\sum_{\substack{m\in C_r^p\\(m,n)\in\Xi}}\gamma_{z,z'}^p(Y_{m,r}^{p}(s))\bigg|\bigg].
\end{split}
\label{I_5}
\end{equation}

Using $(\ref{cond-regul})$, the exchangeability of the variables $\{Y^p_{n,j}, n\in C_j^p\}$ and the boundedness of the functions $\gamma_{z,z'}^p$, we easily show that the right-hand side of $(\ref{I_5})$ goes to zero. Indeed, the $j$-th term satisfies
\begin{equation}
\begin{split}
\bigg|\frac{1}{deg(n)+1}-\frac{q_{jj}}{N^p_j}\bigg|\mathbb{E}\bigg|\sum_{\substack{m\in C_j^p}}\gamma_{z,z'}^p(Y_{m,j}^{p}(s))\bigg|&=\bigg|\frac{N^p_j}{deg(n)+1}-q_{jj}\bigg|\frac{1}{N_j^p}\mathbb{E}\bigg|\sum_{\substack{m\in C_j^p}}\gamma_{z,z'}^p(Y_{m,j}^{p}(s))\bigg|\\
                                                       &\leq\bigg|\frac{N^p_j}{deg(n)+1}-q_{jj}\bigg|\bar{\gamma},
\end{split}
\label{I_5-1}
\end{equation}
 and thus goes to zero by $(\ref{cond-regul})$. Using the same steps, we obtain for $1\leq i\leq r$ with $i\neq j$ that
\begin{equation}
\begin{split}
\bigg|\frac{1}{deg(n)+1}-\frac{q_{ji}}{M_i^n}\bigg|\mathbb{E}\bigg|\sum_{\substack{m\in C_i^p\\(m,n)\in\Xi}}\gamma_{z,z'}^p(Y_{m,i}^{p}(s))\bigg|&\leq\bigg|\frac{M_i^n}{deg(n)+1}-q_{ji}\bigg|\bar{\gamma},
\end{split}
\label{I_5-2}
\end{equation}
which also  vanishes as $N\rightarrow\infty$ by $(\ref{cond-regul})$, so does $\mathcal{I}_4$. In order to bound $\mathcal{I}_3$, we use again the moment inequality $(\ref{moment-bound})$ which straightforwardly gives to us
\begin{align}
I_3\leq \frac{\theta_1 q_{j1}}{\sqrt{M_1^n}}+\cdots+\frac{\theta_j q_{jj}}{\sqrt{N_j^p}}+\cdots+\frac{\theta_r q_{jr}}{\sqrt{M_r^n}},
\label{I_3-final}
\end{align}
where $\theta_1,\cdots,\theta_r$ are positive constants.

Now, by $(\ref{ineq-p2})$, $(\ref{ineq-p4-1})$, $(\ref{ineq-p4-2})$, $(\ref{ineq-p6})$, $(\ref{I_5-1})$, $(\ref{I_5-2})$ and  $(\ref{I_3-final})$ we obtain
\begin{equation}
\begin{split}
\mathbb{E}\left[\|X_{n,j}^{p}-Y_{n,j}^p \|_t\right]\leq &K |\mathcal{E}|\int_0^t\bigg[\bar{\gamma}\mathbb{E}\big\|X^{p}_{n,j}-Y^{p}_{n,j}\big\|_s+\frac{N_j^c}{deg(n)+1}L_{\gamma}\max_{m\in C_j^c}\mathbb{E}\left\|X_{m,j}^{c}-Y_{m,j}^{c}\right\|_s   +\frac{\alpha_j^c\kappa_3}{\sqrt{N_j^c}}+\\
&+\left|\frac{N_j^c}{deg(n)+1}-\alpha_j^c\right|\bar{\gamma}+L_{\gamma}\bigg(\frac{M^n_1}{deg(n)+1}\max_{m\in C_1^p}\mathbb{E}\left\|X_{m,1}^{p}-Y_{m,1}^{p}\right\|_s+\cdots\\
&\qquad+\frac{N^p_j}{deg(n)+1}\max_{m\in C_j^p}\mathbb{E}\left\|X_{m,j}^{p}-Y_{m,j}^{p}\right\|_s+\cdots+\frac{M^n_r}{deg(n)+1}\max_{m\in C_r^p}\mathbb{E}\left\|X_{m,r}^{p}-Y_{m,r}^{p}\right\|_s\bigg)\\
&+\bigg|\frac{M_1^n}{deg(n)+1}-q_{j1}\bigg|\bar{\gamma}+\cdots+\bigg|\frac{N^p_j}{deg(n)+1}-q_{jj}\bigg|\bar{\gamma}+\cdots+\bigg|\frac{M_r^n}{deg(n)+1}-q_{jr}\bigg|\bar{\gamma}\\
&+\frac{\theta_1 q_{j1}}{\sqrt{M_1^n}}+\cdots+\frac{\theta_j q_{jj}}{\sqrt{N_j^p}}+\cdots+\frac{\theta_r q_{jr}}{\sqrt{M_r^n}}\bigg]ds.
\end{split}
\label{final-ineq-periph}
\end{equation}

Taking the maximum over $n\in C_j^p$  and $1\leq j\leq r$ we get,
\begin{equation}
\begin{split}
\max_{1\leq j\leq r}\max_{n\in C_j^p}\mathbb{E}\left[\|X_{n,j}^{p}-Y_{n,j}^p \|_t\right]\leq K |\mathcal{E}|\int_0^t\bigg[&\frac{N_{j^{\diamond}}^c}{deg(n^{\diamond})+1}L_{\gamma}\max_{1\leq j\leq r}\max_{m\in C_j^c}\mathbb{E}\left\|X_{m,j}^{c}-Y_{m,j}^{c}\right\|_s+\\
&+\left(\bar{\gamma}+L_{\gamma}\right)\max_{1\leq j\leq r}\max_{n\in C_j^p}\mathbb{E}\big\|X^{p}_{n,j}-Y^{p}_{n,j}\big\|_s+\\
&+\frac{\alpha_{j^{\diamond}}^c\kappa_3}{\sqrt{N_{j^{\diamond}}^c}}+\left|\frac{N_{j^{\diamond}}^c}{deg(n^{\diamond})+1}-\alpha_{j^{\diamond}}^c\right|\bar{\gamma}+\bigg|\frac{M_1^{n^{\diamond}}}{deg(n^{\diamond})+1}-q_{j^{\diamond}1}\bigg|\bar{\gamma}+\cdots\\
&+\cdots+\bigg|\frac{N^p_{j^{\diamond}}}{deg(n^{\diamond})+1}-q_{j^{\diamond}j^{\diamond}}\bigg|\bar{\gamma}+\cdots+\bigg|\frac{M_r^{n^{\diamond}}}{deg(n^{\diamond})+1}-q_{j^{\diamond}r}\bigg|\bar{\gamma}\\
&+\frac{\theta_1 q_{j^{\diamond}1}}{\sqrt{M_1^{n^{\diamond}}}}+\cdots+\frac{\theta_{j^{\diamond}} q_{j^{\diamond}j^{\diamond}}}{\sqrt{N_{j^{\diamond}}^p}}+\cdots+\frac{\theta_r q_{j^{\diamond}r}}{\sqrt{M_r^{n^{\diamond}}}}\bigg]ds,\\
\end{split}
\label{fin-ineq-peri-max}
\end{equation}
where $(n^{\diamond},j^{\diamond})=\argmax\limits_{\substack{1\leq j\leq r, n\in C_j^p}}\mathbb{E}\|X_{n,j}^{p}-Y_{n,j}^p \|_t$.\\
 Adding side by side the two inequalities $(\ref{ineq-fin-cent-max})$ and $(\ref{fin-ineq-peri-max})$ gives
\begin{equation}
\begin{split}
\max_{1\leq j\leq r}\max_{n\in C_j^c}\mathbb{E}\|X_{n,j}^{c}-Y_{n,j}^{c}(t) \|_t+\max_{1\leq j\leq r}\max_{n\in C_j^p}\mathbb{E}\|X_{n,j}^{p}-Y_{n,j}^p \|_t\leq &K |\mathcal{E}|\int_0^t\bigg(C_1(N)\max_{1\leq j\leq r}\max_{n\in C_j^c}\mathbb{E}\|X_{n,j}^{c}-Y_{n,j}^{c}(t) \|_s \\
&+C_2(N)\max_{1\leq j\leq r}\max_{n\in C_j^p}\mathbb{E}\|X_{n,j}^{p}-Y_{n,j}^p \|_s+C_3(N)\bigg),
\end{split}
\label{final-ineq}
\end{equation}
where, with a slight abuse of notations, the functions $C_1(N),C_2(N)$ and $C_3(N)$ are defined by
\begin{align*}
C_1(N)&= \bar{\gamma}+L_{\gamma}\left(\frac{N_{j^*}^c}{N_{j^*}}+\frac{N^c_{j^{\diamond}}}{deg(n^{\diamond})+1}\right),\\
C_2(N)&=\bar{\gamma}+\frac{N_{j^*}^p}{N_{j^*}}L_{\gamma}+L_{\gamma},\\
C_3(N)&= \bar{\gamma} \bigg|\frac{N_{j^*}^c}{N_{j^*}}-p_{j^*}^c\bigg|+\bar{\gamma} \bigg|\frac{N_{j^*}^p}{N_{j^*}}-p_{j^*}^p\bigg|+\frac{\kappa_1 p_{j^*}^c}{\sqrt{N_{j^*}^c}}+\frac{\kappa_2 p_{j^*}^p}{\sqrt{N_{j^*}^p}}+
\frac{\alpha_{j^{\diamond}}^c\kappa_3}{\sqrt{N_{j^{\diamond}}^c}}+\left|\frac{N_{j^{\diamond}}^c}{deg(n^{\diamond})+1}-\alpha_{j^{\diamond}}^c\right|\bar{\gamma},\\
&\qquad+\bigg|\frac{M_1^{n^{\diamond}}}{deg(n^{\diamond})+1}-q_{j^{\diamond}1}\bigg|\bar{\gamma}+\cdots+\bigg|\frac{N^p_{j^{\diamond}}}{deg(n^{\diamond})+1}-q_{j^{\diamond}j^{\diamond}}\bigg|\bar{\gamma}+\cdots+\bigg|\frac{M_r^{n^{\diamond}}}{deg(n^{\diamond})+1}-q_{j^{\diamond}r}\bigg|\bar{\gamma},\\
&\qquad+\frac{\theta_1 q_{j^{\diamond}1}}{\sqrt{M_1^{n^{\diamond}}}}+\cdots+\frac{\theta_{j^{\diamond}} q_{j^{\diamond}j^{\diamond}}}{\sqrt{N_{j^{\diamond}}^p}}+\cdots+\frac{\theta_r q_{j^{\diamond}r}}{\sqrt{M_r^{n^{\diamond}}}}.
\end{align*}

Therefore, applying Gr\"{o}nwall's lemma to $(\ref{final-ineq})$ we obtain
\begin{equation}
\begin{split}
\max_{1\leq j\leq r}\max_{n\in C_j^c}\mathbb{E}\|X_{n,j}^{c}-Y_{n,j}^{c}(t) \|_t+\max_{1\leq j\leq r}\max_{n\in C_j^p}\mathbb{E}\|X_{n,j}^{p}-Y_{n,j}^p \|_t\leq &K |\mathcal{E}|C_3(N)t \exp\bigg\{ \int_0^t C_4(N) ds\bigg\},
\end{split}
\label{gron-2}
\end{equation}
with
\begin{align*}
C_4(N)=C_{1}(N)+C_{2}(N)= 2\bar{\gamma}+2 L_{\gamma}+L_{\gamma}\frac{N^c_{j^{\diamond}}}{deg(n^{\diamond})+1}.
\end{align*}

Thus, under Condition $\ref{Cond-prin}$, we can easily see that  $C_4(N)\rightarrow 2\bar{\gamma}+2 L_{\gamma}+L_{\gamma}\alpha^c_{j^{\diamond}}$ and $C_3(N)\rightarrow 0$, which proves $(\ref{Asymp-close})$.

We are now ready to conclude the proof of the theorem. The relation in $(\ref{cont-init})$ shows that the solution of the limiting SDE $(\ref{limit-syst})$ is continuous with respect to the initial condition. In addition, Theorem $\ref{theo1}$ ensures the uniqueness of the solution. Hence, the process
$Y(t)=\big((Y_n^{c}(t),Y_m^{p}(t)),n\in C_j^c,m\in C_j^p;1\leq j\leq r\big)$ is $P_{\bar{X}}=\mu^c_{1}\otimes\mu_{1}^p\otimes\cdots \mu_{r}^c\otimes\mu_{r}^p$-multi-chaotic. Therefore, by the relation in $(\ref{Asymp-close})$, we conclude that the sequence of processes $\big((X_n^{c,N}(t),X_m^{p,N}(t),t\geq 0),n\in C_j^c,m\in C_j^p;1\leq j\leq r\big)$ is also  $P_{\bar{X}}=\mu^c_{1}\otimes\mu_{1}^p\otimes\cdots \mu_{r}^c\otimes\mu_{r}^p$-multi-chaotic. The theorem is proved.
\carre

The following laws of large numbers are  immediate consequences of Theorem $\ref{theo2}$.

\begin{coro}
Suppose that the conditions of Theorem $\ref{theo2}$ hold true. Denote $\mu_j^{c}=\mathcal{L}( \bar{X}^c_{n,j}),\mu_j^{p}=\mathcal{L}( \bar{X}^p_{n,j})$ for $1\leq j\leq r$, where $\big((\bar{X}^c_{n,j}(t), \bar{X}^p_{n,j}(t),t\geq 0),1\leq j\leq r\big)$ is the solution of the McKean-Vlasov limiting system in $(\ref{limit-syst})$ with initial distribution $\nu^{1,c}\otimes\nu^{1,p}\cdots \nu^{r,c}\otimes\nu^{r,p}$. Then, for each $1\leq j\leq r$, as $N\rightarrow\infty$,

\begin{align}
\mu_j^{c,N}=\frac{1}{N_j^c}\sum_{n\in C_j^c}\delta_{X_{n,j}^c}\rightarrow\mu_j^c\quad\text{in}\quad\mathcal{M}_1(\mathcal{D}([0,T],\mathcal{Z}))\quad\text{in probability},
\label{conv-meas-1}
\end{align}
and
\begin{align}
\mu_{j}^{p,N}=\frac{1}{N_j^p}\sum_{n\in C_j^p}\delta_{X_{n,j}^p}\rightarrow\mu_{j}^p,
\label{conv-meas-2}
\end{align}
for the weak topology on $\mathcal{M}_1(\mathcal{D}([0,T],\mathcal{Z}))$ with $\mathcal{D}([0,T],\mathcal{Z})$ endowed with the Skorohod topology.
\label{coro-1}
\end{coro}

\paragraph{Proof of Corollary \ref{coro-1}.} Given the symmetry between the central nodes within the same block, the proof of $(\ref{conv-meas-1})$ is classical using Theorem $\ref{theo2}$ and standard arguments (cf. \cite[Prop.2.2 ]{Sznitman91}). In the following, we prove $(\ref{conv-meas-2})$.

Let $\bar{\mu}_j^{N,p}=\frac{1}{N_j^p}\sum_{n\in C_j^p}\delta_{\bar{X}_{n,j}^p}$. Recall that the bounded-Lipschitz metric $d_{BL}$ metrizes the weak convergence on $\mathcal{M}_1(\mathcal{D}([0,T],\mathcal{Z}))$. Therefore, in order to prove the convergence $(\ref{conv-meas-2})$, it suffices to prove that $d_{BL}(\mu_j^{p,N},\bar{\mu}_j^{p,N})\Rightarrow 0$ and that $\bar{\mu}_j^{N,p}\Rightarrow \mu_j^p$.  First, note that

\begin{equation}
\begin{split}
\mathbb{E}\bigg[d_{BL}\big(\mu_j^{p,N},\bar{\mu}_j^{p,N}\big)\bigg]&=\mathbb{E}\bigg[\sup_{g\in Lip(\mathcal{Z})}\big|\langle \mu_j^{p,N},g\rangle-\langle \bar{\mu}_j^{p,N},g\rangle\big|\bigg]\\
                                                                                                  &=\mathbb{E}\bigg[\sup_{g\in Lip(\mathcal{Z})}\bigg|\frac{1}{N_j^p}\sum_{n\in C_j^p}\big(g(X_{n,j}^{p})-g(\bar{X}_{n,j}^p)\big)\bigg|\bigg] \\
                                                                                                  &\leq \frac{1}{N_j^p}\sum_{n\in C_j^p}\mathbb{E}\big[\big|X_{n,j}^{p}-\bar{X}_{n,j}^p\big|_T\big] \\
                                                                                                  &\leq \max_{n\in C_j^p}\mathbb{E}\big|X_{n,j}^{p}-\bar{X}_{n,j}^p\big|_T,
\end{split}
\label{chaos-ineq7}
\end{equation}
which goes to zero according to $(\ref{L_1-norm})$. Thus, $d_{BL}(\mu_j^{p,N},\bar{\mu}_j^{p,N})\Rightarrow 0$ as $N\rightarrow\infty$. It remains to show that $\bar{\mu}_j^{N,p}\Rightarrow \mu_j^p$ as $N\rightarrow\infty$. Since the stochastic processes $\{\bar{X}_{n,j}^p, n\in C_j^p\}$ are i.i.d., for any continuous and bounded function $g\in C_b(\mathcal{Z})$ one finds

\begin{equation}
\begin{split}
\mathbb{E}\bigg( \langle \bar{\mu}_j^{p,N},g\rangle-\langle \mu_j^p,g\rangle\bigg)^2&=\mathbb{E}\bigg(\frac{1}{N_j^p}\sum_{n\in C_j^p}(g(\bar{X}^p_{n,j})-\langle \mu_j^p,g\rangle)\bigg)^2\\
                                                                                                                 &=\mathbb{E}\bigg(\frac{1}{(N_j^p)^2}\sum_{n\in C_j^p}(g(\bar{X}^p_{n,j})-\langle \mu_j^p,g\rangle)^2\bigg)\\
                                                                                                                 &\leq\frac{1}{N_j^p}4\|g\|_{\infty}^2,
\end{split}
\label{chaos-ineq8}
\end{equation}
which goes to zero given the boundedness of $g$. Therefore, $\bar{\mu}_j^{p,N}$ converges weakly to $\mu_j^p$ as $N\rightarrow\infty$. Thus, combining the two convergence results we conclude that $\mu_j^{p,N}$ converges weakly to $\mu_j^{p}$ as $N\rightarrow\infty$. The corollary is proved.       \carre


\section{Large deviations over a finite time duration}
\label{large-dev-sec}
We investigate in this section the large deviation principles of the interacting particle system introduced in Section \ref{model} over finite time duration. For sake of simplicity, we restrict ourselves to the case of a complete peripheral subgraph, that is, where all peripheral nodes in the system are connected. We follow the classical approach developed in \cite{Daw+Gart87} and adapted to the context of jump processes in \cite{Leonard95}. See also \cite{Feng94(1),Feng94(2)} for a related approach.

We make the following assumptions throughout the section.

\begin{ass}
\begin{enumerate}
\item The peripheral subgraph is complete, that is, for any two peripheral nodes $n,m\in\bigcup\limits_{1\leq j\leq r} C_{j}^p$, there exists an edge $(n,m)\in\Xi$ connecting $n$ and $m$.
\item The mappings $\lambda^c_{z,z'}$ and $\lambda^p_{z,z'}$ introduced in $(\ref{lambda-c-func})$ and $(\ref{lambda-p-func})$ are uniformly bounded away from zero, that is, there exists $c > 0$ such that, for all $\nu,\mu_1,\ldots,\mu_r\in\mathcal{M}_1(\mathcal{Z})$ and all $(z, z')\in\mathcal{E}$, we have $\lambda^c_{z,z'}(\nu,\mu_j)\geq c$ and $\lambda_{z,z'}^p(\nu,\mu_1,\ldots,\mu_r)\geq c$.
\item As $N\rightarrow\infty$, $(\ref{p-regul})$ holds and
\begin{align}
\frac{N_j}{N}\rightarrow \alpha_j,
\label{converg-propo}
\end{align}
for some $\alpha_j\in (0,1)$, where we recall that $N_j$ is the number of nodes in the $j$-th block and $N_j^c$ (resp. $N_j^p$) is the number of central (resp. peripheral) nodes in the $j$-th block.
\end{enumerate}
\label{bound-ass}
\end{ass}

\begin{rem}
\begin{enumerate}
\item From  $(\ref{lambda-c-func})$ and $(\ref{lambda-p-func})$,  together with Remark $\ref{rem-regul-cond}$, the functions $\lambda^c_{z,z'}$ and $\lambda^p_{z,z'}$ are Lipschitz.
\item Since $\mathcal{M}_1(\mathcal{Z})$ is compact and the rate functions $\lambda^c_{z,z'}$  and $\lambda^p_{z,z'}$ are continuous and Lipschitz, the rates are uniformly bounded from above, that is, there exists a constant $C <\infty$ such that for all $\nu,\mu_1,\ldots,\mu_r\in\mathcal{M}_1 (\mathcal{Z})$, and all $(z,z')\in\mathcal{E}$, we have $\lambda^c_{z,z'}(\nu,\mu_j)\leq C$ and $\lambda^p_{z,z'}(\nu,\mu_1,\ldots,\mu_r)\leq C$.
\item To facilitate the reading, we omitted subscripts indicating the dependence of rate functions $\lambda^c_{z,z'}$ and $\lambda^p_{z,z'}$ on the proportions.
\item We use again throughout this section the convention that $N$ goes to infinity when both $\min_{1\leq j\leq r}N_j^c$ and $\min_{1\leq j\leq r}N_j^p$ goes to infinity.
\end{enumerate}
\end{rem}

 Let $M^N\in\mathcal{M}_1(\mathcal{D}([0,T],\mathcal{Z}))\times\cdots\times\mathcal{M}_1(\mathcal{D}([0,T],\mathcal{Z}))$ denote the vector of empirical measures  defined by
\begin{equation}
\begin{split}
M^N&=\bigg(M_1^{c,N},M_1^{p.N},\cdots,M_r^{c,N},M_r^{p,N}\bigg)=\bigg(\frac{1}{N_1^c}\sum_{n\in C^c_1}\delta_{X_n},\frac{1}{N_1^p}\sum_{n\in C^p_1}\delta_{X_n},\ldots, \frac{1}{N^c_r}\sum_{n\in C^c_r}\delta_{X_n}, \frac{1}{N^p_r}\sum_{n\in C^p_r}\delta_{X_n}  \bigg),
\label{bar-vector}
\end{split}
\end{equation}
where $X^N=(X_1,\cdots,X_N)\in\mathcal{D}([0,T],\mathcal{Z}^N)$ denotes the full description of the $N$ particles and $M_j^{c,N}$ (resp.  $M_j^{p,N}$) is the empirical measure of the central (resp. peripheral) nodes of the $j$-th block, for $1\leq j\leq r$. With a slight abuse of notations, denote by $G_N$ the mapping that takes the full description $X^N$ to the empirical measures vector $M^N$, that is,
\begin{align*}
G_N: (X_n, 1\leq n\leq N)\in\mathcal{D}([0,T],\mathcal{Z}^N)\rightarrow\bigg(\frac{1}{N_1^c}\sum_{n\in C^c_1}\delta_{X_n},\frac{1}{N_1^p}\sum_{n\in C^p_1}\delta_{X_n},\ldots, \frac{1}{N^c_r}\sum_{n\in C^c_r}\delta_{X_n}, \frac{1}{N^p_r}\sum_{n\in C^p_r}\delta_{X_n}  \bigg).
\end{align*}
Thus, $M^N=G_N(X^N)$. Denote by $\mathbb{P}_{z^N}^N$ the law of $X^N$, where $z^N=(z_1,\cdots,z_N)$ is the initial condition. Note that the distribution of the empirical vector $M^N$ depends on the initial condition only through its empirical vector defined by
\begin{align}
\nu_N=\left(\nu_N^{1,c},\nu_N^{1,p},\ldots,\nu_N^{r,c},\nu_N^{r,p} \right)=\bigg(\frac{1}{N_1^c}\sum_{n\in C_1^c}\delta_{z_n},\frac{1}{N_1^c}\sum_{n\in C_1^p}\delta_{z_n},\ldots,\frac{1}{N_r^c}\sum_{n\in C_r^c}\delta_{z_n},\frac{1}{N_r^p}\sum_{n\in C_r^p}\delta_{z_n}\bigg).
\label{init-emp-vect}
\end{align}

 Moreover, define by $P_{\nu^N}^N$ the distribution  of $M^N$ which is the pushforward of $\mathbb{P}_{z^N}^N$ under  the mapping $G_N$, that is,  $P_{\nu^N}^N=\mathbb{P}_{z^N}^N\circ G_N^{-1}$. 
 
Let introduce now the $\mathcal{M}_1(\mathcal{Z})\times\cdots\times\mathcal{M}_1(\mathcal{Z})$-valued empirical process
\begin{equation}
\begin{split}
\mu^N: t\in [0,T]\longrightarrow \mu^N(t)&=\left(\mu_1^{c,N}(t),\mu_1^{p,N}(t),\cdots,\mu_r^{c,N}(t),\mu_r^{p,N}(t)\right)\\
                                     &=\bigg(\frac{1}{N_1^c}\sum_{n\in C^c_1}\delta_{X_n(t)},\frac{1}{N_1^p}\sum_{n\in C^p_1}\delta_{X_n(t)},\ldots, \frac{1}{N^c_r}\sum_{n\in C^c_r}\delta_{X_n(t)}, \frac{1}{N^p_r}\sum_{n\in C^p_r}\delta_{X_n(t)}  \bigg),
                                     \label{emp-proc}
\end{split}
\end{equation}
 and define by $\gamma_N$ the corresponding mapping  that takes a full description $X^N \in\mathcal{D}([0,T], \mathcal{Z}^N)$ of the $N$ particles of the system to the empirical process vector $\mu^N$, that is,
\begin{align*}
\gamma_N : (X_n, 1\leq n\leq N )\in \mathcal{D}([0,T], \mathcal{Z}^N)\rightarrow \mu^N: [0,T] \rightarrow \mathcal{M}_1(\mathcal{Z})\times\cdots\mathcal{M}_1(\mathcal{Z}).
\end{align*}

Observe that $\mu^N(0)=\nu_N$ and that $\mu^N(t)$ is the projection $\pi_t(M^N)$ at time $t$, that is,
\begin{align*}
\mu^N=\pi(M^N)=\pi(G_N(X^N))=\gamma_N(X^N),
\end{align*}
where the notation $\pi$ denotes, again with a slight abuse of notation, both the vector projection
\begin{align*}
\pi:\left(\mathcal{M}_1(\mathcal{D}([0,T],\mathcal{Z}))\right)^{2r}\rightarrow\left(\mathcal{D}([0,T],\mathcal{M}_1(\mathcal{Z}))\right)^{2r},
\end{align*}
 and the component projection
 \begin{align*}
\pi:\mathcal{M}_1(\mathcal{D}([0,T],\mathcal{Z}))\rightarrow\mathcal{D}([0,T],\mathcal{M}_1(\mathcal{Z})).
\end{align*}

Finally, denote by $p_{\nu_N}^N$ the distribution of $\mu^N$ which is the pushforward $p_{\nu_N}^N=\mathbb{P}_{z^N}^N\circ\gamma_N^{-1}$. Note that, since $\mu^N=\pi(M^N)$, we can also write $p_{\nu_N}$ as the pushforward $p_{\nu_N}^N=P_{\nu_N}^N\circ\pi^{-1}$.

The goal of this section is to study the large deviation principles for the sequences of probability measures $(P_{\nu^N}^N,N\geq 1)$ and $(p_{\nu_N}^N,N\geq 1)$. The two main results are Theorems $\ref{large-dev-meas}$ and $\ref{large-dev-emp-proc}$.

\subsection{Large deviation principle for the empirical measure}

We start by investigating the large deviation principles of the sequence $(P_{\nu^N}^N,N\geq 1)$.

\subsubsection{The Radon-Nikodym derivative}
Consider first the hypothetical non-interacting case. Suppose that all the nodes are independent of each other  and that the color of each node changes with a constant rate equal to 1 for all allowed transitions $(z,z')\in\mathcal{E}$, and all other transition rates are zero. Denote by $P_{z_0}$ the marginal law on $\mathcal{D} ([0,T],\mathcal{Z})$ of this process with initial condition $z_0$.  Therefore, $P_{z_0}$ is the unique solution to the martingale problem in $\mathcal{D}([0,T],\mathcal{Z})$, associated with the generator $\mathcal{L}^0$ operating on bounded measurable functions $\phi$ on $\mathcal{Z}$ according to

\begin{align*}
\mathcal{L}^0\phi(z)=\sum_{z':(z,z')\in\mathcal{E}}1.(\phi (z')-\phi (z)),
\end{align*}
and the initial condition $z_0$. Given that the transition rates are upper bounded and that
\begin{align*}
\sup_{z\in\mathcal{E}}\sum_{z':(z,z')\in\mathcal{E}}|z'-z|<\digamma (1+z)
\end{align*}
 for some constant $\digamma$, then there exists a unique solution of the martingale problem for $(\mathcal{L}^0,z_0)$ (cf. \cite[Prob.~4.11.15]{Eth+Kurt86}).

For any $\eta,\rho_1,\ldots, \rho_r$ in $\mathcal{D}([0,T], \mathcal{M}_1 (\mathcal{Z}))$, let $R_{z_0}^c(\eta,\rho_j)$ be the unique solution of the martingale problem in $\mathcal{D}([0,T],\mathcal{Z})$, associated with the time-varying generator

\begin{align}
\mathcal{L}^c_{\eta (t),\rho_j(t)}\phi (z)=\sum_{z':(z,z')\in\mathcal{E}}\lambda^c_{z,z'}(\eta (t),\rho_j(t))(\phi (z')-\phi (z)),
\label{gener-c}
\end{align}
and the initial condition $z_0$. Similarly, let $R_{z_0}^p(\eta,\rho_1,\ldots,\rho_r)$ be the unique solution of the martingale problem in $\mathcal{D}([0,T],\mathcal{Z})$, associated with the time varying generator
\begin{align}
\mathcal{L}^p_{\eta (t),\rho_1 (t),\ldots,\rho_r(t)}\phi (z)=\sum_{z':(z,z')\in\mathcal{E}}\lambda^p_{z,z'}(\eta (t),\rho_1 (t),\ldots,\rho_r(t))(\phi (z')-\phi (z)),
\label{gener-p}
\end{align}
and the initial condition $z_0$. Again by the upper boundedness of $\lambda^c_{z,z'}$ and $\lambda^p_{z,z'}$, the uniqueness of $R_{z_0}^c(\eta,\rho_j)$ and $R_{z_0}^p(\eta,\rho_1,\ldots,\rho_r)$ follows (see again \cite[Prob.~4.11.15]{Eth+Kurt86}). Therefore, the density of $R_{z_0}^c(\eta,\rho_j)$ and $R_{z_0}^p(\eta,\rho_1,\ldots,\rho_r)$  with respect to $P_{z_0}$ can be written as (see \cite[eqn.~(2.4)]{Leonard95})

\begin{align}
\frac{d R_{z_0}^c(\eta,\rho_j)}{d P_{z_0}}(x)=\exp\{h_1(x,\eta,\rho_j)\}\quad\text{and}\quad \frac{d R_{z_0}^p(\eta,\rho_1,\ldots,\rho_r)}{d P_{z_0}}(x)=\exp\{h_2(x,\eta,\rho_1,\ldots,\rho_r)\},
\label{dens-Rfunc}
\end{align}
where
\begin{align}
h_1(x,\eta,\rho_j)&= \sum_{0\leq t\leq T}\mathds{1}_{\{x_t\neq x_{t-}\}}\log \bigg(\lambda^c_{x_{t-},x_t}(\eta (t-),\rho_j(t-))\bigg)\nonumber \\
       &\qquad-\int_0^T\bigg(\sum_{z:(x_t,z)\in\mathcal{E}}\lambda^c_{x_t,z}(\eta(t),\rho_j(t))-1\bigg)dt,
       \label{h1-func}
\end{align}
and
\begin{align}
h_2(x,\eta,\rho_1,\ldots,\rho_r)&=\sum_{0\leq t\leq T}\mathds{1}_{\{x_t\neq x_{t-}\}}\log \bigg(\lambda^p_{x_t,x_{t-}}(\eta (t-),\rho_1(t-),\ldots,\rho_r(t-)\bigg)\nonumber \\
             &\qquad-\int_0^T\bigg(\sum_{z:(x_t,z)\in\mathcal{E}}\lambda^p_{x_t,z}(\eta (t),\rho_1(t),\ldots,\rho_r(t))-1\bigg)dt.
             \label{h2-func}
\end{align}

Consider now the system of the $N$ non-interacting particles where the $n$-th particle's law is $P_{z_n}$ with the initial condition being $z_n$. The law of such a system is the product distribution $\mathbb{P}_{z^n}^{0,(N)}=\otimes_{n=1}^NP_{z_n}$. Moreover, the distribution of the corresponding empirical vector is given by $P_{\nu^N}^{0,N}=\mathbb{P}_{z^N}^{0,N}\circ G_N^{-1}$ where  $\nu_N$ is the initial empirical vector $(\ref{init-emp-vect})$. Therefore, by applying a analogous argument of the Cameron-Martin-Girsanov formula in the case of stochastic integrals with respect to point processes (see e.g. \cite[Lem.~3.7]{Daw+Zhen91} or \cite[eqn.~(2.8)]{Leonard95}), one can compute the Radon-Nikodym derivative $dP_{\nu^N}^{N}/dP_{\nu^N}^{0,N}$ at any $\mathbf{Q}=(Q_{j}^c,Q_{j}^p,\cdots, Q_{r}^c,Q_{r}^p)\in \mathcal{M}_1(\mathcal{D}([0,T],\mathcal{Z}))\times \cdots\times  \mathcal{M}_1(\mathcal{D}([0,T],\mathcal{Z}))$ as follows:
\begin{equation}
\begin{split}
\frac{dP_{\nu^N}^{N}}{dP_{\nu^N}^{0,N}}(\mathbf{Q})&=\exp\bigg\{\sum_{j=1}^r\bigg[N_j^c\int_{D([0,T],\mathcal{Z})}h_1(x, \pi(Q_{j}^c),\pi(Q_{j}^p))Q_{j}^c(dx)\\
&\qquad\qquad+N_j^p\int_{D([0,T],\mathcal{Z})}h_2(x, \pi(Q_{j}^c), \pi(Q_{1}^p),\ldots, \pi(Q_{r}^p))Q_{j}^p(dx)\bigg]\bigg\}\\
                                                   &=\exp\big\{N h(\mathbf{Q})\big\},
\label{rad-nik}
\end{split}
\end{equation}
where
\begin{equation}
\begin{split}
h(\mathbf{Q})&=\sum_{j=1}^r\bigg[\frac{N_j^c}{N}\int_{D([0,T],\mathcal{Z})}h_1(x,\pi(Q_{j}^c),\pi(Q_{j}^p))Q_{j}^c(dx)\\
&\qquad\qquad+\frac{N_j^p}{N}\int_{D([0,T],\mathcal{Z})}h_2(x,\pi(Q_{j}^c), \pi(Q_{1}^p),\ldots, \pi(Q_{r}^p))Q_{j}^p(dx)\bigg].
\label{h-func}
\end{split}
\end{equation}

\subsubsection{The spaces and topologies of interest}
 
 We introduce here the spaces and topologies of interest following \cite{Leonard95} and \cite{Bork+Sund2012}. Consider the Polish space $(\mathcal{X},d)$ where
\begin{align*}
\mathcal{X}=\bigg\{x\in \mathcal{D}([0,T], \mathcal{Z})| &\sum_{0\leq t\leq T}\mathds{1}_{x_t\neq x_{t-}}<+\infty,\\
&\text{ and for each $t\in(0, T ]$ with $x_t\neq x_{t-}$, we have $(x(t-),x(t))\in\mathcal{E}$}\bigg\},
\end{align*}
and the metric $d$ is defined by
\begin{align}
d(x,y)= d_{Sko}(x,y)+|(\varphi (x)-\varphi (y)|,\quad x,y\in\mathcal{X},
\label{metr-X}
\end{align}
with $\varphi(x)=\sum_{0\leq t\leq T}\mathds{1}_{x_t\neq x_{t-}}$ denoting the number of jumps and $d_{Sko}$ standing for the  Skorokhod complete metric (see \cite[Sec.~12]{Billi99}). For this topology, the function $\varphi$ is continuous and two paths are close to each other if they have the same number of jumps and if they are Skorokhod-close  \cite[p.~299]{Leonard95}. For any function $f:\mathcal{X}\rightarrow\mathbb{R}$ define

\begin{align}
\|f\|_{\varphi}=\sup_{x\in\mathcal{X}}\frac{f(x)}{1+\varphi (x)},
\label{var-norm}
\end{align}
and denote
\begin{align}
C_{\varphi}(\mathcal{X}) = \big\{f | f : \mathcal{X}\rightarrow\mathbb{R}\text{ is continuous and }\|f\|_{\varphi}<\infty\big\},
\label{C-varphi}
\end{align}
\begin{align}
\mathcal{M}_{1,\varphi}(\mathcal{X})=\bigg\{Q\in\mathcal{M}_1(\mathcal{X})\big|\int_{\mathcal{X}}\varphi dQ<+\infty \bigg\}.
\label{M-varphi}
\end{align}

 We endow the set $\mathcal{M}_{1,\varphi}(X)$ with the weak$^*$ topology $\sigma(\mathcal{M}_{1,\varphi}(\mathcal{X}), C_{\varphi}(\mathcal{X}))$, the weakest topology under which $Q_N\rightarrow Q$ as $N\rightarrow +\infty$ if and only if
\begin{align*}
\int_{\mathcal{X}} f dQ_N\rightarrow \int_{\mathcal{X}} f dQ \quad\text{for each $f\in C_{\varphi}(\mathcal{X})$}.
\end{align*}

For a measure $\nu=(\nu^{1,c},\nu^{1,p},\ldots,\nu^{r,c},\nu^{r,p})\in\mathcal{M}_1(\mathcal{Z})\times\cdots\mathcal{M}_1(\mathcal{Z})$ we define, for all $1\leq j\leq r$ and $\iota\in\{c,p\}$, the mixture
\begin{align*}
dP_{j,\iota}(x)=\sum_{z_0\in\mathcal{Z}}\nu^{j,\iota}(z_0)dP_{z_0}(x).
\end{align*}

 Moreover, let $R^c(\eta,\rho_j)$ and $R^p(\eta,\rho_1,\ldots,\rho_r)$ be the mixtures given by
\begin{equation}
\begin{split}
dR^c(\eta,\rho_j)(x)&=\sum_{z_0\in\mathcal{Z}}\nu^{j,c}(z_0)dR^c_{z_0}(\eta,\rho_j)(x),\\
dR^p(\eta,\rho_1,\ldots,\rho_r)(x)&=\sum_{z_0\in\mathcal{Z}}\nu^{j,p}(z_0)dR^p_{z_0}(\eta,\rho_1,\ldots,\rho_r)(x).
\end{split}
\end{equation}

Finally, let introduce the relative entropy $H: \mathcal{M}_{1,\varphi} (\mathcal{X})\rightarrow [0,+\infty]$ of $Q$ with respect to $P$ as follows: 
\begin{align}
H(Q|P)=\left\{\begin{array}{rcl}
\int_{\mathcal{X}}\log (\frac{dQ}{dP})dQ & \mbox{if $Q \ll P$}, & \\
+\infty&\text{otherwise}. & 
\end{array}
\right.
\label{relat-entro}
\end{align}

\subsubsection{Large deviation principle for the empirical measure vector}

The next theorem gives the large deviation principle for the sequence $(P_{\nu_N}^N,N\geq 1)$.

\begin{theo}
Let the space $\mathcal{M}_{1,\varphi}(\mathcal{X})$ be equipped with the weak$^*$ topology $\sigma(\mathcal{M}_{1,\varphi}(\mathcal{X}), C_{\varphi}(\mathcal{X}))$. Moreover, suppose that the initial condition: $\nu_N\rightarrow\nu$ weakly as $N\rightarrow\infty$.  Then, the sequence $(P_{\nu_N}^N,N\geq 1)$ satisfies the large deviation principle in the space $\mathcal{M}_{1,\varphi}(\mathcal{X})\times\cdots\mathcal{M}_{1,\varphi}(\mathcal{X})$, endowed with the product topology, with speed $N$ and the good rate function $I(\mathbf{Q}) = L(\mathbf{Q})-h(\mathbf{Q})$, where the function $h(\mathbf{Q})$ is given by $(\ref{h-func})$ and $L:\mathcal{M}_{1,\varphi}(\mathcal{X})\times\cdots\times \mathcal{M}_{1,\varphi}(\mathcal{X})\rightarrow [0,\infty]$ is defined as
\begin{equation}
\begin{split}
L(\textbf{Q})=\alpha_1p_1^c J^{1,c}(Q_{1}^c)+ \alpha_1p_1^p J^{1,p}(Q_{1}^p)+\cdots+\alpha_rp_r^c J^{r,c}(Q_{r}^c)+ \alpha_rp_r^p J^{r,p}(Q_{r}^p),
\label{rate-func-noninter}
\end{split}
\end{equation}
with, for each $1\leq j\leq r$, $\iota\in\{c,p\}$ and $Q\in\mathcal{M}_{1,\varphi}(\mathcal{X})$,
\begin{align}
J^{j,\iota}(Q)= \sup_{f\in C_{\varphi}(\mathcal{X})}\bigg[\int_{\mathcal{X}}fdQ-\sum_{z_0\in\mathcal{Z}}\nu^{j,\iota} (z_0)\log\int_{\mathcal{X}}e^fdP_{z_0}\bigg],
\label{Daw-Gar-rate}
\end{align}
and $\alpha_j,p_j^c,p_j^p$ being given in $(\ref{converg-propo})$ and $(\ref{p-regul})$.  Furthermore, for each $\mathbf{Q}\in\mathcal{M}_{1,\varphi}(\mathcal{X})\times\cdots\mathcal{M}_{1,\varphi}(\mathcal{X})$, the rate function $I(\mathbf{Q})$ admits the representation 
\begin{align}
& I(\mathbf{Q})= \nonumber \\ 
& \left\{ \begin{array}{ll}
\sum_{j=1}^r\bigg[\alpha_jp_j^cH\bigg(Q_{j}^c\big|R^c(\pi(Q_{j}^c),\pi(Q_{j}^p))\bigg) +\alpha_jp_j^pH\bigg(Q_{j}^p\big| R^p(\pi(Q_{j}^c),\pi(Q_{1}^p),\ldots,\pi(Q_{r}^p))\bigg)\bigg],  & \text{if $\mathbf{Q}\circ\pi_0^{-1}=\nu$}, \\
+\infty, &\text{otherwise}.
\end{array}
\right.
\label{empi-meas-rate}
\end{align}
\label{large-dev-meas}
\end{theo}

\begin{rem}
This is a generalization of \cite[Th.~2.1]{Leonard95} to our multi-population setting. Also, while \cite{Leonard95} studied the case where $z_n=z_0$ for some fixed $z_0$ so that $\nu_N=\delta_{z_0}$, we consider, as in \cite[Th.~3.1]{Bork+Sund2012}, more general starting points for each particle, provided that the initial empirical vector $\nu_N$ converges weakly to $\nu=\left(\nu^{1,c},\nu^{1,p},\cdots,\nu^{r,c},\nu^{r,p}\right)$. Moreover, similar to \cite{Bork+Sund2012}, we consider here the case where not all transitions are allowed, but only those in $\mathcal{E}$, the set of directed edges in the graph $(\mathcal{Z},\mathcal{E})$. \\
\indent Note that, from Definition $\ref{defi-multi-chaos}$, the weak convergence of the initial empirical vector $\nu_N$ towards $\nu$ is amount to the the assertion that the initial conditions $(X_n^{c,N}(0),X_m^{p,N}(0),n\in C_j^c,m\in C_j^p;1\leq j\leq r)$ are  $\nu^{1,c}\otimes\nu^{1,p}\cdots \nu^{r,c}\otimes\nu^{r,p}$-multi-chaotic (cf. \cite{Sznitman91}).
\end{rem}

\paragraph{Proof of Theorem \ref{large-dev-meas}.} The proof of Theorem \ref{large-dev-meas} is based on  the generalization of Sanov's theorem for empirical measures on Polish spaces due to Dawson and G\"artner \cite{Daw+Gart87}, the Girsanov transformation, and the Laplace-Varadhan principle \cite{Vara84}. We proceed through several lemmas. We follow \cite[Th.~2.1]{Leonard95} and \cite[Th.~3.1]{Bork+Sund2012}.

\subparagraph*{Large deviation for the non-interacting case.}

We first establish a large deviation principle in the non-interacting case.

\begin{lem}
Suppose that the initial condition $\nu_N\rightarrow\nu$  weakly. Let $\mathcal{M}_{1,\varphi}(\mathcal{X})$ be endowed with the  weak$^*$ topology $\sigma(\mathcal{M}_{1,\varphi} (\mathcal{X}), C_{\varphi} (\mathcal{X}))$. Then, the sequence $(P_{\nu_N}^{0,N}, N\geq 1)$ satisfies a large deviation principle in $\mathcal{M}_{1,\varphi}(\mathcal{X})\times\cdots\times \mathcal{M}_{1,\varphi}(\mathcal{X})$, endowed with the product topology, with speed $N$ and the action functional $L:\mathcal{M}_{1,\varphi}(\mathcal{X})\times\cdots\times \mathcal{M}_{1,\varphi}(\mathcal{X})\rightarrow [0,\infty]$, given by $(\ref{rate-func-noninter})$.
\label{non-inter-LDP}
\end{lem}

\proof Fix a given block $1\leq j\leq r$. Denote by $\big(P_{\nu_{N}^{j,c}}^{0,N_j^c}, N_j^c\geq 1\big)$ and $\big(P_{\nu_{N}^{j,p}}^{0,{N_j^p}}, N_j^p\geq 1\big)$ the sequences of probability distributions of the local empirical measures $M_j^{N,c}$ and $M_j^{N,p}$ of the central and peripheral nodes of the $j$-th block, respectively. Note that in the non-interacting case, the transition rate from any state to any other state is bounded by $1$. Therefore, the family of probability measures $\{P_z:z\in\mathcal{Z}\}$ is a subset of $\mathcal{M}_{1,\varphi}(\mathcal{X})$. Moreover, for any continuous function $F\in C_{\varphi}(\mathcal{X})$, the integral $\int F(y)P_{z_0}(dy)$ depends continuously upon $z_0$ and then,  $\{P_{z_0}:z_0\in\mathcal{Z}\}$  is a Feller continuous family of probability measures on $\mathcal{X}$. Now, since $\nu^{j,c}_N\rightarrow\nu^{j,c}$ and $\nu^{j,p}_N\rightarrow\nu^{j,p}$, by applying the generalization of Sanov's theorem  \cite[Th.~3.5]{Daw+Gart87}, we find that  both the sequences $\big(P_{\nu_{N}^{j,c}}^{0,N_j^c}, N_j^c\geq 1\big)$ and $\big(P_{\nu_{N}^{j,p}}^{0,{N_j^p}}, N_j^p\geq 1\big)$ satisfy the large deviation principle in $\mathcal{M}_{1,\varphi} (\mathcal{X})$, endowed with the weak$^*$ topology $\sigma(\mathcal{M}_{1,\varphi} (\mathcal{X}), C_{\varphi} (\mathcal{X}))$, with speeds $N_j^c$ and $N_j^p$, respectively, and good rate functions $J^{j,c}(Q)$ and $J^{j,p}(Q)$ defined by (\ref{Daw-Gar-rate}). Let $\mathcal{K}^c_1,\mathcal{K}^p_1,\ldots,\mathcal{K}_r^c,\mathcal{K}_r^p\in\mathcal{B}(\mathcal{M}_{1,\varphi}(\mathcal{X}))$ be closed Borelian sets. By independence, one has

\begin{equation}
\begin{split}
P_{\nu_{N}}^{0,{N}}\bigg\{M^N\in \prod_{j=1}^r(\mathcal{K}^c_j\times \mathcal{K}^p_j)\bigg\}= \prod_{j=1}^r\left( P_{\nu_{N_j^p}}^{0,{N_j^p}}\left\{M_j^{N,c}\in \mathcal{K}_j^c\right\}\times P_{\nu_N^{j,p}}^{0,{N_j^p}}\left\{M_j^{N,p}\in \mathcal{K}_j^p\right\}\right).
\end{split}
\end{equation}

Therefore, by Assumption \ref{bound-ass} we get

\begin{equation}
\begin{split}
\limsup_{N\rightarrow\infty}\frac{1}{N}\log P_{\nu_{N}}^{0,{N}}\bigg(\prod_{j=1}^r\mathcal{K}^c_j\times \mathcal{K}^p_j\bigg)&=\limsup_{N\rightarrow\infty}\frac{1}{N}\log\bigg(\prod_{j=1}^rP_{\nu_N^{j,c}}^{0,{N_j^c}}(\mathcal{K}_j^c)P_{\nu_N^{j,p}}^{0,{N_j^p}}(\mathcal{K}_j^p)  \bigg)\\
                &= \limsup_{N\rightarrow\infty}\sum_{j=1}^r\bigg(\frac{N_j}{N}\frac{N_j^c}{N_j}\frac{1}{N_j^c}\log P_{\nu_N^{j,c}}^{0,{N_j^c}}(\mathcal{K}_j^c)+\frac{N_j}{N}\frac{N_j^p}{N_j}\frac{1}{N_j^p}\log P_{\nu_N^{j,p}}^{0,{N_j^p}}(\mathcal{K}_j^p)\bigg)\\
                &\leq \sum_{j=1}^r\bigg(\alpha_jp_j^c\limsup_{N_j^c\rightarrow\infty}\frac{1}{N_j^c}\log P_{\nu_N^{j,c}}^{0,{N_j^c}}(\mathcal{K}_j^c)+\alpha_jp_j^p\limsup_{N_j^p\rightarrow\infty}\frac{1}{N_j^p}\log P_{\nu_N^{j,p}}^{0,{N_j^p}}(\mathcal{K}_j^p)\bigg)\\
                &\leq \sum_{j=1}^r\bigg( -\alpha_jp_j^c\inf_{Q_{j}^c\in \mathcal{K}_j^c}J^{j,c}(Q_{j}^c)-\alpha_jp_j^p\inf_{Q_{j}^p\in \mathcal{K}_j^p}J^{j,p}(Q_{j}^p)\bigg)\\
                &=-\inf_{\substack{Q_{j}^c\in \mathcal{K}_1^c\\ Q_{j}^p\in \mathcal{K}_1^p\\ \vdots\\ Q_{r}^c\in \mathcal{K}_r^c\\ Q_{r}^p\in \mathcal{K}_r^p }}\sum_{j=1}^r\left( \alpha_jp_j^c J^{j,c}(Q_{j}^c)+\alpha_jp_j^p J^{j,p}(Q_{j}^p)\right).
\end{split}
\end{equation}

Similar arguments allow to prove the LDP lower bound, which conclude the proof.  \carre
\\
\\
The next result gives a characterization of the space containing the probability measures satisfying $L(\textbf{Q})<\infty$.

\begin{lem}
If, for a given $\textbf{Q}=(Q_{1}^c,Q_{1}^p,\ldots,Q_{r}^c,Q_{r}^p)\in\mathcal{M}_1(\mathcal{X})\times\cdots\times\mathcal{M}_1(\mathcal{X}) $,  the action functional $L(\mathbf{Q})<\infty$, then:
\begin{enumerate}
\item  $\textbf{Q}\in\mathcal{M}_{1,\varphi}(\mathcal{X})\times\cdots\times\mathcal{M}_{1,\varphi}(\mathcal{X})$.
\item  $ \textbf{Q}\circ\pi_0^{-1}=\nu$. Thus, $\left(\pi_0^{-1}(Q_{1}^c),\pi_0^{-1}(Q_{1}^p),\ldots,\pi_0^{-1}(Q_{r}^c),\pi_0^{-1}(Q_{r}^p)\right)=\left(\nu^{1,c},\nu^{1,p},\cdots,\nu^{r,c},\nu^{r,p}\right)$.
\end{enumerate}
\label{Q-space}
\end{lem}

\proof This is a generalization of \cite[Lem.~5.2]{Bork+Sund2012}. Recall that the function $\varphi (x)=\sum_{0\leq t\leq T}\mathds{1}_{x(t-)\neq x(t)}$ denotes the number of jumps of $x$ in the interval $[0,T]$. From $(\ref{var-norm})$ we have that $\|\varphi\|_{\varphi}\leq 1$. Moreover, $\varphi$ is continuous in the topology induced by the metric $d$ defined in $(\ref{metr-X})$. Hence $\varphi\in C_{\varphi}(\mathcal{X})$. Furthermore, $L(\textbf{Q})<\infty$ implies that, for all $1\leq j\leq r$,

\begin{align}
\int_{\mathcal{X}}\varphi dQ_{j}^c-\sum_{z_0\in\mathcal{Z}}\nu^{j,c}(z)\log\int_{\mathcal{X}}e^{\varphi} dP_{z_0}<\infty,
\label{L_c-ineq}
\end{align}

and

\begin{align}
\int_{\mathcal{X}}\varphi dQ_{j}^p-\sum_{z_0\in\mathcal{Z}}\nu^{j,p}(z)\log\int_{\mathcal{X}}e^{\varphi} dP_{z_0}<\infty.
\label{L_p-ineq}
\end{align}

Now, note that under the non-interacting distribution $P_{z_0}$, the transition rates are bounded by one. Therefore, since the number of allowed transitions from any state is at most equal to $K-1$, $\varphi$ is thus stochastically dominated by a Poisson random variable of rate $(K-1)T$. Therefore, for any initial condition $z_0\in\mathcal{Z}$, we have $1\leq\int_{\mathcal{X}}e^{\varphi} dP_{z_0}<\infty$. It follows from $(\ref{L_c-ineq})$ and $(\ref{L_p-ineq})$ that $\int_{\mathcal{X}}\varphi dQ_{j}^c<\infty$ and $\int_{\mathcal{X}}\varphi dQ_{j}^p<\infty$ for each $1\leq j\leq r$ and so $\textbf{Q}=(Q_{1}^c,Q_{1}^p,\cdots,Q_{r}^c,Q_{r}^p)\in\mathcal{M}_{1,\varphi}(\mathcal{X})\times\cdots\times\mathcal{M}_{1,\varphi}(\mathcal{X})$, which proves the first claim. 

In order to prove the second point, we proceed by contraposition. Suppose that for a given measure $\textbf{Q}$,  $L(\textbf{Q})<\infty$ and $\textbf{Q}\circ\pi_0^{-1}=\nu_{\textbf{Q}}\neq\nu$. Consider the bounded continuous functions $f_1^c(x),f_1^p(x),\ldots,f_r^c(x),f_r^p(x)$ defined on $\mathcal{X}$ and depending on $x$ only through the initial condition, that is, there exist functions $g_1^c,g_1^p,\ldots,g_r^c,g_r^p$ such that, for all $1\leq j\leq r$,
\begin{align*}
f_j^c(x)=g_j^c(\pi_0(x))\quad\text{and}\quad f_j^p(x)=g_j^p(\pi_0(x)).
\end{align*}

Since $\nu_{\textbf{Q}}\neq \nu$, the above functions satisfy the following claim: either
\begin{align}
\sum_zg_j^c(z)\nu_{\textbf{Q}}^{j,c}(z)-\sum_zg_j^c(z)\nu^{j,c}(z)\neq 0,
\label{ineq-fc}
\end{align}
or
\begin{align}
\sum_zg_j^p(z)\nu_{\textbf{Q}}^{j,p}(z)-\sum_zg_j^p(z)\nu^{j,p}(z)\neq 0,
\end{align}
 for at least one $1\leq j\leq r$. Therefore, one can always find, for at least one $j$, an arbitrary large $a_j^c>0$ (or $a_j^p>0$) such that  $\sum_zg_j^c(z)\nu_{\textbf{Q}}^{j,c}(z)-\sum_zg_j^c(z)\nu^{j,c}(z)=a_j^c$ (or $\sum_zg_j^p(z)\nu_{\textbf{Q}}^{j,p}(z)-\sum_zg_j^p(z)\nu^{j,p}(z)=a_j^p$). Indeed, this can be done by flipping the sign of $f_j^c$ (or $f_j^p$) if necessary and scaling the functions. Note that, by the assumption, $f_j^c,f_j^p\in C_{\varphi}(\mathcal{X})$ since they are bounded continuous. Suppose, without loss of generality that, for a  given $j$, $(\ref{ineq-fc})$ is satisfied, then by direct calculations we obtain,

\begin{align*}
\int_{\mathcal{X}}f_j^cdQ_{j}^c-\sum_{z_0\in\mathcal{Z}}\nu^{j,c} (z)\log\int_{\mathcal{X}}e^{f_j^c}dP_{z_0}&=\int_{\mathcal{X}}g_j^c(\pi_0(x))Q_{j}^c(dx)-\sum_{z_0\in\mathcal{Z}}\nu^{j,c} (z)\log\int_{\mathcal{X}}\exp\{g_j^c(\pi_0(x))\}dP_{z_0}\\
                                    &=\sum_{z_0}g_j^c(z_0)\nu_{\textbf{Q}}^{j,c}(z_0)-\sum_{z_0}g_j^c(z)\nu^{j,c}(z_0)=a.
\end{align*}

Hence, since $a>0$ is arbitrary large one  gets that  $J(Q_{j}^c)=\infty$ and then $L(\textbf{Q})=\infty$, which contradicts the condition of the lemma and then proves the second claim. \carre

\subparagraph*{Conditions of applications of the Laplace-Varadhan Lemma.} We established in Lemma $\ref{non-inter-LDP}$ the large deviation principle for the sequence $(P_{\nu_N}^{0,N},N\geq 1)$ in the topological space $\mathcal{M}_{1,\varphi}(\mathcal{X})\times\cdots\times \mathcal{M}_{1,\varphi}(\mathcal{X})$. Moreover, the Radon-Nikodym derivative is given by

\begin{equation}
\begin{split}
\frac{dP_{\nu^N}^{N}}{dP_{\nu^N}^{0,N}}(\mathbf{Q})=\exp\big\{N h(\mathbf{Q})\big\},
\end{split}
\end{equation}
where the function $h(\mathbf{Q})$ is given by $(\ref{h-func})$. Therefore, in order to find the large deviation principle for $(P_{\nu_N}^{N},N\geq 1)$, one can apply the Laplace-Varadhan principle (cf. \cite[Prop.~2.5]{Leonard95}) to $(P_{\nu_N}^{0,N},N\geq 1)$. To this end, two conditions have to be verified: the continuity of the function $h(\mathbf{Q})$ and that, for any $\alpha>0$,
\begin{align}
\limsup_{N\rightarrow\infty}\frac{1}{N}\log\int_{\mathcal{M}_{1,\varphi (\mathcal{X})}\times\cdots\times\mathcal{M}_{1,\varphi (\mathcal{X})}}\exp\left\{ N\alpha |h|   \right\} dP_{\nu_N}^{0,N}<\infty.
\label{vara-cond-1}
\end{align}

The next four lemmas are dedicated to the verification of these two conditions. First, we establish a regularity property for all the probability measures $\textbf{Q}$ satisfying $L(\textbf{Q})<\infty$. This result is a generalization of \cite[Lem.~5.7]{Bork+Sund2012}.

\begin{lem}
Let $\textbf{Q}=(Q_{1}^c,Q_{1}^p,\cdots,Q_{r}^c,Q_{r}^p)\in\mathcal{M}_1(\mathcal{X})\times\cdots\times\mathcal{M}_1(\mathcal{X})$ such that $L(\textbf{Q})<\infty$. Moreover, suppose that the random vector $\textbf{X}=(X_1^c,X_1^p,\cdots,X_r^c,X_r^p)$ is distributed according to $\textbf{Q}$. Then,

\begin{align}
\sup_{t\in [0,T]}\mathbb{E}\bigg[\sup_{u\in \left[t-\alpha ,t+\alpha\right]\cap [0,T]}\left\{\mathds{1}_{\textbf{X}(u)\neq\textbf{X}(u-)}\right\}\bigg]\rightarrow 0\quad\text{as $\alpha\downarrow 0$}.
\label{regul-Q}
\end{align}
\label{regul-lem}
\end{lem}

\proof Note that $\textbf{X}(u)\neq\textbf{X}(u-)$ if $X_j^c(u)\neq X_j^c(u-)$ or $X_j^p(u)\neq X_j^p(u-)$ for at least one $1\leq j\leq r$. Therefore, one obtains, for each $t\in [0,T]$,
\begin{equation}
\begin{split}
\mathbb{E}&\bigg[\sup_{u\in \left[t-\alpha ,t+\alpha\right]\cap [0,T]}\left\{\mathds{1}_{\textbf{X}(u)\neq\textbf{X}(u-)}\right\}\bigg]\leq \mathbb{E}\bigg[\sup_{u\in \left[t-\alpha ,t+\alpha\right]\cap [0,T]}\{\mathds{1}_{X_1^c(u)\neq X_1^c(u-)}\}\bigg]+ \mathbb{E}\bigg[\sup_{u\in \left[t-\alpha ,t+\alpha\right]\cap [0,T]}\{\mathds{1}_{X_1^p(u)\neq X_1^p(u-)}\}\bigg] \\
& +\cdots+ \mathbb{E}\bigg[\sup_{u\in \left[t-\alpha ,t+\alpha\right]\cap [0,T]}\{\mathds{1}_{X_r^c(u)\neq X_r^c(u-)}\}\bigg]+ \mathbb{E}\bigg[\sup_{u\in \left[t-\alpha ,t+\alpha\right]\cap [0,T]}\{\mathds{1}_{X_r^p(u)\neq X_r^p(u-)}\}\bigg].
\end{split}
\label{local-1}
\end{equation}

Moreover, since $L(\textbf{Q})<\infty$ one gets that $J(Q_{j}^c)<\infty$ and  $J(Q_{j}^p)<\infty$ for all $1\leq j\leq r$. Hence, applying \cite[Lem.~5.7]{Bork+Sund2012} to each of the $X_j^c$ and $X_j^p$ with respective marginal distributions $Q_{j}^c$ and $Q_{j}^p$ gives us that, for each $1\leq j\leq r$,

\begin{align}
\sup_{t\in [0,T]}\mathbb{E}\bigg[\sup_{u\in \left[t-\alpha ,t+\alpha\right]\cap [0,T]}\left\{\mathds{1}_{X_j^c(u)\neq X_j^c(u-)}\right\}\bigg]\rightarrow 0\quad\text{as $\alpha\downarrow 0$},
\label{local-2}
\end{align}
 and
\begin{align}
\sup_{t\in [0,T]}\mathbb{E}\bigg[\sup_{u\in \left[t-\alpha ,t+\alpha\right]\cap [0,T]}\left\{\mathds{1}_{X_j^p(u)\neq X_j^p(u-)}\right\}\bigg]\rightarrow 0\quad\text{as $\alpha\downarrow 0$}.
\label{local-3}
\end{align}

Combining $(\ref{local-1})$, $(\ref{local-2})$ and $(\ref{local-3})$ leads to $(\ref{regul-Q})$. \carre

The next lemma  establishes the continuity of the projection $\pi$, which is needed to establish to continuity of the function $h(\mathbf{Q})$.

\begin{lem}
Let $\mathcal{M}_1(\mathcal{D}([0,T],\mathcal{Z}))$ be equipped with its usual weak topology and let $\mathcal{D}([0,T],\mathcal{M}_1(\mathcal{Z}))$ be equipped with the metric,
\begin{align}
\rho_T (\mu,\nu)=\sup_{0\leq t\leq T}\rho_0 (\mu_t,\nu_t),\quad\mu,\nu\in\mathcal{D}([0,T],\mathcal{M}_1(\mathcal{Z})),
\end{align}
 where $\rho_0(\cdot,\cdot)$  is a metric on $\mathcal{M}_1(\mathcal{Z})$ which generates the weak topology $\sigma(\mathcal{M}_1(\mathcal{Z}),C_b(\mathcal{Z}))$. Moreover, let $\mathcal{M}_1(\mathcal{D}([0,T],\mathcal{Z}))\times\cdots\times\mathcal{M}_1(\mathcal{D}([0,T],\mathcal{Z}))$ be endowed with the product topology induced by the product metric. Equivalently, let $\mathcal{D}([0,T],\mathcal{M}_1(\mathcal{Z}))\times\cdots\times\mathcal{D}([0,T],\mathcal{M}_1(\mathcal{Z}))$ be equipped with the product topology obtained from the product metric $N(\rho_T,\cdots,\rho_T)$. Then, the projection:
\begin{align*}
\pi: \mathbf{Q}\in\left(\mathcal{M}_1(\mathcal{D}([0,T],\mathcal{Z}))\right)^{2r}&\rightarrow\pi(\mathbf{Q})=(\mathbf{Q}_t)_{0\leq t\leq T}\in\left(\mathcal{D}([0,T],\mathcal{M}_1(\mathcal{Z}))\right)^{2r}
\end{align*}
is continuous at each $\textbf{Q}\in\left(\mathcal{M}_1(\mathcal{D}([0,T],\mathcal{Z}))\right)^{2r}$ where $L(\textbf{Q})<\infty$.
\label{pi-cont}
\end{lem}

\proof The statement of our lemma resembles the statement of \cite[Lem.~2.8]{Leonard95}. The difference here is that our spaces of interest are the product spaces $\mathcal{M}_1(\mathcal{D}([0,T],\mathcal{Z}))\times\cdots\times\mathcal{M}_1(\mathcal{D}([0,T],\mathcal{Z}))$ and $\mathcal{D}([0,T],\mathcal{M}_1(\mathcal{Z}))\times\cdots\times\mathcal{D}([0,T],\mathcal{M}_1(\mathcal{Z}))$ endowed with product metrics. Moreover, the rate $J(Q)$ in \cite[Lem.~2.8]{Leonard95} is here replaced by $L(\mathbf{Q})$. Therefore,  replacing the norm $|\cdot|$ by the product norm $\|\cdot\|$ adapted to our product spaces context, the proof of our lemma follows verbatim the proof of \cite[Lem.~2.8]{Leonard95} provided that we can prove \cite[eqn.~(2.14)]{Leonard95}. This is done in Lemma \ref{regul-lem}. This concludes the proof.
 \carre
\\
\\
We now state the continuity of the function $h$.
\begin{lem}
The function $h:\mathcal{M}_{1,\varphi}(\mathcal{X})\times\cdots\times\mathcal{M}_{1,\varphi}(\mathcal{X})\rightarrow\mathbb{R}$ defined at $(\ref{h-func})$ is continuous at any $\textbf{Q}$ such that $L(\textbf{Q})<\infty$.
\label{h-continuity}
\end{lem}

\proof This is a generalization of \cite[Lem.~2.9]{Leonard95}. For any $\textbf{Q}\in\mathcal{M}_{1,\varphi}(\mathcal{X})\times\cdots\times\mathcal{M}_{1,\varphi}(\mathcal{X})$ define,
\begin{align}
\theta_{\textbf{Q}}^{j,c}(x)&=\sum_{0\leq t\leq T}\mathds{1}_{x_t\neq x_{t-}}\log\bigg(\sum_{(x_{t-},x(t))\in\mathcal{E}}\lambda_{x_{t-},x_t}^c\left(Q^c_j(t-),Q^p_j(t-)\right)\bigg),\\
\theta_{\textbf{Q}}^{j,p}(x)&=\sum_{0\leq t\leq T}\mathds{1}_{x_t\neq x_{t-}}\log\bigg(\sum_{(x_{t-},x_t)\in\mathcal{E}}\lambda_{x_{t-},x_t}^p(Q^c_j(t-),Q^p_1(t-),\ldots,Q^p_r(t-) )\bigg),\\
\gamma_{\textbf{Q}}^{j,c}(x)&=\int_0^T\bigg(\sum_{z:(x_t,z)\in\mathcal{E}}\lambda^c_{x_t,z}(Q_{j}^c (t),Q_{j}^p (t))-1\bigg)dt,\\
\gamma_{\textbf{Q}}^{j,p}(x)&=\int_0^T\bigg(\sum_{z:(x_t,z)\in\mathcal{E}}\lambda^p_{x_t,z}(Q_{j}^c (t),Q_{1}^p (t),\ldots,Q_{r}^p (t))-1\bigg)dt.
\end{align}

Note that the function $h$ given by $(\ref{h-func})$ can be rewritten using the functions $\theta_{\textbf{Q}}^{j,c}(x)$,$\theta_{\textbf{Q}}^{j,p}(x)$, $\gamma_{\textbf{Q}}^{j,c}(x)$ and $\gamma_{\textbf{Q}}^{j,p}(x)$ as follows:
\begin{equation}
    h(\mathbf{Q}) =\sum_{j=1}^r\bigg[\frac{N_j^c}{N}\int_{\mathcal{X}}\bigg(\theta_{\textbf{Q}}^{j,c}(x)-\gamma_{\textbf{Q}}^{j,c}(x)\bigg)Q_{j}^c(dx)
 + \frac{N_j^p}{N}\int_{\mathcal{X}}\bigg(\theta_{\textbf{Q}}^{j,p}(x)-\gamma_{\textbf{Q}}^{j,p}(x)\bigg)Q_{j}^p(dx)\bigg].
\label{h-theta-gamma}
\end{equation}

Therefore, in order to show the continuity of $h(\mathbf{Q})$, we show that, for any $1\leq j\leq r$, the functions
\begin{align*}
\textbf{Q}\rightarrow\int_{\mathcal{X}}\theta_{\textbf{Q}}^{j,c}(x)Q_{j}^c(dx),&\qquad \textbf{Q}\rightarrow\int_{\mathcal{X}}\theta_{\textbf{Q}}^{j,p}(x)Q_{j}^p(dx),\\
\textbf{Q}\rightarrow\int_{\mathcal{X}}\gamma_{\textbf{Q}}^{j,c}(x)Q_{j}^c(dx),&\qquad \textbf{Q}\rightarrow\int_{\mathcal{X}}\gamma_{\textbf{Q}}^{j,p}(x)Q_{j}^p(dx),\\
\end{align*}
 are continuous at any $\textbf{Q}$ where $L(\textbf{Q})<\infty$. First, from Assumption $\ref{bound-ass}$, there exists a positive constant $C>0$ such that, for each $1\leq j\leq r$,

\begin{equation}
\begin{split}
\left|\theta_{\textbf{Q}}^{j,c}(x)\right|&\leq\sup_{\xi,\zeta}\bigg(\bigg|\log\bigg(\sum_{(x_{t-},x_t)\in\mathcal{E}}\lambda_{x_{t-},x_t}^c\left(\xi,\zeta\right)\bigg)\bigg|\bigg)\varphi (x)\\
                 &\leq C(1+\varphi (x)),\quad\forall x\in\mathcal{X},
\label{theta1-bound}
\end{split}
\end{equation}
and
\begin{equation}
\begin{split}
\left|\theta_{\textbf{Q}}^{j,p}(x)\right|&\leq\sup_{\xi,\zeta_1,\ldots,\zeta_r}\bigg(\bigg|\log\bigg(\sum_{(x_{t-},x_t)\in\mathcal{E}}\lambda_{x_{t-},x_t}^p\left(\xi,\zeta_1,\ldots,\zeta_r\right)\bigg)\bigg|\bigg)\varphi (x)\\
                 &\leq C(1+\varphi (x)),\quad\forall x\in\mathcal{X}.
\end{split}
\label{theta2-bound}
\end{equation}

Similarly, by Assumption \ref{bound-ass} we have that, for each $1\leq j\leq r$,
\begin{equation}
\left|\gamma_{\textbf{Q}}^{j,c}(x)\right|\leq\sup_{\xi,\zeta}\bigg|\int_0^T\bigg(\sum_{z:(x_t,z)\in\mathcal{E}}\lambda^c_{x_t,z}(\xi,\zeta)-1\bigg)dt\bigg|<\infty,\quad\forall x\in\mathcal{X},
\label{gamma1-bound}
\end{equation}
 and
\begin{equation}
\left|\gamma_{\textbf{Q}}^{j,p}(x)\right|\leq\sup_{\xi,\zeta_1,\ldots,\zeta_r}\bigg|\int_0^T\bigg(\sum_{z:(x_t,z)\in\mathcal{E}}\lambda^c_{x_t,z}(\xi,\zeta_1,\ldots,\zeta_r)-1\bigg)dt\bigg|<\infty,\quad\forall x\in\mathcal{X}.
\label{gamma2-bound}
\end{equation}

Take $\mathbf{Q'}\in\mathcal{M}_{1,\varphi}(\mathcal{X})\times\cdots\times\mathcal{M}_{1,\varphi}(\mathcal{X})$ in the neighborhood of $\mathbf{Q}$. Note that,
\begin{equation}
\begin{split}
\big|h(\mathbf{Q})-h(\mathbf{Q'})\big|\leq\sum_{j=1}^r\bigg[&\frac{N_j^c}{N}\bigg(\big|\langle \theta_{\textbf{Q}}^{j,c}, Q_{j}^c \rangle-\langle \theta_{\textbf{Q}'}^{j,c}, Q_{j}'^c \rangle\big|+ \big|\langle \gamma_{\textbf{Q}}^{j,c}, Q_{j}^c \rangle-\langle \gamma_{\textbf{Q}'}^{j,c}, Q_{j}'^c \rangle\big|
\bigg)\\
                    &+\frac{N_j^p}{N}\bigg(\big|\langle \theta_{\textbf{Q}}^{j,p}, Q_{j}^p \rangle-\langle \theta_{\textbf{Q}'}^{j,p}, Q_{j}'^p \rangle\big|+ \big|\langle \gamma_{\textbf{Q}}^{j,p}, Q_{j}^p \rangle-\langle \gamma_{\textbf{Q}'}^{j,c}, Q_{j}'^p \rangle\bigg|.
\bigg)
\bigg]
\end{split}
\end{equation}

In addition, for each $1\leq j\leq r$, the following inequalities hold
\begin{align}
\big|\langle \theta_{\textbf{Q}}^{j,c}, Q_{j}^c \rangle-\langle \theta_{\textbf{Q}'}^{j,c}, Q_{j}'^c \rangle\big|&\leq \big|\langle \theta_{\textbf{Q}}^{j,c}, Q_{j}^c-Q_{j}'^c \rangle\big| +\big|\langle \theta_{\textbf{Q}}^{j,c}-\theta_{\textbf{Q}'}^{j,c}, Q_{j}'^c \rangle\big|,
\label{theta-1}
\end{align}
\begin{align}
\big|\langle \theta_{\textbf{Q}}^{j,p}, Q_{j}^p \rangle-\langle \theta_{\textbf{Q}'}^{j,p}, Q_{j}'^p \rangle\big|&\leq \big|\langle \theta_{\textbf{Q}}^{j,p}, Q_{j}^p-Q_{j}'^p \rangle\big| +\big|\langle \theta_{\textbf{Q}}^{j,p}-\theta_{\textbf{Q}'}^{j,p}, Q_{j}'^p \rangle\big|,
\end{align}
\begin{align}
\big|\langle \gamma_{\textbf{Q}}^{j,c}, Q_{j}^c \rangle-\langle \gamma_{\textbf{Q}'}^{j,c}, Q_{j}'^c \rangle\big|&\leq \big|\langle \gamma_{\textbf{Q}}^{j,c}, Q_{j}^c-Q_{j}'^c \rangle\big| +\big|\langle \gamma_{\textbf{Q}}^{j,c}-\gamma_{\textbf{Q}'}^{j,c}, Q_{j}'^c \rangle\big|,
\end{align}
\begin{align}
\big|\langle \gamma_{\textbf{Q}}^{j,p}, Q_{j}^p \rangle-\langle \gamma_{\textbf{Q}'}^{j,p}, Q_{j}'^p \rangle\big|&\leq \big|\langle \gamma_{\textbf{Q}}^{j,p}, Q_{j}^p-Q_{j}'^p \rangle\big| +\big|\langle \gamma_{\textbf{Q}}^{j,p}-\gamma_{\textbf{Q}'}^{j,p}, Q_{j}'^p \rangle\big|.
\end{align}

The idea now is to control the right-hand sides of the last four inequalities. We show this for the inequality in $(\ref{theta-1})$. Similar arguments can be used to treat the three other inequalities. First, notice that the function $ \theta_{\textbf{Q}}^{j,c}$ is continuous. Indeed, the topology of $\mathcal{X}$ is built such that the function $x\rightarrow\sum_{0\leq t\leq T}\mathds{1}_{x_t\neq x_{t-}}$ is continuous. Moreover, from Assumption $\ref{bound-ass}$, the functions $\lambda^c_{z,z'}$ are continuous. Furthermore, from Lemma $\ref{pi-cont}$, the component projection $Q_j^c\rightarrow\pi(Q_j^c)=(Q_{j}^c(t))_{0\leq t\leq T}$ is continuous since $\pi(\mathbf{Q})=(\mathbf{Q}(t))_{0\leq t\leq T}$ is continuous. Finally, the $\log$ function being continuous gives  that $ \theta_{\textbf{Q}}^{j,c}$ is continuous. In addition, from $(\ref{theta1-bound})$ we have that $\theta_{\textbf{Q}}^{j,c}\leq C(1+\varphi (\mathcal{X}))$, thus $\theta_{\textbf{Q}}^{j,c}\in C_{\varphi}(\mathcal{X})$ provided that $L(\mathbf{Q})<\infty$. Therefore, the term $ \big|\langle \theta_{\textbf{Q}}^{j,c}, Q_{j}^c-Q_{j}'^c \rangle\big|$ is as small as desired by taking $\mathbf{Q}'$ close enough to $\mathbf{Q}$ (and thus $Q_{j}'^c$ close enough to $Q_{j}^c$ ). The second term in the right-hand side of $(\ref{theta-1})$ is bounded as follows:
\begin{equation}
\begin{split}
\big|\langle \theta_{\textbf{Q}}^{j,c}-\theta_{\textbf{Q}'}^{j,c}, Q_{j}'^c \rangle\big|\leq \sup_{t}&\bigg|\log\bigg(\sum_{(x_{t-},x_{t})\in\mathcal{E}}\lambda_{x_{t-},x_{t}}^c\left(Q^c_j(t-),Q^p_j(t-)\right)\bigg)\\
                     &-\log\bigg(\sum_{(x_{t-},x_t)\in\mathcal{E}}\lambda_{x_{t-},x_{t}}^c\left(Q'^c_j(t-),Q'^p_j(t-)\right)\bigg)\bigg|\int_{\mathcal{X}}\varphi dQ'^c_j.
\end{split}
\label{theta-2}
\end{equation}

Therefore, using again Assumption \ref{bound-ass}, Lemma \ref{pi-cont} and the continuity of the $\log$ function, the right-hand side of $(\ref{theta-2})$ is controlled for any $\mathbf{Q}'$ in the neighborhood of  $\mathbf{Q}$ in $\mathcal{M}_{1,\varphi}(\mathcal{X})\times\cdots\mathcal{M}_{1,\varphi}(\mathcal{X})$ provided that $L(\mathbf{Q})<\infty$. Thus, the integral $\textbf{Q}\rightarrow\int\theta_{\textbf{Q}}^{j,c} dQ_{j}^c$ is continuous. The exact same steps allow us to show that
\begin{align*}
 \textbf{Q}\rightarrow\int_{\mathcal{X}}\theta_{\textbf{Q}}^{j,p}(x)Q_{j}^p(dx),\qquad \textbf{Q}\rightarrow\int_{\mathcal{X}}\gamma_{\textbf{Q}}^{j,c}(x)Q_{j}^c(dx),\qquad\text{and}\quad \textbf{Q}\rightarrow\int_{\mathcal{X}}\gamma_{\textbf{Q}}^{j,p}(x)Q_{j}^p(dx)
\end{align*}
 are also continuous at any $\textbf{Q}$ where $L(\textbf{Q})<\infty$. Hence, the function $h$ is a linear combination of continuous functions and thus is continuous, which concludes the proof.  \carre
\\
\\
The final step before applying the Laplace-Varadhan principle is to verify that $(\ref{vara-cond-1})$ is satisfied.

\begin{lem}
For any $\alpha>0$,

\begin{align*}
\limsup_{N\rightarrow\infty}\frac{1}{N}\log\int_{\mathcal{M}_{1,\varphi (\mathcal{X})}\times\cdots\times\mathcal{M}_{1,\varphi (\mathcal{X})}}\exp\left\{ N\alpha |h|   \right\} dP_{\nu_N}^{0,N}<\infty.
\end{align*}
\label{vara-bound}
\end{lem}

\proof  First, note that, using the bounds $(\ref{theta1-bound}), (\ref{theta2-bound}), (\ref{gamma1-bound})$ and $(\ref{gamma2-bound})$ we find that, for all $\textbf{Q}\in\mathcal{M}_{1,\varphi}(\mathcal{X})\times\cdots\times\mathcal{M}_{1,\varphi}(\mathcal{X})$,
\begin{equation}
    |h(\mathbf{Q})| \leq\sum_{j=1}^r\bigg[\frac{N_j^c}{N}C\left(1+\int_{\mathcal{X}}\varphi (x) Q_{j}^c(dx)\right)
 + \frac{N_j^p}{N}C\bigg(1+\int_{\mathcal{X}}\varphi (x)Q_{j}^p(dx)\bigg)\bigg].
\label{h-bound}
\end{equation}

 Therefore, in order to show $(\ref{vara-cond-1})$, it is enough to show that, for any $\alpha>0$,
\begin{align}
\limsup_{N\rightarrow\infty}\frac{1}{N}\log\int_{\mathcal{M}_{1,\varphi (\mathcal{X})}\times\cdots\times\mathcal{M}_{1,\varphi (\mathcal{X})}}\exp\bigg\{ N\alpha \sum_{j=1}^r\bigg[\frac{N_j^c}{N}\int_{\mathcal{X}}\varphi (x) Q_{j}^c(dx)
+\frac{N_j^p}{N}\int_{\mathcal{X}}\varphi (x)Q_{j}^p(dx)\bigg]
    \bigg\} dP_{\nu_N}^{0,N}(\mathbf{Q})<\infty.
 \label{h-bound-pro}
\end{align}

Recall that $P_{\nu^N}^{0,N}=\mathbb{P}_{z^N}^{0,N}\circ G_N^{-1}$ where $\mathbb{P}_{z^n}^{0,N}=\otimes_{n=1}^NP_{z_n}$ and $P_{z_n}$ is the  law of the $n$-th particle in the case of non-interaction, with the initial condition being $z_n$. Hence, by independence, the integral term in the left-hand side of $(\ref{h-bound-pro})$ is equivalent to
\begin{equation}
\begin{split}
\prod_{j=1}^r \bigg(\int_{\mathcal{M}_{1,\varphi (\mathcal{X})}}\exp\bigg\{ N_j^c\alpha\int_{\mathcal{X}}\varphi (x) Q_{j}^c(dx)\bigg\}dP_{\nu_N^{j,p}}^{0,N_j^c}(Q_{j}^c)
\int_{\mathcal{M}_{1,\varphi (\mathcal{X})}}\exp\bigg\{N_j^p\alpha\int_{\mathcal{X}}\varphi (x)Q_{j}^p(dx)\bigg\} dP_{\nu_N^{j,p}}^{0,N_j^p}(Q_{j}^p)\bigg).
\end{split}
\label{prod-int}
\end{equation}

 Now, using \cite[Lem.~2.10]{Leonard95}, we find that, for all $1\leq j\leq r$,
\begin{align}
\limsup_{N_j^c\rightarrow\infty}\frac{1}{N_j^c}\log \int_{\mathcal{M}_{1,\varphi (\mathcal{X})}}\exp\bigg\{ N_j^c\alpha \int_{\mathcal{X}}\varphi (x) Q_{j}^c(dx)\bigg\}dP_{\nu_N^{j,c}}^{0,N_j^c}(Q_{j}^c)<\infty,
\label{lim-sup-fin1}
\end{align}
 and
\begin{align}
\limsup_{N_j^p\rightarrow\infty}\frac{1}{N_j^p}\log \int_{\mathcal{M}_{1,\varphi (\mathcal{X})}}\exp\bigg\{ N_j^p\alpha \int_{\mathcal{X}}\varphi (x) Q_{j}^p(dx)\bigg\}dP_{\nu_N^{j,p}}^{0,N_j^p}(Q_{j}^p)<\infty.
\label{lim-sup-fin2}
\end{align}

Since $N_j^c<N$ and $N_j^p<N$ for all $1\leq j\leq r$, $(\ref{lim-sup-fin1})$, $(\ref{lim-sup-fin2})$ and $(\ref{prod-int})$ lead to $(\ref{h-bound-pro})$, which concludes the proof.  \carre

\subparagraph*{The interacting case.} We are now ready to apply the Laplace-Varadhan principle to the sequence of probability measures $\{P^{0,N}_{\nu_N},N\geq 1\}$. By Lemma $\ref{non-inter-LDP}$, the sequence $\{P^{0.N}_{\nu_N},N\geq 1\}$ obeys a large deviation principle in the the topological space $\mathcal{M}_{1,\varphi (\mathcal{X})}\times\cdots\times\mathcal{M}_{1,\varphi (\mathcal{X})}$ with rate function $L(\mathbf{Q})$, defined by $(\ref{rate-func-noninter})$, and speed $N$. By Lemma $\ref{h-continuity}$, the real function $h$ defined in $(\ref{h-func})$ is continuous at any $\mathbf{Q}$ such that $L(\mathbf{Q})<\infty$. Moreover, using  Lemma $\ref{Q-space}$, the function $h$ is continuous on the set $\big\{\mathbf{Q}\in\mathcal{M}_{1,\varphi (\mathcal{X})}\times\cdots\times\mathcal{M}_{1,\varphi (\mathcal{X})}|L(\mathbf{Q})<\infty\big\}$. Finally, we have seen in Lemma $\ref{vara-bound}$ that $(\ref{vara-cond-1})$ is satisfied. Hence, a straightforward application of \cite[Prop.~2.5]{Leonard95} gives
\begin{equation}
\begin{split}
\frac{1}{N}\log\int_{\mathcal{M}_{1,\varphi (\mathcal{X})}\times\cdots\times\mathcal{M}_{1,\varphi (\mathcal{X})}}\exp\left\{ N h   \right\} dP_{\nu_N}^{0,N}\longrightarrow\sup_{\mathbf{Q}\in\mathcal{M}_{1,\varphi (\mathcal{X})}\times\cdots\times\mathcal{M}_{1,\varphi (\mathcal{X})}}\big[h(\mathbf{Q})-L(\mathbf{Q})\big],
\end{split}
\label{conc-1}
\end{equation}
 as $N\rightarrow\infty$, and the sequence
\begin{equation}
\begin{split}
\bigg\{ \frac{\exp(Nh)}{\int_{\mathcal{M}_{1,\varphi (\mathcal{X})}\times\cdots\times\mathcal{M}_{1,\varphi (\mathcal{X})}}\exp(Nh)dP^{0,N}_{\nu_N}}\cdot P^{0,N}_{\nu_N}, N\geq 1    \bigg\}
\end{split}
\label{conc-2}
\end{equation}
 obeys a large deviation principle with speed $N$ and rate function
\begin{align}
\mathbf{Q}\rightarrow L(\mathbf{Q})-h(\mathbf{Q})-\inf_{\mathbf{Q}'\in\mathcal{M}_{1,\varphi (\mathcal{X})}\times\cdots\times\mathcal{M}_{1,\varphi (\mathcal{X})}}[L(\mathbf{Q}')-h(\mathbf{Q}')].
\end{align}

Now, from $(\ref{rad-nik})$ we have
\begin{align}
 \frac{dP_{\nu^N}^{N}}{dP_{\nu^N}^{0,N}}(\mathbf{Q})=\exp\big\{N h(\mathbf{Q})\big\}.
\end{align}
Since $P_{\nu^N}^{N}$ is a probability measure we obtain
\begin{align}
\int_{\mathcal{M}_{1,\varphi (\mathcal{X})}\times\cdots\times\mathcal{M}_{1,\varphi (\mathcal{X})}}\exp(Nh)dP^{0,N}_{\nu_N}=\int_{\mathcal{M}_{1,\varphi (\mathcal{X})}\times\cdots\times\mathcal{M}_{1,\varphi (\mathcal{X})}}dP_{\nu^N}^{N}=1.
\end{align}
Thus, the left side of $(\ref{conc-1})$ is always zero and so
\begin{align}
\sup_{\mathbf{Q}\in\mathcal{M}_{1,\varphi (\mathcal{X})}\times\cdots\times\mathcal{M}_{1,\varphi (\mathcal{X})}}\big[h(\mathbf{Q})-L(\mathbf{Q})\big]=0,
\end{align}
 which gives that
\begin{align}
\inf_{\mathbf{Q}'\in\mathcal{M}_{1,\varphi (\mathcal{X})}\times\cdots\times\mathcal{M}_{1,\varphi (\mathcal{X})}}[L(\mathbf{Q}')-h(\mathbf{Q}')]=0.
\end{align}

We then conclude that the sequences $\{P_{\nu^N}^{N},N\geq 1\}$ obeys a large deviation principle in the topological space $\mathcal{M}_{1,\varphi (\mathcal{X})}\times\cdots\times\mathcal{M}_{1,\varphi (\mathcal{X})}$, with speed $N$ and rate function
\begin{align}
I(\mathbf{Q})=L(\mathbf{Q})-h(\mathbf{Q}).
\end{align}

In order to obtain the representation $(\ref{empi-meas-rate})$, we proceed as follows: First, from $(\ref{h-bound})$ we have that, for $\mathbf{Q}\in\mathcal{M}_{1,\varphi}(\mathcal{X})\times\cdots\mathcal{M}_{1,\varphi}(\mathcal{X})$, $h(\mathbf{Q})<\infty$. Moreover, from \cite[Lem.~5.6]{Bork+Sund2012}, the functions $J^{j,\iota}(Q)$ defined by $(\ref{Daw-Gar-rate})$ takes the following representation, 
\begin{align}
J^{j,\iota}(Q)=\left\{\begin{array}{ll}
       H(Q|P_{j,\iota}), & \mbox{if $Q\circ\pi_0^{-1}=\nu^{j,\iota}$},\\
       +\infty, & \mbox{Otherwise},
\end{array}\right.
\end{align}
 where $H(Q|P_{j,\iota})$ is the relative entropy defined by $(\ref{relat-entro})$. Therefore, if either $Q\circ\pi_0^{-1}\neq\nu^{j,\iota}$ or $Q$ is not absolutely continuous with respect to $P_{j,\iota}$, one can immediately observe that $J^{j,\iota}(Q)=\infty$, thus $L(\mathbf{Q})=\infty$ and finally $I(\mathbf{Q})=\infty$.

Now, assume that, for all $1\leq j\leq r$, $Q_{j}^c\circ\pi_0^{-1}=\nu^{j,c}$, $Q_{j}^c\ll P$ and $Q_{j}^p\circ\pi_0^{-1}=\nu^{j,p}$, $Q_{j}^p\ll P$, then
\begin{align*}
L(\mathbf{Q})=\alpha_1p_1^c H(Q_{1}^c|P_{1,c})+ \alpha_1p_1^p H(Q_{1}^p|P_{1,p})+\cdots+\alpha_rp_r^c H(Q_{r}^c|P_{r,c})+ \alpha_rp_r^p H(Q_{r}^p|P_{r,p}).
\end{align*}
 Furthermore, one can observe from $(\ref{dens-Rfunc})$ that the densities $\exp\{h_1(x,\eta,\rho_j)\}$ and $\exp\{h_2(x,\eta,\rho_1,\ldots,\rho_r)\}$ do not depend on the initial condition $z_0$. Therefore, for each $1\leq j\leq r$, the densities of $R^c(\pi(Q_{j}^c),\pi(Q_{j}^p))$ and $R^p(\pi(Q_{j}^c)),\pi(Q_{1}^p),\ldots,\pi(Q_{r}^p)))$ with respect to the mixtures $P_{j,c}$ and $P_{j,p}$  are given by, respectively,
\begin{align*}
\frac{d R^c(\pi(Q_{j}^c),\pi(Q_{j}^p))}{d P_{j,c}}(x)=\exp\{h_1(x,\pi(Q_{j}^c),\pi(Q_{j}^p))\},
\end{align*}
 and
\begin{align*}
 \frac{d R^p(\pi(Q_{j}^c)),\pi(Q_{1}^p),\ldots,\pi(Q_{r}^p)))}{d P_{j,p}}(x)=\exp\{h_2(x,\pi(Q_{j}^c)),\pi(Q_{1}^p),\ldots,\pi(Q_{r}^p)))\}.
\end{align*}

Replacing $h_1()$ and $h_2()$ in $(\ref{h-func})$ by the two last representations we find,
\begin{equation}
\begin{split}
h(\mathbf{Q})&=\sum_{j=1}^r\bigg[\frac{N_j^c}{N}\int_{D([0,T],\mathcal{Z})}dQ_{j}^c\log\frac{d R^c(\pi(Q_{j}^c),\pi(Q_{j}^p))}{d P_{j,c}}\\
&\qquad\qquad+\frac{N_j^p}{N}\int_{D([0,T],\mathcal{Z})}dQ_{j}^p\log\frac{d R^p(\pi(Q_{j}^c),\pi(Q_{1}^p),\ldots,\pi(Q_{r}^p)))}{d P_{j,p}}\bigg].
\end{split}
\end{equation}

Finally, using the assumption $(\ref{converg-propo})$ we find that, as $N\rightarrow\infty$,
\begin{equation}
\begin{split}
L(\mathbf{Q})-h(\mathbf{Q})&=\sum_{j=1}^r\bigg[\alpha_jp_j^c\bigg(\int_{D([0,T],\mathcal{Z})}dQ_{j}^c\log\frac{d Q_{j}^c}{d P_{j,c}}-\int_{D([0,T],\mathcal{Z})}dQ_{j}^c\log\frac{d R^c(\pi(Q_{j}^c),\pi(Q_{j}^p))}{d P_{j,c}}\bigg)\\
&\qquad\qquad+\alpha_jp_j^p\bigg(\int_{D([0,T],\mathcal{Z})}dQ_{j}^p\log\frac{d Q_{j}^p}{d P}-\int_{D([0,T],\mathcal{Z})}dQ_{j}^p\log\frac{d R^p(\pi(Q_{j}^c),\pi(Q_{1}^p),\ldots,\pi(Q_{r}^p)))}{d P}\bigg)\bigg]\\
                          &=\sum_{j=1}^r\bigg[\alpha_jp_j^c\int_{D([0,T],\mathcal{Z})}dQ_{j}^c\log\frac{d Q_{j}^c}{d R^c(\pi(Q_{j}^c),\pi(Q_{j}^p))}\\
&\qquad\qquad+\alpha_jp_j^p\int_{D([0,T],\mathcal{Z})}dQ_{j}^p\log\frac{d Q_{j}^p}{d R^p(\pi(Q_{j}^c),\pi(Q_{1}^p),\ldots,\pi(Q_{r}^p))}\bigg]  \\
                          &=\sum_{j=1}^r\bigg[\alpha_jp_j^cH\bigg(Q_{j}^c\big|R^c(\pi(Q_{j}^c),\pi(Q_{j}^p))\bigg) +\alpha_jp_j^pH\bigg(Q_{j}^p\big| R^p(\pi(Q_{j}^c),\pi(Q_{1}^p),\ldots,\pi(Q_{r}^p)\bigg)\bigg].
\end{split}
\end{equation}
This concludes the proof.
   \carre

\subsection{Large deviation principle for the Empirical process}

We investigate in this section the large deviations of the sequence $(p_{\nu_N}^N,N\geq 1)$ where $p_{\nu_N}^N=\mathbb{P}_{z^N}^N\circ\gamma_N^{-1}=\pi(M^N)$ is the distribution of the $\mathcal{M}_1(\mathcal{Z})\times\cdots\times\mathcal{M}_1(\mathcal{Z})$-valued empirical process defined as
\begin{equation*}
\begin{split}
\mu^N: t\in [0,T]\longrightarrow \mu^N(t)&=\left(\mu_1^{c,N}(t),\mu_1^{p,N}(t),\cdots,\mu_r^{c,N}(t),\mu_r^{p,N}(t)\right)\\
                                     &=\bigg(\frac{1}{N_1^c}\sum_{n\in C^c_1}\delta_{X_n(t)},\frac{1}{N_1^p}\sum_{n\in C^p_1}\delta_{X_n(t)},\ldots, \frac{1}{N^c_r}\sum_{n\in C^c_r}\delta_{X_n(t)}, \frac{1}{N^p_r}\sum_{n\in C^p_r}\delta_{X_n(t)}  \bigg).
\end{split}
\end{equation*}

The flow $\mu^N$ takes values in the product space $\big(\mathcal{D}([0,T],\mathcal{M}_1(\mathcal{Z}))\big)^{2r}$. Let again $\mathcal{D}([0,T],\mathcal{M}_1(\mathcal{Z}))$ be equipped with the metric
\begin{align}
\rho_T (\mu,\nu)=\sup_{0\leq t\leq T}\rho_0 (\mu_t,\nu_t),\quad\mu,\nu\in\mathcal{D}([0,T],\mathcal{M}_1(\mathcal{Z})),
\label{rho-metric}
\end{align}
 where $\rho_0(\alpha,\beta)$, $\alpha,\beta\in\mathcal{M}_1(\mathcal{Z})$ is a metric on $\mathcal{M}_1(\mathcal{Z})$ which generates the weak topology on $\mathcal{M}_1(\mathcal{Z})$. Moreover, let the product space $\big(\mathcal{D}([0,T],\mathcal{M}_1(\mathcal{Z}))\big)^{2r}$ be equipped with the product topology induced by the product metric $N(\rho_T,\cdots,\rho_T)$.

For any $\mathbf{\xi}=(\xi_1^c,\xi_1^p,\ldots,\xi_r^c,\xi_r^p)\in\mathcal{M}_1(\mathcal{Z})\times\cdots\mathcal{M}_1(\mathcal{Z})$, let us introduce the rate matrices
\begin{align*}
A_{\xi_j^c}=\left(\lambda^c_{z,z'}(\xi_j^c,\xi_j^p)\right)_{(z,z')\in\mathcal{Z}\times\mathcal{Z}}\qquad\text{and}\quad A_{\xi_j^p}=\left(\lambda^p_{z,z'}(\xi_j^c,\xi_1^p\ldots,\xi^p_r)\right)_{(z,z')\in\mathcal{Z}\times\mathcal{Z}},
\end{align*}
 with $\lambda^c_{z,z}(\xi_j^c,\xi_j^p)=-\sum_{z'\neq z}\lambda^c_{z,z'}(\xi_j^c,\xi_j^p)$ and $\lambda^p_{z,z}(\xi_j^c,\xi_1^p,\ldots,\xi_r^p)=-\sum_{z'\neq z}\lambda^p_{z,z'}(\xi_j^c,\xi_1^p,\ldots,\xi_r^p)$. From the law of large numbers given in Corollary \ref{coro-1}, one can deduces that, as $N\rightarrow\infty$, the sequence $(\mu^N,N\geq 1)$ converges weakly, for converging initial conditions, towards the solution $\mu$ of the following McKean-Vlasov system 
\begin{equation}
\begin{split}
\left\{ \begin{array}{lcl}\dot{\mu}_j^c(t)=A_{\mu_j^c(t)}^{*}\mu_j^c(t), & &  \\
 \dot{\mu}_j^p(t)=A_{\mu_j^p(t)}^{*}\mu_j^p(t),  & & \\
 \mu_j^c(0)=\nu_j^c,\mu_j^p(0)=\nu_j^p, & & \\
 1\leq j\leq r, t\in [0,T], & &
 \end{array}\right.
\end{split}
\label{McKean-Vlas-syst}
\end{equation}
 where $A^*$ is the adjunct/transpose of the matrix $A$ and $\dot{\mu}(t)= \frac{\partial}{\partial t} \mu(t)$. Note that the Lipschitz property of the functions $\lambda^c_{z,z'}$ and $\lambda^p_{z,z'}$ assures that $(\ref{McKean-Vlas-syst})$ is well-posed. Also, one can notice that the representation $(\ref{McKean-Vlas-syst})$ is consistent with the the infinitesimal generators $\mathcal{L}^c_{\xi,\eta_j}$ and $\mathcal{L}_{\xi,\eta_1,\dots,\eta_r}^p$ introduced in $(\ref{gener-c})$ and $(\ref{gener-p})$.  Indeed, by considering $\phi$, $\mathcal{L}^c_{\xi,\eta_j}\phi$  and $\mathcal{L}_{\xi,\eta_1,\dots,\eta_r}^p\phi$ as column vectors, the right-hand sides of $(\ref{gener-c})$ and $(\ref{gener-p})$ are results of right multiplying, respectively, the rates matrices $A_{j,c}$  and $A_{j,p}$ by the vector $\phi$.

Denote by $\tau$  the log-Laplace transform of the centered Poisson distribution with parameter $1$ given by $\tau(u)=e^u-u-1$, and let $\tau^*$ be its Legendre transform defined by
\begin{align*}
\tau^*(u)=\left\{\begin{tabular}{lll}
$(u+1) \log (u+1)-u$& \text{if} & $u>-1$, \\
$1$               & \text{if} & $u=-1$,\\
$+\infty$         & \text{if}  & $u<-1$. \\
\end{tabular}
\right.
\end{align*}

We recall now  the notion of absolute continuity introduced in \cite[Def.~4.1]{Daw+Gart87}. Denote by $\mathcal{S}$ the Schwartz space of test functions $\mathbb{R}^d\rightarrow\mathbb{R}$ having a compact support and possessing continuous derivatives of all orders. We endow $\mathcal{S}$ with the usual inductive topology. Let $\mathcal{S}'$ be the corresponding space of real distributions. For each compact set $K\subset\mathbb{R}^d$, $\mathcal{S}_K$ will denote the subspace of $\mathcal{S}$ consisting of all test
functions, the support of which is contained in $K$. Finally, let $\langle\nu,f\rangle$ denote the application of the test function $f$ to the distribution $\nu$.

\begin{defi}
Let $I$ be an interval of the real line. A map $\nu (\cdot): I\rightarrow\mathcal{S}'$ is called absolutely continuous if for each compact set $K\subset\mathbb{R}^d$, there exist a neighborhood $U_K$ of $0$ in $\mathcal{S}_K$ and an absolutely
continuous function $H_K: I\rightarrow{R}$ such that
\begin{align*}
|\langle\nu(u),f\rangle|-|\langle\nu(v),f\rangle|\leq H_K(u)-H_K(v),
\end{align*}
for all $u,v\in I$ and $f\in U_K$.
\label{abs-cont-def}
\end{defi}

Finally, define,  for any $\theta\in\mathcal{M}(\mathcal{Z})$,
\begin{align*}
|||\theta|||^{j,c}_{\mu(t)}&=\sup_{\Phi:\mathcal{Z}\rightarrow\mathbb{R}} \bigg\{ \sum_{z\in\mathcal{Z}}\theta(z)\cdot\Phi(z)-\sum_{z':(z,z')\in\mathcal{E}} \tau \big(\Phi(z')- \Phi(z)\big)\cdot\mu_j^c(t)(z)\cdot\lambda^c_{zz'}\big(\mu_j^c(t),\mu_j^p(t)\big) \bigg\},\\
|||\theta|||^{j,p}_{\mu(t)}&=\sup_{\Phi:\mathcal{Z}\rightarrow\mathbb{R}} \bigg\{ \sum_{z\in\mathcal{Z}}\theta(z)\cdot\Phi(z)-\sum_{z':(z,z')\in\mathcal{E}} \tau \big(\Phi(z')- \Phi(z)\big)\cdot\mu_j^p(t)(z)\cdot\lambda^p_{zz'}\big(\mu_j^c(t),\mu_1^p(t),\ldots,\mu_r^p(t)\big) \bigg\},
\end{align*}
 and let introduce, for each $\nu\in(\mathcal{M}_1(\mathcal{Z}))^{2r}$, and according to \cite[eqn.~(4.9)]{Daw+Gart87}, the functional $S(\mu|\nu)$ defined from  $(\mathcal{D}([0, T ], \mathcal{M}_1 (Z)))^{2r}$ to $[0,\infty]$ by setting
\begin{equation}
\begin{split}
S_{[0,T]}(\mu|\nu)= \sum_{j=1}^r\bigg[\alpha_jp_j^c \int_0^T ||| \dot{\mu}_j^c(t)-A_{\mu_j^c(t)}^{*}\mu_j^c(t)|||_{\mu(t)}dt  +\alpha_jp_j^p\int_0^T ||| \dot{\mu}_j^p(t)-A_{\mu_j^p(t)}^{*}\mu_j^p(t)|||_{\mu(t)}dt \bigg]
\label{rate-emp-proc}
\end{split}
\end{equation}
 if $\mu(0)=\nu$ and $\mu_j^c,\mu_j^p$ are absolutely continuous in the sense of Definition $\ref{abs-cont-def}$ for all $1\leq j\leq r$, and $S_{\nu}(\mu)=+\infty$ otherwise.

We are now ready to formulate our large deviations result.

\begin{theo}
Suppose that $\nu_N\rightarrow\nu$ weakly. The sequence of probability measures $(p^N_{\nu_N},N\geq 1)$ obeys a large deviation principle in the space $\big(\mathcal{D}([0,T],\mathcal{M}_1(\mathcal{Z}))\big)^{2r}$, with speed $N$, and rate function $S_{[0,T]}(\mu|\nu)$ given by  $(\ref{rate-emp-proc})$. \\
\indent Moreover, if a path $\mu\in\big(\mathcal{D}([0,T],\mathcal{M}_1(\mathcal{Z}))\big)^{2r}$ satisfies $S_{[0,T]}(\mu|\nu)<\infty$, then  $\mu_j^c$ and $\mu_j^p$  are absolutely continuous and there exist rate families $(l_{z,z'}^{j,c}(t), t \in [0, T ], (z, z')\in\mathcal{E})$ and $(l_{z,z'}^{j,p}(t), t \in [0, T ], (z, z')\in\mathcal{E})$ such that, for all $1\leq j\leq r$,
\begin{align*}
\dot{\mu}_j^c(t)&={L_{j,c}(t)}^{*}\mu_j^c(t),\\
\dot{\mu}_j^p(t)&={L_{j,p}(t)}^{*}\mu_j^p(t),
\end{align*}
where $L_{j,c}(t)$ (resp. $L_{j,p}(t)$) is the rate matrix associated with the time-varying rates $(l_{z,z'}^{j,c}(t), (z, z')\in\mathcal{E})$ (resp. $(l_{z,z'}^{j,p}(t), (z, z')\in\mathcal{E})$) and $L_{j,c}(t)^*$ (resp. $L_{j,p}(t)^*$) is its adjoint. Furthermore, in this case, the good rate function $S_{[0,T]}(\mu|\nu)$ is given by
\begin{equation}
\begin{split}
\sum_{j=1}^r\bigg[&\alpha_jp_j^c \int_0^T \bigg( \sum_{(z,z')\in\mathcal{E}}(\mu_j^c(t)(z))\lambda^c_{z,z'}\left(\mu_j^c(t),\mu_j^p(t)\right) \tau^*\bigg(\frac{l_{z,z'}^{j,c}(t)}{\lambda^c_{z,z'}\left(\mu_j^c(t),\mu_j^p(t)\right)}-1\bigg) \bigg)dt \\
 &+\alpha_jp_j^p\int_0^T  \bigg( \sum_{(z,z')\in\mathcal{E}}(\mu_j^p(t)(z))\lambda^p_{z,z'}\left(\mu_j^c(t),\mu_1^p(t),\ldots,\mu_r^p(t)\right) \tau^*\bigg(\frac{l_{z,z'}^{j,p}(t)}{\lambda^p_{z,z'}\left(\mu_j^c(t),\mu_1^p(t),\ldots,\mu_r^p(t)\right)}-1\bigg) \bigg) dt \bigg].
\label{rate-emp-proc-2}
\end{split}
\end{equation}
\label{large-dev-emp-proc}
\end{theo}

\proof We first use a contraction argument to derive a large deviation result for the sequence $(p_{\nu_N}^N,N\geq 1)$. From Theorem \ref{large-dev-meas}, the sequence $(P_{\nu_N}^N,N\geq 1)$ obeys a large deviation principle with speed $N$ and rate function $I(\mathbf{Q})$ given by 
\begin{align*}
I(\mathbf{Q})= \left\{ \begin{array}{lll}
\sum_{j=1}^r\bigg[\alpha_jp_j^cH\bigg(Q_{j}^c\big|R^c(\pi(Q_{j}^c),\pi(Q_{j}^p))\bigg) +\alpha_jp_j^pH\bigg(Q_{j}^p\big| R^p(\pi(Q_{j}^c),\pi(Q_{1}^p),\ldots,\pi(Q_{r}^p)\bigg)\bigg] & \text{if $\mathbf{Q}\circ\pi_0^{-1}=\nu$},& \\
+\infty &\text{Otherwise}.&
\end{array}
\right.
\end{align*}

Moreover, from Lemma \ref{pi-cont} the projection:
\begin{align*}
\pi: \big(\mathcal{M}_1(\mathcal{D}([0,T],\mathcal{Z}))\big)^{2r}&\rightarrow\big(\mathcal{D}([0,T],\mathcal{M}_1(\mathcal{Z}))\big)^{2r}
\end{align*}
is continuous at each $\textbf{Q}\in\big(\mathcal{M}_1(\mathcal{D}([0,T],\mathcal{Z}))\big)^{2r}$ where $L(\textbf{Q})<\infty$ and thus at any $\mathbf{Q}$ such that $I(\mathbf{Q})<\infty$. The latter corresponds to the effective domain $\mathcal{D}_I=\{\mathbf{Q}:I(\mathbf{Q})<\infty\}$  of the rate function $I$ (see \cite[p. 4]{Dem+Zeit2010}). Therefore, by applying the contraction principle to the large deviation principle of  $(P_{\nu_N}^N,N\geq 1)$ (see \cite[Th.~4.2.1: Remark (c)]{Dem+Zeit2010} ) with rate $I$, we deduce that the family of probability measures $(P_{\nu_N}^N\circ\pi^{-1},N\geq 1)$ obeys a large deviation principle in $\big(\mathcal{D}([0,T],\mathcal{M}_1(\mathcal{Z}))\big)^{2r}$ with rate function defined, for any $\mu\in\big(\mathcal{D}([0,T],\mathcal{M}_1(\mathcal{Z}))\big)^{2r}$, by
\begin{align}
V(\mu)=\inf\bigg\{ I(\mathbf{Q}), \mathbf{Q}\in\big(\mathcal{M}_1(D([0,T],\mathcal{Z}))\big)^{2r},\pi(\mathbf{Q})=\mu  \bigg\}.
\label{V-rate-func}
\end{align}

We now derive another representation for the rate function $V$ following \cite{Daw+Gart87} and \cite{Leonard95}. Fix $\mu=\big(\mu_j^c,\mu_j^p,1\leq j\leq r\big)\in\big(\mathcal{D}([0,T],\mathcal{M}_1(\mathcal{Z}))\big)^{2r}$. Note that writing $\pi(\mathbf{Q})=\mu$, with $\mathbf{Q}\in\big(\mathcal{M}_1(D([0,T],\mathcal{Z}))\big)^{2r}$, is amount to $\pi(Q_{j}^c)=\mu_j^c$ and $\pi(Q_{j}^p)=\mu_j^p$ for all $1\leq j\leq r$. Therefore $V$ can be rewritten as
\begin{align*}
V(\mu)=\inf\bigg\{ \sum_{j=1}^r\bigg[\alpha_jp_j^cH\bigg(Q_{j}^c\big|R^c(\mu_j^c,\mu_j^p)\bigg) &+\alpha_jp_j^pH\bigg(Q_{j}^p\big| R^p(\mu_j^c,\mu_1^p,\ldots, \mu_r^p) \bigg)\bigg],\\
&\qquad\qquad \mathbf{Q}\in\big(\mathcal{M}_1(D([0,T],\mathcal{Z}))\big)^{2r}, \pi(\mathbf{Q})=\mu  \bigg\}.
\end{align*}

Fix $1\leq j\leq r$. Let $\left(X^{(i)}_{j,c}\right)_{i\geq 1}$ and $\left(X^{(i)}_{j,p}\right)_{i\geq 1}$ be sequences of i.i.d. processes with common law $R^c(\mu_j^c,\mu_j^p)$ and $R^p(\mu_j^c,\mu_1^p,\ldots, \mu_r^p)$, respectively. Therefore, by Sanov's theorem, the empirical measures $\frac{1}{N_j^c}\sum_{i=1}^{N_j^c}X^{(i)}_{j,c}$ and $\frac{1}{N_j^p}\sum_{i=1}^{N_j^p}X^{(i)}_{j,p}$  obey  large deviation principles as $N_j^c\rightarrow\infty$ and $N_j^p\rightarrow\infty$ with speed $N_j^c$ and $N_j^p$, respectively, and rate functions given by, respectively,
\[
    Q\in \mathcal{M}_1(\mathcal{D}([0,T],\mathcal{Z}))\rightarrow H(Q|R^c(\mu_j^c,\mu_j^p)),
\]
and
\[
    Q\in \mathcal{M}_1(\mathcal{D}([0,T],\mathcal{Z}))\rightarrow H(Q|R^p(\mu_j^c,\mu_1^p,\ldots,\mu_r^p)).
\]

Using the same arguments as in the proof of Lemma \ref{pi-cont}, one can show that the projection $\pi$ is continuous at any  $\mathbf{Q}\in\big(\mathcal{M}_1(D([0,T],\mathcal{Z}))\big)^{2r}$ such that
\begin{align*}
\sum_{j=1}^r\left[\alpha_jp_j^cH\left(Q_{j}^c\big|R^c(\mu_j^c,\mu_j^p)\right) +\alpha_jp_j^pH\left(Q_{j}^p\big| R^p(\mu_j^c,\mu_1^p,\ldots, \mu_r^p \right)\right]<\infty.
\end{align*}
Thus, the component projections $\pi(Q_{j}^c)$ and $\pi(Q_{j}^p)$ are also continuous. Hence, using the contraction principle (\cite[Th.~4.2.1]{Dem+Zeit2010}), the sequences $\left(t\in [0,T]\rightarrow \frac{1}{N_j^c}\sum_{i=1}^{N_j^c}X^{(i)}_{j,c}(t)\right)$ and $\left(t\in [0,T]\rightarrow \frac{1}{N_j^p}\sum_{i=1}^{N_j^p}X^{(i)}_{j,p}(t)\right)$ obey large deviation principles with rate functions, respectively,
\begin{align*}
&\qquad\eta\in\mathcal{D}([0,T],\mathcal{M}_1(\mathcal{Z})) \rightarrow S^{j,c}_{\mu}(\eta)=\inf\left\{ H\left(Q\big|R^c(\mu_j^c,\mu_j^p)\right), Q\in\mathcal{M}_1(D([0,T],\mathcal{Z})),\pi(Q)=\eta  \right\},\\
&\text{and}\\
&\qquad\eta\in\mathcal{D}([0,T],\mathcal{M}_1(\mathcal{Z}))\rightarrow S^{j,p}_{\mu}(\eta)=\inf\left\{ H\left(Q\big|R^p(\mu_j^c,\mu_1^p,\ldots,\mu_r^p)\right), Q\in\mathcal{M}_1(D([0,T],\mathcal{Z})),\pi(Q)=\eta  \right\}.
\end{align*}

Note that, by using independence argument and following the same steps as in the proof of Lemma \ref{non-inter-LDP}, one can show that the sequence
\begin{align*}
\bigg(t\in [0,T]\rightarrow \bigg(\frac{1}{N_j^c}\sum_{i=1}^{N_j^c}X^{(i)}_{j,c}(t),\frac{1}{N_j^c}\sum_{i=1}^{N_j^p}X^{(i)}_{j,p}(t),\ldots, \frac{1}{N_j^c}\sum_{i=1}^{N_r^c}X^{(i)}_{r,c}(t),\frac{1}{N_r^c}\sum_{i=1}^{N_r^p}X^{(i)}_{r,p}(t) \bigg)\bigg)_{N\geq 1}
\end{align*}
 obeys a large deviation principle with speed $N$ and rate function
\begin{align*}
\mathbf{\eta}=(\eta_1^c,\eta_1^p,\ldots,\eta_r^c,\eta_r^p)\in\big(\mathcal{D}([0,T],\mathcal{M}_1(\mathcal{Z}))\big)^{2r}\rightarrow S_{\mu}(\eta)=\sum_{j=1}^r\bigg[\alpha_jp_j^c S^{j,c}_{\mu}(\eta_j^c) +\alpha_jp_j^pS^{j,p}_{\mu}(\eta_j^p))\bigg].
\end{align*}

In addition, the vector
\begin{align*}
\left( \frac{1}{N_j^c}\sum_{i=1}^{N_j^c}X^{(i)}_{j,c},\frac{1}{N_j^c}\sum_{i=1}^{N_j^p}X^{(i)}_{j,p},\ldots, \frac{1}{N_j^c}\sum_{i=1}^{N_r^c}X^{(i)}_{r,c},\frac{1}{N_r^c}\sum_{i=1}^{N_r^p}X^{(i)}_{r,p} \right)_{N\geq 1}
\end{align*}
obeys a large deviation principle with rate $I(Q)$. Therefore, by contraction argument and using again the continuity of the projection we find that
\begin{align*}
\bigg(t\in [0,T]\rightarrow\bigg( \frac{1}{N_j^c}\sum_{i=1}^{N_j^c}X^{(i)}_{j,c}(t),\frac{1}{N_j^c}\sum_{i=1}^{N_j^p}X^{(i)}_{j,p}(t),\ldots, \frac{1}{N_j^c}\sum_{i=1}^{N_r^c}X^{(i)}_{r,c}(t),\frac{1}{N_r^c}\sum_{i=1}^{N_r^p}X^{(i)}_{r,p}(t) \bigg)\bigg)_{N\geq 1}
\end{align*}
obeys a large deviation principle with rate $V(\mu)$. Hence, by the uniqueness of the rate function (cf. \cite[Lem.~2.1.1]{Deus+Strook89}), we find
\begin{align}
V(\mu)=S_{\mu}(\mu).
\label{V-S-equal}
\end{align}

We next derive another representation for $S_{\mu}(\nu)$. For any $1\leq j\leq r$ and $\nu\in\mathcal{D}([0,T],\mathcal{M}_1(\mathcal{Z}))$ we have from \cite[p.~319]{Leonard95} that
\begin{equation}
\begin{split}
S^{j,c}_{\mu}(\nu)&=U^{j,c}_{\mu}(\nu_t(dz)dt), \\
S^{j,p}_{\mu}(\nu)&=U^{j,p}_{\mu}(\nu_t(dz)dt),
\label{S-U-equal}
\end{split}
\end{equation}
 where, for all $\tilde{\nu}\in\mathcal{M}_1([0,T[\times\mathcal{Z})$, $U^{j,c}_{\mu}(\tilde{\nu})$ and $U^{j,p}_{\mu}(\tilde{\nu})$ are given by (see \cite[eqn.~(3.14)]{Leonard95})
\begin{equation}
\begin{split}
U^{j,c}_{\mu}(\tilde{\nu})=\sup_{f\in C_1^c}\bigg\{\int_0^T\bigg\langle-&\bigg(\frac{\partial}{\partial t}+A_{\mu_j^c(t)}\bigg)f(t,z)\\
&-\sum_{z':(z,z')\in\mathcal{E}}\tau\big(f(t,z')-f(t,z)\big)\lambda^c_{zz'}\big(\mu_j^c(t),\mu_j^p(t)\big),\nu_t (dz)\bigg\rangle\bigg\},
\label{U-c-func}
\end{split}
\end{equation}

\begin{equation}
\begin{split}
U^{j,p}_{\mu}(\tilde{\nu})=\sup_{f\in C_1^c}\bigg\{\int_0^T\bigg\langle-&\bigg(\frac{\partial}{\partial t}+A_{\mu_j^p(t)}\bigg)f(t,z)\\
&-\sum_{z':(z,z')\in\mathcal{E}}\tau\big((f(t,z')-f(t,z)\big)\lambda^p_{zz'}\big(\mu_j^c(t),\mu_1^p(t),\ldots,\mu_r^p(t)\big),\nu_t (dz)\bigg\rangle\bigg\},
\label{U-p-func}
\end{split}
\end{equation}
where $C_1^c$ stands for the set of all continuous functions with compact support on $[0,T[ \times\mathcal{Z}$ which are
 $t$-differentiable. Using $(\ref{S-U-equal})$, $(\ref{U-c-func})$ and $(\ref{U-p-func})$ together with  \cite[Lem.~3.2]{Leonard95} we obtain
\begin{equation}
\begin{split}
S^{j,c}_{\mu}(\mu)&=\int_0^T ||| \dot{\mu}_j^c(t)-A_{\mu_j^c(t)}^{*}\mu_j^c(t)|||_{\mu(t)}dt, \\
S^{j,p}_{\mu}(\mu)&=\int_0^T ||| \dot{\mu}_j^p(t)-A_{\mu_j^p(t)}^{*}\mu_j^p(t)|||_{\mu(t)}dt.
\end{split}
\end{equation}

Finally, using $(\ref{V-S-equal})$, we deduce that $(p^N_{\nu_N},N\geq 1)$ obeys a large deviation principle with rate $N$ and good rate function $(\ref{rate-emp-proc})$. The representation $(\ref{rate-emp-proc-2})$ follows immediately from \cite[Lem.~3.2]{Leonard95} and the statement about absolute continuity follows from \cite[Th.~3.1]{Leonard95}. The theorem is proved. \carre
\\
\\
The following result shows that the large deviation principle for $(p_{\nu}^{N},N\geq 1)$ holds uniformly in the initial condition.
\begin{coro}
For any compact set $K\subset\big(\mathcal{M}_1(\mathcal{Z}))\big)^{2r}$, any closed set $F\subset\big(\mathcal{D}([0,T],\mathcal{M}_1(\mathcal{Z}))\big)^{2r}$, and any open set $G\subset\big(\mathcal{D}([0,T],\mathcal{M}_1(\mathcal{Z}))\big)^{2r}$, we have
\begin{align}
\limsup_{N\rightarrow\infty}\frac{1}{N} \log \sup_{\nu\in K} p_{\nu}^{N} \big({\mu_N\in F}\big) \leq - \inf_{\substack{\nu\in K\\ \mu\in F}} S_{[0,T]} (\mu|\nu),
\end{align}

\begin{align}
\limsup_{N\rightarrow\infty}\frac{1}{N} \log \sup_{\nu\in K} p_{\nu}^{N} \big({\mu_N\in G}\big) \geq - \sup_{\nu\in K}\inf_{\substack{ \mu\in G}} S_{[0,T]} (\mu|\nu),
\end{align}
\label{coro-ldp-emp-proc}
\end{coro}

  \proof This follows immediately from \cite[Coro.~5.6.15]{Dem+Zeit2010} and Theorem \ref{large-dev-emp-proc}. \carre

\section*{Acknowledgment}
This research was supported by the Natural Sciences and Engineering Research Council of Canada Discovery Grants and by Carleton University.

\bibliographystyle{livre} 
\bibliography{biblio}
\end{document}